\documentclass[a4paper,reqno]{amsart}
\usepackage[left=1in,right=1in,bottom=1.2in,top=1.2in]{geometry}
\usepackage{tikz,amsthm,amsmath,amstext,amssymb,amscd,epsfig,euscript, mathrsfs, dsfont,pspicture,multicol, mathtools,graphpap,graphics,graphicx,times,enumerate,subfig,sidecap,wrapfig,color}
\usepackage{hyperref}
\DeclareMathAlphabet{\mathbbb}{U}{bbold}{m}{n}
\hypersetup{
colorlinks,
citecolor=green,
filecolor=black,
linkcolor=blue,
urlcolor=blue,
}
\usepackage{bm}
\usepackage{amsmath}
\usepackage{mathrsfs}  
\usepackage{amsthm}
\usepackage{amssymb}

\def \rr {{\mathbb R}}
\def \cc {{\mathbb C}}
\def \nn {{\mathbb N}}

\def \zz {{\mathbb Z}}

\DeclareGraphicsExtensions{.png}
 \numberwithin{equation}{section}

\def\var{\varphi}

\def\g{\gamma}

\def\G{\Gamma}

\def\O{\Omega}

\def\S{\Sigma}

\def\sS{{\mathcal{S}}} 
\def\cP{{\mathcal{P}}}

\def \ty {\tilde{y}}

\def \cc {{\mathbb{C}}}
\newtheorem{definition}{Definition}[section]
\newtheorem{theorem}{Theorem}[section]
\newtheorem{proposition}{Proposition}[section]
\newtheorem{lemma}{Lemma}[section]
\newtheorem{corollary}{Corollary}[section]
\newtheorem{remark}{Remark}[section]

	\title[On the approximation of the Dirac operator coupled with $\delta$-shell interactions]{On the approximation of the Dirac operator coupled with confining Lorentz scalar $\delta$-shell interactions}
\author[Mahdi ZREIK]{Mahdi ZREIK\textsuperscript{1}}
\address{\textsuperscript{1}Institut de Math\'ematiques de Bordeaux, UMR 5251, Universit\' e de Bordeaux 33405 Talence Cedex, FRANCE\\ and Departamento de Matem\'aticas, Universidad del Pa\' is Vasco, Barrio Sarriena s/n 48940 Leioa, SPAIN.}
	\email{\textsuperscript{1} mahdi.zreik@math.u-bordeaux.fr or mzreik@bcamath.org.}
	\date{}
\subjclass[2010]{81Q10 , 81V05, 35P15, 58C40}
	
\begin{document}
	\begin{abstract}
Let $\Omega_+\subset\mathbb{R}^{3}$ be a fixed bounded domain with boundary $\Sigma = \partial\Omega_{+}$. We consider $\mathcal{U}^\varepsilon$ a tubular neighborhood of the surface $\Sigma$ with a thickness parameter $\varepsilon>0$, and we define the perturbed Dirac operator $\mathfrak{D}^{\varepsilon}_{M}=D_m +M\beta \mathbbb{1}_{\mathcal{U}^{\varepsilon}},$ with $D_m$ the free Dirac operator, $M>0$, and $\mathbbb{1}_{\mathcal{U }^{\varepsilon}}$ the characteristic function of $\mathcal{U}^{\varepsilon}$. Then, in the norm resolvent sense, the Dirac operator $\mathfrak{D}^{\varepsilon}_M$ converges to the Dirac operator coupled with Lorentz scalar $\delta$-shell interactions as $\varepsilon = M^{-1}$ tends to $0$, with a convergence rate of $\mathcal{O}(M^{-1})$.
	\end{abstract}
\maketitle	
\tableofcontents
\section{Introduction and Main results}
The aim of this work is to approximate the Dirac operator coupled with a singular $\delta$-interactions, supported on a closed surface. More precisely, our main goal in this article is to approximate the Dirac operator coupled with confining Lorentz scalar $\delta$-shell interactions (\emph{i.e.,} when $\eta =0$ and $\mu=\pm 2$ in \eqref{deltainter}, below) by a perturbed Dirac operator $\mathfrak{D}^{\varepsilon}_M=D_m + M\beta \mathbbb{1}_{\mathcal{U}^\varepsilon}$, where $D_m$ is the free Dirac operator, and $M$ is a large mass supported on a tubular neighborhood, $\mathcal{U}^\varepsilon$, with thickness $\varepsilon>0$. Working with this type of massive potential leads to the appearance of what we've seen in \cite{BBZ}, called Dirac operators with MIT bag boundary conditions, when the mass $M$ becomes large. In this paper we interested in establishing the convergence (for suitable relation between $\varepsilon$ and $M$: $\varepsilon=M^{-1}$, as $\varepsilon$ goes to $0$) of such perturbations to a direct sum of two MIT bag operators, which we denote by $D_{\mathrm{MIT}}^{\Omega_+}(m)$ and $D_{\mathrm{MIT}}^{\Omega_-}(m)$ (see Section \ref{DefMIT} for the exact notations), acting in the domains $\Omega_+$ and $\Omega_-:=\mathbb{R}^3\setminus \overline{\Omega_+}$, respectively. This decoupling of these MIT bag Dirac operators can be linked to the confining version of the Dirac operator coupled with purely Lorentz scalar $\delta$-shell interaction supported on the surface $\Sigma:=\partial\Omega_+$, which will be discussed briefly in the following part of the current paper.\\

\begin{figure}[h]
\includegraphics[scale=0.12]{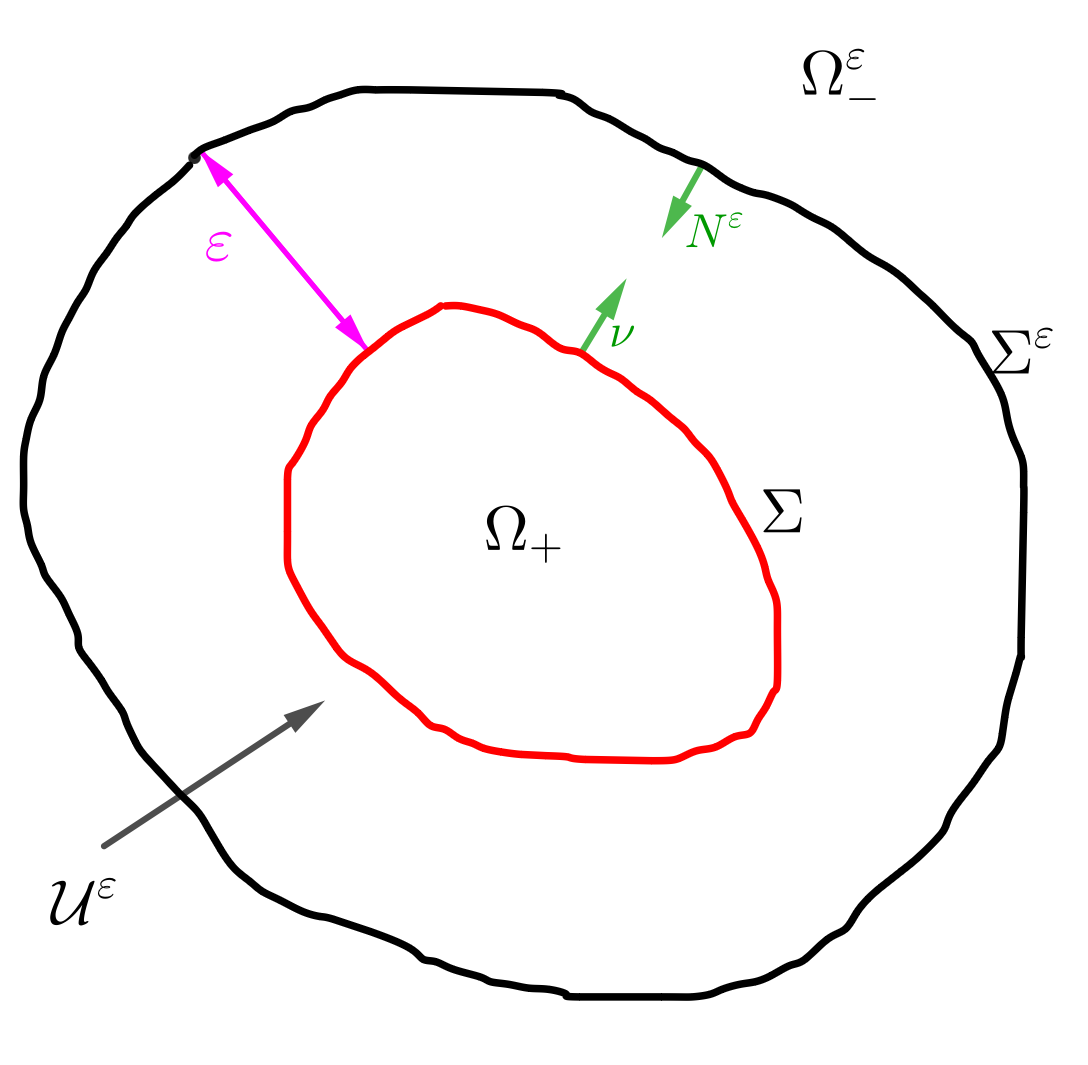} 
    \caption{Domain}
    \label{figure1}
    \end{figure}

The convergence of $\mathfrak{D}_M^\varepsilon$ to the MIT bag operator was established in \cite[Section 6]{BBZ}, in the norm resolvent sense, when $M$ tends to $+\infty$, and $\varepsilon$ fixed. However, in \cite{BBZ}, the mass $M$ is supported on an unbounded domain, which has only one boundary. Whereas, in the current work, $M$ is supported on a bounded domain with two boundaries, whose distance between them is the thickness $\varepsilon$, as shown in Figure \ref{figure1}. Thus, it is then natural to address the following question: Let $M$ be a large mass supported on a tubular vicinity of surface $\Sigma$. What happens when the thickness of the tubular tends to zero with $M^{-1}$?\\

The methodology followed, as in the problem of \cite{BBZ} study the pseudodifferential properties of the Poincaré-Steklov (PS) operators for the Dirac operator (\emph{i.e.}, an analogue of the Diricklet-to-Neumann operators for the Laplace operator). The complexity in the current problem is that these operators take a pair of functions with respect to $\partial\mathcal{U}^\varepsilon:= \Sigma \cup \Sigma^\varepsilon$ such that for all $x_\Sigma \in \Sigma,$ we have $\Sigma^\varepsilon\ni x= x_\Sigma + \varepsilon \nu(x_\Sigma)$, where $\nu$ is the unit normal to the surface $\Sigma$ pointing outside $\Omega_+$. So, we will control these operators by tracking the dependence on the parameter $\varepsilon$, and consequently, the convergence when $\varepsilon$ goes to $0$ and $M$
goes to $+\infty$. 
\\

Now, to give a rigorous definition of the operator we are dealing within this paper and to go into more details, we need to introduce some notations. For $m> 0$, the free Dirac operator $D_m$ in $\rr^3$ is defined by $D_m:=- i \alpha \cdot\nabla + m\beta$, with
\begin{align*}
	\alpha_j&=\begin{pmatrix}
	0 & \sigma_j\\
	\sigma_j & 0
	\end{pmatrix}\quad \text{for } j=1,2,3,
	\quad
	\beta=\begin{pmatrix}
	\mathbb{I}_2 & 0\\
	0 & -\mathbb{I}_2
	\end{pmatrix},\quad 
\mathbb{I}_2 :=\begin{pmatrix}
1& 0 \\
0 & 1
\end{pmatrix}, \\
&\text{and }\,
\sigma_1=\begin{pmatrix}
	0 & 1\\
	1 & 0
	\end{pmatrix},\quad \sigma_2=
	\begin{pmatrix}
	0 & -i\\
	i & 0
	\end{pmatrix} ,
	\quad
	\sigma_3=\begin{pmatrix}
	1 & 0\\
	0 & -1
	\end{pmatrix},
	\end{align*}  
the family of Dirac and Pauli matrices satisfying the anticommutation relations:
 \begin{align}\label{PMD}
   \lbrace \alpha_j,\alpha_k\rbrace = 2\delta_{jk}\mathbb{I}_4,\quad \lbrace \alpha_j,\beta\rbrace = 0,\quad \text{and} \quad \beta^{2} = \mathbb{I}_4,\quad j,k\in \lbrace 1,2,3\rbrace,  
 \end{align}
 where $\lbrace \cdot,\cdot\rbrace$ is the anticommutator bracket. As usual, we use the notation $\alpha \cdot x=\sum_{j=1}^{3}\alpha_j x_j$  for $x=(x_1,x_2,x_3)\in\rr^3$.  We recall that $D_m$ is self-adjoint in $\mathit{L}^2(\rr^3)^4$ with $\mathrm{dom}(D_m)=\mathit{H}^1(\rr^3)^4$ (see, \emph{e.g.}, \cite[ subsection~1.4]{Tha}), and that the spectrum is given by
	\begin{align*}
	\mathrm{Sp}(D_m) =\mathrm{Sp}_{\mathrm{cont}}(D_m)=(-\infty,-m]\cup [m,+\infty).
	\end{align*}

Let $\Omega_+$ be a bounded smooth domain in $\mathbb{R}^3,$ and $\Sigma:=\partial\Omega_+$ its boundary. For $(\eta,\mu)\in\mathbb{R}^2$, the three-dimensional Dirac operator with $\delta$-shell interactions is defined formally by
\begin{align}\label{deltainter}
    \mathbb{D}_{\eta,\mu}:f\mapsto D_m f + (\eta\mathbb{I}_4 + \mu\beta)\delta_{\Sigma} f,
\end{align}
where $\delta_\Sigma$ is the Dirac delta distribution supported on $\Sigma,$ and the constant $\eta$ (resp. $\mu$) measures the strength of the electrostatic (resp. Lorentz scalar) part of the interaction.  In this case, the operator in  \eqref{deltainter} is called the Dirac operator coupled with electrostatic and Lorentz scalar $\delta$-shell interactions.\\

The investigation of the properties of the Dirac operator  $\mathbb{D}_{\eta,\mu}$ goes back to the articles \cite{DES} and \cite{Do}. Furthermore, in \cite{DES}, the authors state that the shell becomes impenetrable if we assume that $\eta^2-\mu^2 = 4$ (known as the confinement case). Physically, this means that a particle such as an electron that is in the region $\Omega_+$ at time $t=0$ cannot cross the surface $\Sigma$ to reach the region $\mathbb{R}^3\setminus \overline{\Omega_+}$ as time progresses (and vice versa). Mathematically, this implies that we can decompose the considered Dirac operator into a direct sum of two operators acting respectively on $\Omega_+$ and $\mathbb{R}^3\setminus \overline{\Omega_+}$, each with the corresponding boundary conditions. If $\eta=0,$ physicists in particular have been aware of this phenomenon since the 1970s, when they considered confinement in hadrons with a model (see \cite{CJJTW} and \cite{KJ}). The mathematical model describing this, using the Dirac operator with MIT boundary conditions, has been extensively studied in mathematical papers such as those mentioned in \cite{BHM}. In our paper we refer to the Dirac operator, with MIT bag boundary conditions as $D^{\Omega_{+-}^{\varepsilon}}_{\mathrm{MIT}}(m)$ (see the beginning of Section \ref{DefMIT} for the exact definitions). 
\\

The approximation of the Dirac operators with regular/singular potential has been the subject of several recent mathematical papers. Therefore, in the one-dimensional case, the analysis is carried out in \cite{PS}, where Šeba showed that convergence in the sense of norm resolvent is true. In 2D case, \cite{CLMT} considered the approximation of Dirac operators with electrostatic, Lorentz scalar, and anomalous magnetic $\delta$-shell potentials on closed and bounded curves, in the non-critical and non-confinement cases. In 3D case, the authors of \cite{MP} showed an approximation of the Dirac operators coupled with $\delta$-shell interactions, however, a smallness assumption for the potential was required to achieve such a result. Finally, in 3D case, I have established in \cite{MZ1} an approximation of the operator $\mathbb{D}_{\eta,\tau}$, in terms of the strong resolvent, in the non-critical and non-confinement cases (\emph{i.e.}, when $\eta^2 - \mu^2  \neq\pm 4$) without the smallness assumed in \cite{MP}. Now, let us describe the main results of the present manuscript.
\subsection*{Description of main results.} Let $\Omega_+$ be a open bounded set in $\mathbb{R}^3$ with a compact smooth boundary $\Sigma:=\partial\Omega_+$,  let $\nu$ be the outward unit normal to  $\Omega_+$. Throughout the current paper, we shall work on the Hilbert space $L^2(\mathbb{R}^{3})^4$ \big(resp. $L^2(\O^{\varepsilon}_\pm)^4$ with $\Omega^{\varepsilon}_+=\Omega_+ \cup \mathcal{U}^{\varepsilon}$ and $\Omega^\varepsilon_- = \mathbb{R}^3\setminus \overline{\Omega^{\varepsilon}_+}\big)$ with respect to the Lebesgue measure, and we will make use of the orthogonal decomposition $L^2(\mathbb{R}^{3})^4=L^2(\O^{\varepsilon}_-)^4\oplus L^2(\O^{\varepsilon}_+)^4$. We denote by $ N^{\varepsilon}$ the outward unit normal with respect to  $\O_{-}^{\varepsilon}$. More precisely, for $\varepsilon_0$ sufficiently small, we assume that $\S$, $\O_{-}^{\varepsilon}$, $\S^{\varepsilon}$ and $\mathcal{U}^{\varepsilon}$ satisfied 
	\begin{equation}\label{conditionsgeo}
  \begin{aligned}
	    &\S^{\varepsilon}:=\lbrace x\in\mathbb{R}^3, \,x = x_{\S}+ \varepsilon \nu(x_{\S}): \, x_{\S}\in\S\rbrace,\\&
	     \Omega^{\varepsilon}_{-} =\lbrace x\in\mathbb{R}^{3},\, \mathrm{dist}(x,\S)>\varepsilon\rbrace,\\&
    \mathcal{U}^\varepsilon:=\lbrace x\in\mathbb{R}^{3},\,\,x=x_{\S}+t\,\nu(x_{\S}):\,x_{\S}\in\S\,\,\text{and}\,\,t\in(0,\varepsilon)\rbrace, \quad  \text{with }\varepsilon\in(0,\varepsilon_0).
	\end{aligned}
 \end{equation}	 
 In other words, the Euclidean space is divided as follows:
\begin{align*}
\rr^3=\O^{\varepsilon}_-\cup\Sigma^{\varepsilon}\cup \mathcal{U}^{\varepsilon}\cup\Sigma\cup\Omega_+.
\end{align*}
  We consider perturbations of the free Dirac operator $D_m$ in the whole space by a large mass $M$ term living in an  $\varepsilon$ neighborhood $\mathcal{U}^\varepsilon$ of $\Sigma$. The perturbed Dirac operator we are interesting on is $\mathfrak{D}_{M}^\varepsilon:= D_m + M\beta \mathbbb{1}_{\mathcal{U}^\varepsilon}$, where $\mathbbb{1}_{\mathcal{U}^\varepsilon}$ is the characteristic function of ${\mathcal{U}^\varepsilon}$ and $\varepsilon$ is the thickness of the tubular region $\mathcal{U}^{\varepsilon}.$ The results of the present article are presented as follows:\\

To establish the main result outlined in Theorem \ref{maintheo}, we must show the following approximations: 
\begin{proposition}\label{coro1.1}
We consider the confining version of the Dirac operator coupled with a purely Lorentz scalar $\delta$-shell interaction, denoted by $\mathscr{D}_{L}:=\mathbb{D}_{0,2}$ \big(\textit{i.e.,} when $\eta=0$ and $\mu=2$ in \eqref{deltainter}\big). Then, for any $z\in \rho(\mathscr{D}_{L})$ and $\varepsilon$ sufficiently small, the following estimate holds:

\begin{align}\label{stat}
\left|\left| e_{\Omega_{+-}^{\varepsilon}}R^{\Omega^\varepsilon_{+-}}_{\mathrm{MIT}} (z)r_{\Omega_{+-}^{\varepsilon}} - R_{L}(z)\right|\right|_{\mathit{L}^2(\rr^3)^4 \rightarrow \mathit{L}^2(\rr^3)^4}& =\mathcal{O}(\varepsilon)\quad\text{as}\quad  \varepsilon\rightarrow 0.
\end{align} 
where $R^{\Omega^\varepsilon_{+-}}_{\mathrm{MIT}}$ is the resolvent of the direct sum of both MIT bag operators, refer to $D^{\Omega_{+-}^{\varepsilon}}_{\mathrm{MIT}}(m)$ and which will be defined rigorously in Section \ref{DefMIT},  $R_L$ is the resolvent of the Dirac operator coupled with purely Lorentz scalar $\delta$-shell interactions, $\mathscr{D}_L$, and $r_{\Omega_{+-}^{\varepsilon}}$ resp.  $e_{\Omega_{+-}^{\varepsilon}}$ is the restriction operator in $\Omega_{+-}^{\varepsilon}:=\Omega_+ \cup\Omega_-^\varepsilon$ resp. its adjoint operator, \textit{i.e.}, the extension by $0$ outside of $\Omega_{+-}^{\varepsilon}$.
\begin{remark}
We mention that the proof of Proposition \ref{coro1.1} is not difficult to realize. Indeed, we establish the above approximation by tracking the dependence on the thickness $\varepsilon$, when $\varepsilon$ goes to $0$. However, what is important to achieve is the proof of the following proposition, for which studies and estimates are required by tracking the dependence on the parameters $\varepsilon$ and $M$, in order to establish such a relationship between the parameters, and prove therefore the main result of Theorem \ref{maintheo}.    
\end{remark}
\end{proposition}
\begin{proposition}\label{Prop1.1}
   Let $K\subset\mathbb{C}\setminus\mathbb{R} $ be a compact set. Then, there is $M_0>0$ such that for all $M>M_0$ and $\varepsilon=M^{-1}$: $K\subset\rho(\mathfrak{D}^\varepsilon_M)$ and for all $z\in K$, the following estimate holds on the whole space
    \begin{align*}
       \big|\big| R_M^\varepsilon(z) - e_{\Omega_{+-}^{\varepsilon}}R^{\Omega^\varepsilon_{+-}}_{\mathrm{MIT}}(z)r_{\Omega_{+-}^{\varepsilon}} \big|\big|_{L^2(\mathbb{R}^3)^4 \rightarrow L^2(\mathbb{R}^3)^4} = \mathcal{O}(M^{-1}).
    \end{align*}
\end{proposition}
The latter proposition means that the Dirac operator $\mathfrak{D}_M^\varepsilon$ is approximated, in the norm resolvent sense, by both MIT bag Dirac operators, acting in $\Omega_{+-}^{\varepsilon}$ with a rate of $\mathcal{O}(M^{-1})$ when $M$ tends to $\infty$. \\

By combining Propositions \ref{coro1.1}, \ref{Prop1.1}, we arrive  at  the following main result:
\begin{theorem}\label{maintheo} Let $z\in\rho(\mathscr{D}_{L})$, then for $M$ sufficiently large, $z\in\rho(\mathfrak{D}_{M}^{\varepsilon})$,  and $\varepsilon=M^{-1}$, the following estimate holds
\begin{align*}
\left|\left| R^{\varepsilon}_M(z)- R_{L}(z)\right|\right|_{\mathit{L}^2(\rr^3)^4}&=\mathcal{O}\left( M^{-1}\right).
\end{align*} \qed
\end{theorem}
The most important ingredient in proving Proposition \ref{coro1.1} is the use of the Krein formula of the resolvents of $\mathscr{D}_{L}$ and both MIT bag operators, $D_{\mathrm{MIT}}^{\Omega_+}$ and $D_{\mathrm{MIT}}^{\Omega_-^\varepsilon}$ (see Section \ref{kreinresolform}), acting in  $L^2(\Omega_+)^4$ and $L^2(\Omega_-^\varepsilon)^4$, respectively. Then, in Proposition \ref{RMITCONV}, we establish that the convergence $(D_{\mathrm{MIT}}^{\Omega_-^\varepsilon} -z)^{-1}$ toward $(D_{\mathrm{MIT}}^{\Omega_-} -z)^{-1}$ holds for any non-real $z$, when $\varepsilon$ goes to $0,$ and we then obtain, in the norm resolvent sense, the convergence of $D_{\mathrm{MIT}}^{\Omega^\varepsilon_{+-}}:=D_{\mathrm{MIT}}^{\Omega_+}\oplus D_{\mathrm{MIT}}^{\Omega_-^\varepsilon}$ to $\mathscr{D}_L = D_{\mathrm{MIT}}^{\Omega_+}\oplus D_{\mathrm{MIT}}^{\Omega_-} $.\\\\ 
The key point to establish the result of Proposition \ref{Prop1.1} is to treat the elliptic problem $(\mathfrak{D}_M^\varepsilon - z)\mathfrak{U} = f \in  L^2(\mathbb{R}^3)^4$ as a transmission problem (where $P_\pm t_{\Sigma}\mathfrak{U}_{|_{\Omega_+}} = P_\pm t_{\Sigma} \mathfrak{U}_{|_{\mathcal{U}^\varepsilon}}$  and $P^\varepsilon_\pm t_{\Sigma^\varepsilon}\mathfrak{U}_{|_{\Omega_-^\varepsilon}} = P^\varepsilon_\pm t_{\Sigma^\varepsilon} \mathfrak{U}_{|_{\mathcal{U}^\varepsilon}}$ are the transmission conditions) and to use the semiclassical properties of the auxiliary operator $\Upsilon_M^\varepsilon (z)$ acting on the boundary $\partial\mathcal{U}^\varepsilon=\Sigma\cup\Sigma ^\varepsilon$, which is constructed by the Poincaré-Steklov operators (see \eqref{kfullreso} for the exact notation). Indeed, in Section \ref{convergences}, we show convergence of the Dirac operator, $\mathfrak{D}_M^\varepsilon$, to both MIT bag operators, $D_{\mathrm{MIT}}^{\Omega_+}$ and $D_{\mathrm{MIT}}^{\Omega_-^\varepsilon}$, with a convergence rate of $\mathcal{O}(M^{-1})$ for $M=\varepsilon^{-1}$ sufficiently large. Consequently, using these ingredients, a kind of convergence can be established in Theorem \ref{maintheo} for $\varepsilon = M^{-1}$.\\

Unlike the application in paper \cite[Theorem 6.1]{BBZ}, we mention that in this problem the operator $\Upsilon_M^\varepsilon$ (which is constructed by the Poincaré-Steklov operators) takes a pair of functions with respect to $\partial\mathcal{U}^\varepsilon$.\\\\
We note that $P_{\pm}^\varepsilon$ and $P_{\pm}$ are the orthogonal projections with respect to $N^\varepsilon$ and $\nu$, respectively, defined by 
\begin{align}\label{projections}
P^\varepsilon_{\pm}:=(\mathbb{I}_{4}\mp i\beta\alpha\cdot N^\varepsilon)/2\quad\text{and}\quad 
    P_{\pm}:=(\mathbb{I}_{4}\mp i\beta\alpha\cdot \nu)/2.
\end{align}

We end this part with the following remark on the projections $P_{\pm}$ and $P_{\pm}^\varepsilon$:
\begin{remark}\label{remark3.1}
We define the diffeomorphism $p:\S\longrightarrow \S^{\varepsilon}$ such that for all $x_{\S}\in\S$, we get $p(x_{\S}):=x_{\S}+ \varepsilon \nu(x_{\S})=x$. Then, we have $$N^{\varepsilon}(x)= - (\nu\circ p^{-1})(x)=-\nu(x_\Sigma),$$ with  
\begin{align*}
    P_{\pm}^{\varepsilon}(x)&=\dfrac{1}{2}\left ( \mathbb{I}_4 \mp i\beta\alpha\cdot N^{\varepsilon}_{+}(x)\right)=\dfrac{1}{2}\left ( \mathbb{I}_4 \pm i\beta\alpha\cdot \nu(x_\Sigma)\right):= P_{\mp}\circ p^{-1}(x) =P_{\mp}(x_\Sigma).
\end{align*}
\end{remark}
\subsection*{Organization of the paper.} The present paper is structured as follows. Section \ref{prelimi} is dedicated to the preliminaries and the MIT bag operators, where we give some notations and definitions, and we recall some basic properties of boundary integral operators associated with $(D_m - z)$. Moreover, in this section we set up some geometric aspects characterizing our domains, define the Dirac operator with MIT bag boundary conditions and give some properties. Section \ref{section4} is devoted to the study of pseudodifferential properties of the Poincaré-Steklov operators, where the main result are Proposition \ref{prop4.6} and Corollary \ref{corollsry5.1}. In Section \ref{app}, we set up a Krein formula connecting the resolvents of $\mathfrak{D}_M^\varepsilon$ with those of $D_{\mathrm{MIT}}^{\Omega_{+-}^{\varepsilon}}$. With its help, in Section \ref{convergences} turns out that a kind of convergence can be achieved for $\varepsilon = M^{-1}$, with a convergence rate of $\mathcal{O}(M^{-1})$ as $M$ becomes large (\emph{i.e.,} $\varepsilon \in (0,\varepsilon_0)$ sufficiently small). Therefore, we show the main results of this paper: in the proof of Proposition \ref{coro1.1}, we approximate the resolvent of MIT bag operators with that of the Dirac operator coupled with purely Lorentz scalar $\delta$-shell interactions, in the norm resolvent sense, with a convergence rate of $\mathcal{O}(\varepsilon)$, and we prove Proposition \ref{Prop1.1} on the convergence of the resolvent of $\mathfrak{D}_M^\varepsilon$ to those of the MIT bag operators, $D_{\mathrm{MIT}}^{\Omega_{+-}^{\varepsilon}}(m)$, for $M$ sufficiently large. 
\section{Setting and $\mathrm{MIT}$ bag operator}\label{prelimi}
In this section we gather some well-known results about boundary integral operators. Before proceeding further, however, we need to introduce some notations that we will use in what follows.\\
We define the unitary Fourier–Plancherel operator $\mathscr{F}: \mathit{L}^2(\rr^{\mathit{d}})^4  \longrightarrow \mathit{L}^2(\rr^{\mathit{d}})^4$ as follows:
\begin{align*}
 \hat{u}(\xi):=\mathscr{F}[u](\xi)=(2\pi)^{-d}\int_{\rr^d}e^{-ix\cdot\xi}u(x)\mathrm{d}x,\quad \forall\xi \in\rr^d.
\end{align*}
For $\overline{x}\in\rr^{d-1}$,  we will abbreviate the partial Fourier transform on the variable  $\overline{x}$ with $\mathscr{F}_{\overline{x}}$. Given $s\in[0,1]$, we define the usual  Sobolev space $\mathit{H}^s(\rr^{\mathtt{d}})^4$  as
\begin{align*}
 \mathit{H}^s(\rr^{\mathit{d}})^4:=\{ u\in\mathit{L}^2(\rr^{\mathit{d}})^4: \int_{\rr^{\mathit{d}}}(1+|\xi|^2)^{s} \left|\mathscr{F}[u](\xi)\right|^2\mathrm{d}\xi<\infty\},
\end{align*}
and for a bounded or unbounded Lipshitz domain $\Omega\subset \rr^3$, we write $\partial\Omega:=\S$ for its boundary and we denote by $\nu$ and
$\sigma$ the outward pointing normal to $\O$ and the surface measure on $\S$, respectively. By $L^{2}(\rr^3)^4:= L^{2}(\rr^3,\cc^4)$ (resp. $L^{2}(\O)^4:= L^{2}(\O,\cc^4)$) we denote the usual $L^{2}-$space over $\rr^3$ (resp. $\O$), and we let $r_{\O}:L^{2}(\rr^3)^4\longrightarrow L^{2}(\O)^4$ be the restriction operator on $\O$ and $e_{\O}:L^{2}(\O)^4\longrightarrow L^{2}(\rr^3)^4$ its adjoint operator, \emph{i.e.}, the extension by 0 outside
of $\O$. Now, we let $H^{1}(\O)^{4}$ to be the first order Sobolev space 
\begin{align*}
\mathit{H}^1(\O)^4=\{ \varphi\in\mathit{L}^2(\O)^4: \text{ there exists } \tilde{\varphi}\in\mathit{H}^1(\rr^{3})^4 \text{ such that }  \tilde{\varphi}|_{\O} =\varphi\}.
\end{align*}
By $\mathit{L}^2(\S)^4:=\mathit{L}^2(\S, \mathrm{d}\sigma)^4$ we denote the usual $\mathit{L}^2$-space over $\S$. The Sobolev space of order $1/2$ along the boundary, $\mathit{H}^{1/2}(\S)^4$,  consists of all functions $g\in\mathit{L}^2(\S)^4$ for which 
	 \begin{align*}
 \|g \|^2_{\mathit{H}^{1/2}(\S)^4}:= \int_{\S}|g(x)|^2\mathrm{d}\sigma(x) + \int_{\S}\int_{\S}\frac{|g(x)-g(y)|^2}{|x-y|^{3}}\mathrm{d}\sigma(y) \mathrm{d}\sigma(x)<\infty.
  \end{align*}
As usual we let $\mathit{H}^{-1/2}(\S)^4$ to be the dual space of $\mathit{H}^{1/2}(\S)^4$. We denote by  $t_{\S}:\mathit{H}^1(\O)^4\rightarrow \mathit{H}^{1/2}(\S)^4$ the classical trace operator, and by  $\mathcal{E}_{\O}:  \mathit{H}^{1/2}(\S)^4\rightarrow  \mathit{H}^1(\O)^4$ the extension operator, that is 
$$t_{\S}\mathcal{E}_{\O}[f]=f,\quad \forall f\in \mathit{H}^{1/2}(\S)^4.$$
\subsection{Boundary integral operators associated with the free Dirac operator}\label{Subintope} The aim of this part is to introduce boundary integral operators associated  to the fundamental solution of $D_m$  and to summarize some of their well-known properties.  In this section, $\O$ is a bounded domain in $\rr^3$ with $\S:=\partial\O$ its boundary and we denote by $\nu$ the outward pointing normal to $\O$. We set $\O_+:=\Omega$ and $\O_-=\rr^3\setminus\overline{\O_{+}}.$\\         
For $z\in\mathbb{C}\backslash(-\infty,-m]\cup[m,+\infty)$, with the convention that $\mathrm{Im}\sqrt{z^2-m^2}>0$,  the fundamental solution of $(D_m-z)$ is given by
\begin{align}\label{defsolo}
	\phi^{z}_{m}(x)=\frac{e^{i\sqrt{z^2-m^2}|x|}}{4\pi|x|}\left(z +m\beta+( 1-i\sqrt{z^2-m^2}|x|)i\alpha\cdot\frac{x}{|x|^2}\right), \quad \forall  x\in\rr^3\setminus\{0\}.
\end{align}
 We define the potential operator $\Phi^{z}_{m}: \mathit{L}^2(\S)^4  \longrightarrow \mathit{L}^2(\rr^3)^4$ by
\begin{align*}\label{}
  \Phi^{z}_{m}[f](x):=\int_\S \phi^{z}_{m}(x-y)f(y)\mathrm{d}\sigma(y), \quad \text{for all } x\in\rr^3\setminus\S.
\end{align*}
Furthermore, $(D_{m} - z)\Phi^{z}_{m}[f]=0$ holds in $\mathcal{D}^{'}(\Omega_{\pm})^4$, for all $f\in L^2(\S)^4$. 
Finally, given $x\in\S$ we define the Cauchy operators $\mathscr{C}^{z}_{m}: \mathit{L}^2(\S)^4  \longrightarrow \mathit{L}^2(\S)^4$ as the singular integral operator acting as
\begin{align}\label{CauchyOpe}
\mathscr{C}^{z}_{m}[f](x)&:=  \lim\limits_{\rho\searrow 0}\int_{|x-y|>\rho}\phi^{z}_{m}(x-y)f(y)\mathrm{d}\sigma(y), \quad \text{for } \mathrm{d}\sigma\emph{-a.e.,}\,x\in\S,\,f\in\mathit{L}^2(\S)^4,
\end{align} 
and the following bounded operator $\mathit{C}^{z}_{\pm,m}: \mathit{L}^2(\S)^4  \longrightarrow \mathit{L}^2(\S)^4$ as follows:
\begin{align*}
    \mathit{C}^{z}_{\pm,m}[f](x):=\lim_{\O_{\pm}\ni y\,\,\overset{nt}{\rightarrow}\,\, x}\Phi^{z}_{m}[f](y),
\end{align*}
where $\O_{\pm}\ni y\,\,\overset{nt}{\rightarrow}\,\, x$ means that $y$ tends to $x$ non-tangentially from $\O_+$ and $\O_-,$  respectively, \emph{i.e.}, for $y\in\O_{\pm},$ we get $|x-y|<(1+a)\mathrm{dist}(y,\S)$ for $a>0$ and $x\in\S$.\\ 

It is well known that $\Phi^{z}_{m}$ and $\mathscr{C}^{z}_{m}$ are bounded and everywhere defined (see  \cite[Section 2]{AMV1}), and that 
\begin{align*}
((\alpha\cdot \nu)\mathscr{C}^{z}_{m})^2=\left(\mathscr{C}^{z}_{m}(\alpha\cdot \nu)\right)^2=-\dfrac{1}{4}\mathbb{I}_4, \quad \forall z\in\rho(D_m),
\end{align*}
holds in $ \mathit{L}^2(\S)^4$,  cf. \cite[Lemma 2.2]{AMV2}. In particular, the inverse $(\mathscr{C}^{z}_{m})^{-1}=-4(\alpha\cdot \nu)\mathscr{C}^{z}_{m}(\alpha\cdot \nu)$ exists and is bounded and everywhere defined. Note that $\phi^{z}_m(y-x)^*=\phi^{\overline{z}}_m(x-y)$, as a consequence $(\mathscr{C}^{z}_{m})^{\ast}=\mathscr{C}^{\overline{z}}_{m}$ holds in $ \mathit{L}^2(\S)^4$. In particular, $\mathscr{C}^{z}_{m}$ is self-adjoint in $\mathit{L}^2(\S)^4$ for all $z\in(-m,m)$.\\

Now, we define the operator  $\Lambda^{z}_{\pm,m}$ by
\begin{align*}
\Lambda^{z}_ {\pm,m}=\frac{1}{2}\beta \pm \mathscr{C}^{z}_{m},\quad\text{ for all }  z\in\rho(D_m),
\end{align*}
which is clearly a bounded operator from $\mathit{L}^{2}(\S)^4$ into itself. 

In the next lemma,  we collect  the main  properties of the operators $\Phi^{z}_{m}$, $\mathscr{C}^{z}_{m}$ and $\Lambda^{z}_{\pm,m}$.
\begin{lemma}{\cite[Lemma 2.1]{BBZ}}.\label{prop of C} Given $z\in\rho(D_m)$ and let $\Phi^{z}_{m}$, $\mathscr{C}^{z}_{m}$ and $\Lambda^{z}_{\pm,m}$ be as above. Then the following holds true:
\begin{itemize}
   \item[($\mathrm{i}$)] The operator $\Phi^{z}_{m}$ is bounded from $\mathit{H}^{1/2}(\S)^4$ to $\mathit{H}^{1}(\O)^4$, and the following Plemelj-Sokhotski jump formula holds that 
   \begin{align*}\label{}
   t_\S \Phi^{z}_{m}|_{\Omega_\pm}[f]=  C_{\pm,m}^z[f]=\left(\mp\frac{i}{2}(\alpha\cdot \nu) + \mathscr{C}_{m}^{z}\right)[f], \quad \forall f\in \mathit{H}^{1/2}(\S)^4.
   \end{align*}
  \item[($\mathrm{ii}$)] The operator $\mathscr{C}^{z}_{m}$ gives rise to a bounded operator $\mathscr{C}_{m}^{z} : \mathit{H}^{1/2}(\S)^4  \longrightarrow \mathit{H}^{1/2}(\S)^4$. 
  \item[($\mathrm{iii}$)]  The operator $\Lambda^{z}_{\pm, m} : \mathit{H}^{1/2}(\S)^4  \longrightarrow \mathit{H}^{1/2}(\S)^4$ is bounded invertible for all  $z\in\rho(D_m)$.\qed 
\end{itemize}
\end{lemma}
The last thing in this section is the definition of the Dirac operator coupled with purely Lorentz scalar $\delta$-interaction. 
\begin{definition}\label{DiracLorentz}
Let $\mu\in\rr\setminus\lbrace 0\rbrace$. The Dirac operator coupled with purely Lorentz scalar $\delta$-shell interaction of strength $\mu$, is the operator $\mathbb{D}_{0,\mu}$, acting in $L^{2}(\rr^3)^{4}$ and defined on the following domain
\begin{align}\label{DOMDL}
    \mathrm{dom}(\mathbb{D}_{0,\mu}):=\lbrace \varphi=u+\Phi^{z}_{m}[g],\,\,u\in H^{1}(\rr^3)^{4},\,g\in L^{2}(\S)^{4},\, t_{\S}u=-\Lambda_{+,m}^{z}[g] \,\, \text{on }\Sigma\rbrace.
\end{align}
Hence, $\mathbb{D}_{0,\mu}$ acts in the sense of distributions as $\mathbb{D}_{0,\mu}(\varphi)=D_m u$, for all $\varphi=u+\Phi_{m}^{z}[g]\in\mathrm{dom}(\mathbb{D}_{0,\mu}).$ 
Consequently, we can identify $\mathbb{D}_{0,\mu}$ as
\begin{align*}
  \mathbb{D}_{0,\mu}&=D_{m}\varphi_{-}\oplus D_{m}\varphi_{+},\\&
    \mathrm{dom}(\mathbb{D}_{0,\mu})=\lbrace w_{\pm}+\Phi_{m,\pm}^{z}[g],\,w_{\pm}\in H^{1}(\Omega_{\pm})^{4},\, g\in L^{2}(\S)^{4},\\& \hspace{6cm}P_{\pm}(t_{\S}w_{\pm}+C_{\pm,m}^{z}[g])=0,\text{ with } t_{\S}w_{\pm}=-\Lambda^{z}_{\pm,m}[g] \text{ on }\Sigma\rbrace,
 \end{align*}
where $\Phi_{m,\pm}^{z}[g]:L^{2}(\S)^{4}\longrightarrow L^{2}(\O_{\pm})^{4}$ is the operator defined by $\Phi_{m,\pm}^{z}[g](x)=\Phi_{m}^{z}|_{\Omega_{\pm}}[g](x),$ for $g\in L^{2}(\S)^{4}$ and $x\in \O_{\pm}.$ \\ Moreover, recall that $\mathbb{D}_{0,\mu}$ is a   self-adjoint operator on $H^{1}(\rr^3)^4$ for all $\mu\in\rr$ (see,  \cite[Section 5.1]{AMV2}), and for all $z\in\mathbb{C}\setminus\mathbb{R}$, the following resolvent formula holds \cite[Proposition 4.1]{BB1}
\begin{align*}\label{}
		(\mathbb{D}_{0,\mu}-z)&= (D_m -z)^{-1} - \Phi_{m}^{z}(\Lambda_{+,m}^{z})^{-1}t_{\S}(D_m -z)^{-1}.
	\end{align*}
\end{definition}
\subsection{Definition and some properties of the MIT bag operator.}\label{DefMIT} 
Recall the definition of the perturbed Dirac operator $\mathfrak{D}_{M}^{\varepsilon}:= D_{m}+M\beta \mathbbb{1}_{\mathcal{U}^{\varepsilon}}$, where $\mathbbb{1}_{\mathcal{U}^{\varepsilon}}$ is the characteristic function of $\mathcal{U}^{\varepsilon}$.
Then, we consider the MIT bag operators, $D^{\Omega_+}_{\mathrm{MIT}}(m)$ and  $D^{\Omega^\varepsilon_-}_{\mathrm{MIT}}(m)$, acting in $\Omega_+$ and $\Omega^\varepsilon_-$, respectively, and defined on the following domains\\
\begin{align*}
D^{\Omega_+}_{\mathrm{MIT}}(m) v_+ = D_m v_+,\quad \forall v_+\in\mathrm{dom}(D^{\Omega_+}_{\mathrm{MIT}}(m))=\lbrace v_{+}\in H^{1}(\Omega_{+})^{4},\quad  P_-t_{\Sigma}v_{+}=0 \text{ on } \Sigma\rbrace,
\end{align*}
\begin{align*}
D_{\mathrm{MIT}}^{\Omega^\varepsilon_-}(m) v^\varepsilon = D_m v^\varepsilon,\quad \forall v^\varepsilon\in\mathrm{dom}(D^{\Omega^\varepsilon_-}_{\mathrm{MIT}}(m))=\lbrace v^\varepsilon\in H^{1}(\Omega_{-}^{\varepsilon})^{4},\quad P^{{\varepsilon}}_{-}t_{\Sigma^{\varepsilon}}v^{\varepsilon}_-=0 \text{ on } \Sigma^\varepsilon\rbrace.
\end{align*}

Then, let the MIT Dirac operator, $D^{\Omega_{+-}^{\varepsilon}}_{\mathrm{MIT}} = D^{\Omega_+}_{\mathrm{MIT}}\oplus D^{\Omega^{\varepsilon}_-}_{\mathrm{MIT}},$ acts in $\Omega_{+-}^{\varepsilon} := \Omega_+ \cup \Omega_-^\varepsilon$, and defined on the following domain
\begin{align*}
\mathrm{dom}(D^{\Omega_{+-}^{\varepsilon}}_{\mathrm{MIT}})=\lbrace v^{\varepsilon}=(v^{\varepsilon}_{-},v_{+})\in H^{1}(\Omega_{-}^{\varepsilon})^{4}\oplus H^{1}(\Omega_{+})^{4},\quad P^{{\varepsilon}}_{-}t_{\Sigma^{\varepsilon}}v^{\varepsilon}_-=0= P_-t_{\Sigma}v_{+}\rbrace,
\end{align*}
with $D_{\mathrm{MIT}}^{\Omega_{+-}^{\varepsilon}}v^{\varepsilon}=(D_{+}\oplus D_{-})v^{\varepsilon} \,\,; \,\,D_{+}=D_{-}=D_{m}$ for all $v^{\varepsilon}\in\mathrm{dom}(D^{\Omega_{+-}^{\varepsilon}}_{\mathrm{MIT}}),$  and where the boundary condition holds in $H^{1/2}(\S^{\varepsilon})^4$ and $H^{1/2}(\S)^4$, respectively. Here, we recall that $P^{\varepsilon}_{\pm}$ and $P_{\pm}$ are the projections given in \eqref{projections}. \\

Finally, on $\mathcal{U}^\varepsilon$, we introduce the following Dirac auxiliary operator
\begin{align*}  D_{\mathrm{MIT}}^{\mathcal{U}^\varepsilon}(m+M)u^{\varepsilon}&= D_{m+M}u^{\varepsilon},\\ u^{\varepsilon}\in\text{dom}\big(D_{\mathrm{MIT}}^{\mathcal{U}^\varepsilon}(m+M)\big)&=\lbrace u^{\varepsilon}\in H^{1}(\mathcal{U}^{\varepsilon})^{4},\, P^{\varepsilon}_{+}t_{\Sigma^{\varepsilon}}u^{\varepsilon}=0=P_{+}t_{\Sigma}u^{\varepsilon}\,\text{on}\,\, \partial\mathcal{U}^{\varepsilon}:=\Sigma\cup\Sigma^{\varepsilon}\rbrace,
\end{align*}
with $D_{m+M}= D_m + M\beta = -i\alpha\cdot\nabla + (m+M)\beta.$ 
We note that  $D^{\mathcal{U}^\varepsilon}_{\mathrm{MIT}}$ is the MIT bag operator on $\mathcal{U}^{\varepsilon}$. 
\begin{theorem}\label{MITLipshictz} The operators $(D^{\Omega_+}_{\mathrm{MIT}},\mathrm{dom}(D^{\Omega_+}_{\mathrm{MIT}}))$ (resp. $(D^{\Omega^\varepsilon_-}_{\mathrm{MIT}},\mathrm{dom}(D^{\Omega^\varepsilon_-}_{\mathrm{MIT}}))$ and  $(D^{\mathcal{U}^\varepsilon}_{\mathrm{MIT}},\mathrm{dom}(D^{\mathcal{U}^\varepsilon}_{\mathrm{MIT}}))$) are self-adjoint and we have 
\begin{align*}
(D^{\Omega_+}_{\mathrm{MIT}}-z)^{-1}= r_{\Omega_{+}}(D_m -z)^{-1}e_{\Omega_{+}} - \Phi_{m,+}^{z}(\Lambda^z_{+,m})^{-1}t_{\Sigma}(D_m -z)^{-1}e_{\Omega_{+}}, \quad \forall z\in\rho(D_m).
\end{align*}
Moreover, the following statements hold true:
\begin{itemize}
\item[(i)]  $\mathrm{Sp}(D^{\Omega_+}_{\mathrm{MIT}})=\mathrm{Sp}_{\mathrm{disc}}(D^{\Omega_+}_{\mathrm{MIT}})\subset\rr\setminus[-m,m]$. \Big(Similarly for $D^{\mathcal{U}^\varepsilon}_{\mathrm{MIT}}$  for $(m+M)$ instead of $m$\Big).
\item[(ii)] $\mathrm{Sp}(D^{\Omega^\varepsilon_-}_{\mathrm{MIT}})=\mathrm{Sp}_{\mathrm{ess}}(D^{\Omega^\varepsilon_-}_{\mathrm{MIT}})=(-\infty,-m]\cup [m,+\infty)$. Moreover, if $\Omega_-^{\varepsilon}$ is connected then $\mathrm{Sp}(D^{\Omega^\varepsilon_-}_{\mathrm{MIT}})$ is purely continuous.
\item[(iii)]	Let $z\in\rho(D^{\Omega^\varepsilon_-}_{\mathrm{MIT}})$ be such that $2|z|<(m+M)$, then for all $f\in\mathit{L}^2(\mathcal{U}^{\varepsilon})^4$, it holds that
 \begin{align*}\label{}
     \left\|(D^{\mathcal{U}^\varepsilon}_{\mathrm{MIT}}-z)^{-1}f\right\|_{\mathit{L}^{2}(\mathcal{U}^{\varepsilon})^4}\lesssim M^{-1} \left\|f\right\|_{\mathit{L}^2(\mathcal{U}^{\varepsilon})^4},
      \end{align*}
\end{itemize}
uniformly with respect to $\varepsilon\in(0,\varepsilon_0).$
\end{theorem}
\textbf{Proof.} The proof of this theorem follows the same arguments as the proof of \cite[Theorem 3.1]{BBZ}, where the estimates are valid uniformly with respect to $\varepsilon$.\qed
\begin{definition}\label{def}
Let  $ z\in\rho(D_m)\cap \rho(D^{\mathcal{U}^\varepsilon}_{\mathrm{MIT}})  )$, 
$g^{\varepsilon}\in P^{{\varepsilon}}_-\mathit{H}^{1/2}(\S^{\varepsilon})^4$, $g_+\in P_-\mathit{H}^{1/2}(\S)^4$ and  $(h^{\varepsilon},h_+)\in P^{{\varepsilon}}_+\mathit{H}^{1/2}(\S^{\varepsilon})^4 \\ {\oplus P_+\mathit{H}^{1/2}(\S)^4}$. We denote by $E_{m}(z): P_-\mathit{H}^{1/2}(\S)^4\rightarrow \mathit{H}^{1}(\O^{+})^4$, respectively, $E^{\varepsilon}_{m}(z): P^{\varepsilon}_-\mathit{H}^{1/2}(\S^{\varepsilon})^4\rightarrow \mathit{H}^{1}(\O^{\varepsilon}_{-})^4$ the unique solution of the boundary value problem:
	\begin{equation}\label{L1b}
	\left\{
	\begin{aligned}
	(D_m-z)v_+&=0, \quad  &\text{ in } \Omega_{+},\\
P_-t_{\S}v_{+}&= g_{+},  \quad &\text{ on } \S,
 \end{aligned}
	\right.
	\end{equation}
 \begin{equation}\label{L1b'}
	\left\{
	\begin{aligned}
 (D_m-z)v_{-}^{\varepsilon}&=0, \quad &\text{ in } \O^{\varepsilon}_{-},\\
 P^{{\varepsilon}}_-t_{\S^{\varepsilon}_{-}}v^{\varepsilon}_-&= g^{\varepsilon},  \quad &\text{ on } \S^{\varepsilon}.
	\end{aligned}
	\right.
	\end{equation}
Similarly, we denote by $\mathcal{E}^{\varepsilon}_{m+M}(z): P^{{\varepsilon}}_+\mathit{H}^{1/2}(\S^{\varepsilon})^4 \oplus P_+\mathit{H}^{1/2}(\S)^4\rightarrow \mathit{H}^{1}(\mathcal{U}^{\varepsilon})^4$   the unique solution of the boundary value problem:  
\begin{equation}\label{L2b}
\left\{
\begin{aligned}
(D_{m+M}-z)u^{\varepsilon}=0, \quad &\text{ in } {\mathcal{U}^{\varepsilon}},\\
P^{{\varepsilon}}_+t_{\Sigma^{\varepsilon}} u^{\varepsilon}= h^{\varepsilon},  \quad &\text{ on } \Sigma^{\varepsilon},\\
P_+t_{\Sigma} u^{\varepsilon}= h_{+},  \quad &\text{ on } \Sigma.
\end{aligned}
\right.
\end{equation}
Define the Poincaré-Steklov operators associated with the above problems by  \begin{equation*}
   \begin{aligned}
\begin{array}{rcl}
       \mathscr{A}_{m}(z): P_- H^{1/2}(\S)^{4}&\to & P_+ H^{1/2}(\S)^{4}\\ 
         g_{+} &\mapsto & \mathscr{A}_{m}(z)g_{+}:=P_+ t_{\S} E_{m}(z)P_-g_{+},
       \end{array}\\
\begin{array}{rcl}
       \mathscr{A}^{\varepsilon}_{m}(z): P^\varepsilon_- H^{1/2}(\S^{\varepsilon})^{4}&\to & P^\varepsilon_+ H^{1/2}(\S^{\varepsilon})^{4}\\ 
         g^{\varepsilon}_- &\mapsto & \mathscr{A}^{\varepsilon}_{m}(z)g^{\varepsilon}:=P^\varepsilon_+ t_{\S} E^\varepsilon_{m}(z)P^\varepsilon_{-}g^{\varepsilon},
       \end{array}
       \end{aligned}
        \end{equation*}
      \begin{align*}\label{}
          \mathcal{A}^{\varepsilon}_{m+M}(z):  P_+H^{1/2}(\S)^{4}\oplus P^\varepsilon_+H^{1/2}(\S^{\varepsilon})^{4}\to  P_-H^{1/2}(\S)^{4}\oplus P^\varepsilon_-H^{1/2}(\S^{\varepsilon})^{4},\quad \text{with }
      \end{align*}
      $\quad\quad\quad\hspace{2cm} \mathcal{A}^{\varepsilon}_{m+M}(h_+,h^{\varepsilon}):=\big( P_- t_{\Sigma} \mathcal{E}^{\varepsilon}_{m+M}(z)P_{+},P^{{\varepsilon}}_- t_{\Sigma^{\varepsilon}} \mathcal{E}^{\varepsilon}_{m+M}(z)P^{{\varepsilon}}_{+}\big)$.\\
      
In particular, for $z\in \rho(D_m)$ we have the following explicit formulas
\begin{align*}
&\mathscr{A}_{m}(z)= -P_+ \beta (\beta/2 + \mathscr{C}_m^{z})^{-1}P_{-} ,\quad\mathscr{A}_{m}^{\varepsilon}(z)= -P_+^\varepsilon \beta (\beta/2 + \mathscr{C}_m^{z,\varepsilon})^{-1}P^\varepsilon_{-}.
\end{align*}
\end{definition}
\begin{remark}\label{remarkA}
    We define the Poincaré-Steklov operator, $\mathtt{A}_{m+M}^\varepsilon$, as a part of the operator $\mathcal{A}_{m+M}^\varepsilon$, which is only associated with $\Sigma^\varepsilon$ as follows:
     \begin{align*}\label{}
          \mathtt{A}^{\varepsilon}_{m+M}(z):  P^\varepsilon_+H^{1/2}(\S^{\varepsilon})^{4}&\to P^\varepsilon_-H^{1/2}(\S^{\varepsilon})^{4}\\ h^\varepsilon &\mapsto \mathtt{A}^{\varepsilon}_{m+M}(z) h^\varepsilon:= P^{{\varepsilon}}_- t_{\Sigma^{\varepsilon}} \mathcal{E}^{\varepsilon}_{m+M}(z)P^{{\varepsilon}}_{+}.
      \end{align*}
      In particular, $\mathtt{A}^{\varepsilon}_{m+M}$ will be used to establish the approximation in Section \ref{section4}.
\end{remark}
\subsection{Some geometric aspect}\label{SubNot}
\begin{definition}{[Weingarten map].}\label{Weingarten}
Let $\S$ be parametrized by the family $\lbrace \phi_{j}, U_{j},V_{j},\rbrace_{j\in J}$ with $J$ a finite set, $U_j\subset \rr^{2},\,V_j \subset \rr^{3},\,\S\subset \bigcup_{j\in J}V_j$ and $\phi(U_j)=V_j$ and $\phi_j=V_j\cap\S$  for all $j\in J.$ For $x=\phi_j (u)\in\S\cap V_j$ with $u\in U_j,$ one defines the Weingarten map (arising from the second fundamental form) as the following linear operator
\begin{align}
\begin{array}{rcl}
W_{x}:=W(x):T_x &\to & T_x\\
\partial_i \phi_j (u) &\mapsto & W(x)[\partial_i \phi_j] (u):=-\partial_i \nu (\phi_j (u)),
\end{array}
\end{align} where $T_x$ denotes the tangent space of $\S$ on $x$ and $\lbrace \partial_i \phi_j (u)\rbrace_{i=1,2}$ is a basis vector of $T_x$.\\

 The eigenvalues $k_{1}(x),....,k_{n}(x)$ of the Weingarten map $W_{x}$ are called principal curvatures of $\S$ at $x$. Then, we have the following proposition:
\end{definition}
\begin{proposition}{[\cite{JAT}, Chapter 9 (Theorem 2), 12 (Theorem 2)].}\label{PWM}
Let ${\S}$ be an $n-$surface in $\mathbb{R}^{n+1}$, oriented by the unit normal vector field $\nu$, and let $x\in {\S}$. The principal curvatures are uniformly bounded  on $\S$.
\end{proposition}
\begin{definition}{[Transformation operator].}\label{Transformationoperator}
Let $\S$, $\S^{\varepsilon}\subset \rr^3$ be as above. We define the diffeomorphism $p:\S\longrightarrow \S^{\varepsilon}$ such that for all $x_{\S}\in\S$, we get $p(x_{\S}):=x_{\S}+ \varepsilon \nu(x_{\S}),\,\varepsilon\in(0,\varepsilon_0).$ Then for $\varepsilon_0$ sufficiently small, we define the transformation operator as an unitary and invertible operator as follows
\begin{align}\label{transformation}
\begin{array}{rcl}
\mathcal{T}_{\varepsilon}:L^{2}(\S)^{4} &\to & L^{2}(\S^{\varepsilon})^{4},\\
\psi &\mapsto & \mathcal{T}_{\varepsilon}[\psi](x)=\dfrac{1}{\mathrm{det}(1-\varepsilon W(x_{\S}))}(\psi\circ p^{-1})(x),\quad x=p(x_\Sigma),
\end{array}
\end{align} 
and its inverse is given by 
\begin{align*}
\begin{array}{rcl}
   \mathcal{T}_{\varepsilon}^{-1}:L^{2}(\S^\varepsilon)^{4} &\to & L^{2}(\S)^{4},\\
\varphi &\mapsto &\mathcal{T}^{-1}_{\varepsilon}[\varphi](x_\Sigma)=\mathrm{det}(1-\varepsilon W(x_{\S}))(\varphi\circ p)(x_\Sigma). 
\end{array}
\end{align*}
\end{definition}
We also introduce the projection $P_{\S}:\mathcal{U}^{\varepsilon}\longrightarrow \S$ given by 
\begin{align*}
   P_{\S}(x_{\S}+t\nu(x_{\S})):=x_{\S},\quad\forall\,\,x_{\S}\in \S\,\,\text{and}\,\,t\in(0,\varepsilon].
\end{align*}
\section{Parametrix for the Poincaré-Steklov operators (large mass limit)}\label{section4}
Set $\kappa:= (M+m)$. This section is devoted to study the (classical and semiclassical) pseudodifferential properties of the Poincaré-Steklov operator, $\mathcal{A}_\kappa^\varepsilon$, in order to use it in the application of Section \ref{app}. The main goal of this section is to study the Poincar\'e-Steklov operator, $\mathcal{A}^{\varepsilon}_\kappa$,  as a  $\kappa$-dependent pseudodifferential operator when $\kappa$ is large enough.  Roughly speaking, we will look for a local approximate formula for the solution of \eqref{L2b}. The approximation in this section follows the steps of the one in paper \cite[Section 5]{BBZ}, but since our elliptic problem \eqref{L2b}, defined on the domain $\mathcal{U}^\varepsilon$, has two different boundary ($\partial\mathcal{U}^\varepsilon = \Sigma\cup\Sigma^\varepsilon$), and we have to take into account the dependence in $\varepsilon$, so we prefer to study rigorously the construction of the approximation. Once this is done, we use the regularization property of the resolvent of the MIT bag operator to catch the semiclassical principal symbol of $\mathcal{A}^{\varepsilon}_\kappa$. Throughout this section, we assume that $z\in\rho(D_{\mathrm{MIT}}^{\mathcal{U}^\varepsilon}(\kappa))$.\\

We see that $\mathcal{U}^{\varepsilon}$ has two boundaries, $\Sigma$ and $\Sigma^\varepsilon$. Since the approximation with respect to $\Sigma$ has already been established in \cite[Section 4]{BBZ}, and we therefore have this result in the present problem, it is then sufficient to establish the approximation of $\mathcal{A}^{\varepsilon}_\kappa$ just with respect to $\Sigma^\varepsilon$. For this purpose, and for simplicity of notation, we set $\mathcal{A}^h:=\mathtt{A}^\varepsilon_\kappa$ with $\varepsilon\equiv h:=\kappa^{-1} \in (0,1]$ as the semiclassical parameter, where $\mathtt{A}^\varepsilon_\kappa$ is defined in Remark \ref{remarkA}.
\subsection{Symbol classes  and Pseudodifferential operators}\label{section5} We recall here the basic facts concerning the classes of pseudodifferential operators that will serve in the rest of the paper. Let $\mathscr{M}_{4}(\cc)$ be the set of $4\times 4$ matrices over $\cc$.  For $d\in\mathbb{N}^{\ast}$ we let $\sS^{m}(\rr^d\times\rr^d)$ be the standard symbol class  of order $m\in\rr$ whose elements are matrix-valued functions $a$ in the space $\mathit{C}^{\infty}(\rr^d\times\rr^d;\mathscr{M}_{4}(\cc))$ such that 
	\begin{align*}
		|\partial^{\alpha}_{x}\partial^{\beta}_{\xi}a(x,\xi) |\leqslant C_{\alpha\beta}(1+|\xi|^2)^{m-|\beta|}, \quad \forall (x,\xi)\in\rr^d\times\rr^d, \,\,\forall \alpha\in\nn^{d},\,\,\forall \beta\in\nn^{d}.
	\end{align*}
	Let $\mathscr{S}(\rr^d)$ be the Schwarz class of functions. Then,  for each  $a\in  \sS^{m}(\rr^d\times\rr^d)$ and any  $h \in (0,1]$,  we associate a semiclassical pseudodifferential operator $Op^h(a):\mathscr{S}(\rr^d)^4\rightarrow \mathscr{S}(\rr^d)^4$ via the standard formula
	$$Op^h(a)u(x)=\frac{1}{(2\pi)^{d}}\int_{\rr^d}e^{i\xi\cdot x}a(x,h\xi) \hat{u} (\xi)\mathrm{d}\xi, \quad \forall u\in \mathscr{S}(\rr^d)^4.$$
	If $a\in \sS^0 (\rr^d\times\rr^d)$, then Calder\'on-Vaillancourt theorem's (see, \emph{e.g.}, \cite{CAVA}) yields that $Op^h(a)$ extends to a bounded operator from $\mathit{L}^{2}(\rr^d)^4$ into itself, and there exists $C, N_C>0$  such that 
	\begin{align}\label{Calderon-Vaillancourt}
		\left| \left| Op^{h}(a)\right|\right|_{\mathit{L}^{2}\rightarrow \mathit{L}^{2}}\leqslant C \max_{ |\alpha+\beta|\leqslant N_C}\left| \left| \partial_{x}^{\alpha}\partial_{\xi}^{\beta}a \right|\right|_{\mathit{L}^\infty}.
	\end{align}
	By definition, a semiclassical pseudodifferential operator $Op^{h}(a)$, with $a\in \sS^0 (\rr^d\times\rr^d)$, can also be considered as a classical pseudodifferential operator $Op^{1}(a_h)$ with $a_h=a(x,h\xi)$ which is bounded with respect to $h\in(0,h_0)$, where $h_0>0$ is fixed. Thus the Calder\'on-Vaillancourt theorem also provides the boundedness of these operators in Sobolev spaces $\mathit{H}^s(\rr^d)^4=\langle D_x\rangle^{-s}\mathit{L}^2(\rr^d)^4$ where $\langle D_x\rangle=\sqrt{-\Delta +\mathbb{I}}$. Indeed,  we have 
	\begin{align}\label{Calderon-VaillancourtHS}
		\left| \left| Op^{1}(a_h)\right|\right|_{\mathit{H}^{s}\rightarrow \mathit{H}^{s}}= \left| \left| \langle D_x\rangle^{s}Op^{1}(a_h)\langle D_x\rangle^{-s}\right|\right|_{\mathit{L}^{2}\rightarrow \mathit{L}^{2}},
	\end{align}
	and since $\langle D_x\rangle^{s}Op^{1}(a_h)\langle D_x\rangle^{-s}$ is a classical pseudodifferential operator with a uniformly bounded symbol in $\sS^0$, we deduce that $Op^{h}(a)$ is uniformly bounded with respect to $h$ from $\mathit{H}^{s}$ into itself.
\subsection{Reduction to local coordinates}
Let us consider $\mathbb{A}=\{ (U_{\varphi_j}, V_{\varphi_ j},\var_j): j\in\{1,\cdots, N\} \}$ an atlas of $\S$ and $(U_{\varphi}, V_\varphi,\var) \in \mathbb{A}$.  We consider also the case where $U_{\varphi}$ is the graph of a smooth function $\chi$, and we assume that $\O_-^\varepsilon$ corresponds locally to the side $x_3 > \chi(x_1,x_2)$. Then,  for
\begin{align*} 
U_{\varphi} =  & \{ (x^1_{\S},x^2_{\S}, \chi(x^1_{\S},x^2_{\S})) ; \, (x^1_{\S},x^2_{\S}) \in V_\varphi \}; \quad \varphi((x^1_{\S},x^2_{\S}, \chi(x^1_{\S},x^2_{\S})) = (x^1_{\S},x^2_{\S}), \\
\mathcal{V}_{\varphi,\eta}:= & \{ (y_1,y_2, y_3 + \chi(y_1,y_2)) ; \, (y_1,y_2,y_3) \in V_\varphi \times (0, \eta) \} \subset \O_+,
	\end{align*} 
with $\eta$ sufficiently small, we have the following homeomorphism:
\begin{align*} 
\phi :  \mathcal{V}_{\varphi,\eta} & \longrightarrow  V_\varphi \times (\varepsilon, \eta)\\
 (x_{\S}^1,x^2_{\S}, x_{\S}^{3})  & \mapsto (x^{1}_\S,x_{\S}^2, x^{3}_{\S} - \chi(x_{\S}^{1},x^2_{\S})),
	\end{align*} 
and the  pull-back 
\begin{align*} 
\phi^*:  C^\infty(V_\varphi \times (\varepsilon, \eta)) & \longrightarrow   C^\infty(\mathcal{V}_{\varphi,\eta} ) \\
v  & \mapsto \phi^*v:= v  \circ \phi .
	\end{align*} 
Now, using the coordinates in \eqref{conditionsgeo}, we let the diffeomorphism $\phi_{\varepsilon}: C^{\infty}(\mathcal{V}_{\varphi,\eta}) \longrightarrow (\mathcal{V}_{\varphi,\eta}^{\varepsilon})$ defined by follows: 
\begin{align*} 
\phi_{\varepsilon}(x_1,x_2,x_3) := \phi(x^1_{\S},x^2_{\S}, x^3_{\S}) + \varepsilon \nu(\phi(x_{\Sigma}))= \big(x^1_{\S}+\varepsilon\nu_{1},x^2_{\S}+\varepsilon\nu_{2}, x^3_{\S}+\varepsilon\nu_{3} - \chi(x^1_{\S},x^2_{\S})\big),
	\end{align*} 
with $\tilde{y}=(y_1 , y_2)$ and  $\nu$ the outward pointing normal to $\Omega_+$. Now, let $\nu^\varphi=(\varphi^{-1})^*\nu$ be the pull-back of the outward pointing normal to $\O_+$ restricted on $V_\varphi$:
\begin{equation*}
\nu^\varphi(\ty) =\frac{1}{\sqrt{1+ | \nabla \chi |^2}} \begin{pmatrix}
-\partial_{x_1} \chi  \\   -\partial_{x_2} \chi   \\  1 
\end{pmatrix} (y_1,y_2) =:
 \begin{pmatrix}
\nu^{\varphi}_1  \\   \nu^{\varphi}_2   \\  \nu^{\varphi}_3
\end{pmatrix}.
\end{equation*}
Then, the pull-back $(\phi_{\varepsilon}^{-1})^{*}$ transforms the differential operator $D_m$ restricted on $\mathcal{V}_{\varphi, \eta}$ into the following operator on $V_\varphi \times (0, \eta)$:
\begin{align*} 
\widetilde{D}^\varphi_m&:=  (\phi^{-1}_{\varepsilon} )^* D_m (\phi_{\varepsilon} )^* \\& = 
 -i\left( \alpha_1 \partial_{y_1} +   \alpha_2 \partial_{y_2} - ( -\alpha_1 \partial_{x_1} \chi  -\alpha_2 \partial_{x_2} \chi + \alpha_3)\partial_{y_3}\right) + m \beta -i\varepsilon\left[c_1\partial_{y_{1}}+c_2\partial_{y_{2}}+c_3\partial_{y_{3}} \right]\\
 & =  -i( \alpha_1 \partial_{y_1}  +   \alpha_2 \partial_{y_2}) + \sqrt{1+ | \nabla \chi |^2}(i \alpha \cdot \nu^\varphi)(\ty )\partial_{y_3}  -i\varepsilon\left[c_1\partial_{y_{1}}+c_2\partial_{y_{2}}+c_3\partial_{y_{3}} \right] + m \beta,
 \end{align*} 
 where $c_{\bullet}$ are $4\times4$ matrices having the form $c_\bullet=(\alpha_1\partial_{x_1}+\alpha_2\partial_{x_2})\nu^{\varphi}_\bullet,$ for $\bullet=1,2$,$3.$
\\\\
Thus, in the variable $y \in V_\varphi \times (\varepsilon, \eta)$ for $0<\varepsilon<\eta$, the system \eqref{L2b}  
becomes:
\begin{equation}
	\left\{
\begin{aligned}\label{T11}	
&(\widetilde{D}^\varphi_\kappa  -z )u	=0,  \quad  &\text{ in } & V_\varphi \times (\varepsilon, +\infty) ,\\
	&\G^\varphi_- u= g^\varphi= g \circ \varphi^{-1},  \quad    &\text{ on }& V_\varphi \times \{\varepsilon\}, 
	\end{aligned}
	\right.
	\end{equation}
where $ \G^\varphi_\pm = P^\varphi_{\pm}t_{\{ y_3 = \varepsilon \}}$. \\\\
By isolating the derivative with respect to $y_3$, and using that $(i\alpha \cdot \nu^\varphi)^{-1} = - i \alpha \cdot \nu^\varphi$, we get
\begin{align*}
    &\partial_{y_3} u =  \left(\mathbb{I}_{4} - \dfrac{\varepsilon(\alpha\cdot \nu^{\varphi}c_{3})}{\sqrt{1+| \nabla \chi|^{2}}} \right)^{-1} \dfrac{i \alpha \cdot \nu^\varphi (\ty)}{\sqrt{1+ | \nabla \chi (\ty) |^2}} 
	\Big( -i \alpha_1 \partial_{y_1}  - i \alpha_2 \partial_{y_2}   + m \beta -z - i\varepsilon c_1 \partial_{y_{1}} - i\varepsilon c_2\partial_{y_2} \Big) u .
 \end{align*}
Since, $\dfrac{(\alpha\cdot \nu^{\varphi}c_3)}{\sqrt{1+|\nabla \chi |^{2}}}$ is a bounded linear operator, then for $\varepsilon\in (0,\varepsilon_0)$ with $\varepsilon_0$ sufficiently small, the following Neumann series converges
\begin{align*}
    \left(\mathbb{I}_{4} - \dfrac{\varepsilon(\alpha\cdot \nu^{\varphi}c_{3})}{\sqrt{1+| \nabla \chi|^{2}}} \right)^{-1}=\sum_{k=0}^{+\infty} \varepsilon^{k}\left(\dfrac{\alpha\cdot\nu^{\varphi}c_3}{\sqrt{1+ | \nabla \chi (\ty) |^2}}\right)^{k},
\end{align*}
and we obtain
\begin{equation*}
	\left\{
\begin{aligned}\label{}	
&\partial_{y_3} u	=\\& \sum_{k=0}^{+\infty}\varepsilon^{k}\Bigg(\dfrac{\alpha\cdot\nu^{\varphi} c_{3}}{\sqrt{1+|\nabla\chi(\tilde{y})|^{2}}}\Bigg)^{k+1} 
	\Big( -i \alpha_1 \partial_{y_1}  - i \alpha_2 \partial_{y_2}   +  \kappa\beta  - i\varepsilon c_1 \partial_{y_{1}} - i\varepsilon c_{2} \partial_{y_{2}} - z \Big) u, \quad  &\text{in }& V_\varphi \times (\varepsilon, +\infty) ,\\
	&\G^\varphi_- u= g^\varphi ,  \quad    &\text{on }& V_\varphi \times \{\varepsilon\}.
	\end{aligned}
	\right.
	\end{equation*}
Let us now introduce the matrices-valued symbols 
\begin{equation}\label{defL0L1}
L_0(\ty, \xi) := \dfrac{i \alpha \cdot \nu^\varphi (\ty)}{\sqrt{1+ | \nabla \chi (\ty) |^2}} 
	\Big( \alpha \cdot \xi    +  \beta \Big), \quad \text{and} \quad L_1(\ty):= \dfrac{i \alpha \cdot \nu^\varphi (\ty)}{\sqrt{1+ | \nabla \chi (\ty) |^2}}\Big( c \cdot \xi    -   z \Big),
\end{equation}
with $\xi = (\xi_1, \xi_2)\in \mathbb{R}^{2}$ identified with $(\xi_1, \xi_2,0)\in\mathbb{R}^{3}$ and $c=(c_1,c_2)$.
Then for $\varepsilon=h:=1/m$, the system \eqref{T11} becomes:
\begin{equation}
	\left\{
	\begin{aligned}\label{T1bh}	
&h \partial_{y_3} u^{h}	= L_0(\ty, h D_{\ty})  u^{h} + h L_1(\ty, h D_{\ty})  u^{h}  \\& +\sum_{k=1}^{+\infty}h^{k}\dfrac{(\alpha\cdot\nu^{\varphi}c_3)^{k}}{(1+|\nabla \chi |^{2})^{k/2}}\left ( L_0(\ty, h D_{\ty})  u^{h} + h L_1(\ty, h D_{\ty}) \right)  u^{h} 
 , \,  &\text{ in } &\mathbb{R}^{2} \times (\varepsilon, +\infty) ,\\ 
 \vspace{0.5cm}
	&P^{\varphi}_+ t_{\lbrace y_{3}=\varepsilon \rbrace}u^{h} = g^\varphi ,  \,    &\text{ on } &\rr^{2} \times \{\varepsilon\}.
	\end{aligned}
	\right.
	\end{equation}
 \begin{remark}\label{Remark3.1}
     In this remark, we clarify the first difference in the approximation of this section compared to that of \cite[Section 5]{BBZ}. Indeed, according to the formula of $L_1$ from \eqref{defL0L1}, we observe that the term $c\cdot\xi$ appears in our case, whereas it was absent in the case of \cite{BBZ}. Moreover, we mention that this difference plays an important role in the subsequent progression of this approximation, exerting a significant impact on the symbol class of the solution $u^h$.
 \end{remark}
Before constructing an approximate solution of the system \eqref{T1bh}, let us give some properties of $L_0$. Besides, we mention that $L_1$ also verifies these properties.
\begin{lemma}\label{properties3.1} 
 Recall the projections  $P^\varphi_\pm := (\mathbb{I}_4\mp i\beta \, \alpha \cdot \nu^\varphi(\ty))/2$, and set
 \begin{align}\label{spinangular}
 \gamma_5 :=-i\alpha_1\alpha_2\alpha_3=\begin{pmatrix}
0& \mathbb{I}_2 \\
\mathbb{I}_2  & 0
\end{pmatrix}\quad\text{and}\quad S\cdot X= -\gamma_5 (\alpha\cdot X), \quad \forall X\in\rr^3.
\end{align}
Using the anticommutation relations of the Dirac's matrices we easily get the following identities
	\begin{align*}\label{}
	\begin{split}
	i (\alpha \cdot X ) (\alpha \cdot Y  ) &= i X \cdot Y + S \cdot (X \wedge Y),\\
	\left\{S \cdot X,   \alpha \cdot Y\right\} &= -(X \cdot Y) \g_5,\quad [S \cdot X,\beta]=0, \quad \forall X,Y\in\rr^3.
	\end{split}
	\end{align*}
Let  $\nu^\varphi$ and $\xi$ be as above. Then, for any  $z \in \cc$ and any $\tau\in\rr^3$ such that  $\tau \perp \nu^\varphi$, the following identities hold:
	\begin{align*}
	 \left(S \cdot \tau   - i m \beta (\alpha \cdot \nu^\varphi (\ty))  \right)^2 = \left(|\tau|^2  + m^2  \right) \mathbb{I}_4,
	 \end{align*}
	 \begin{align*}
 P_\pm^{\varphi} (S\cdot\tau)=(S\cdot\tau)P^{\varphi}_\mp\quad\text{ and }\quad P^{\varphi}_\pm(i\alpha\cdot \nu^{\varphi})= (i\alpha\cdot \nu^{\varphi})P^{\varphi}_\mp .
	  \end{align*}
		\end{lemma}
		The next proposition gathers the main properties of the operator $L_0$.
	\begin{proposition}{\cite[Proposition 5.1]{BBZ}.}\label{expo de q} Let $L_0(\ty,\xi)$ be as in \eqref{defL0L1}, then we have 
        \begin{align*}
	L_0(\ty,\xi) &= \frac{1 }{\sqrt{1+ | \nabla \chi (\ty)|^2}} 
	\Big( i \xi \cdot \nu^\varphi(\ty)  + S \cdot (\nu^\varphi(\ty)  \wedge \xi)   - i  \beta (\alpha \cdot \nu^\varphi(\ty) )    \Big)\\
	 &=i \xi \cdot \tilde{\nu}^\varphi (\ty) + \frac{\lambda (\ty,\xi) }{\sqrt{1+ | \nabla \chi (\ty)|^2}}\Pi_+ (\ty,\xi) - \frac{\lambda (\ty,\xi) }{\sqrt{1+ | \nabla \chi (\ty)|^2}}\Pi_-  (\ty,\xi),
	\end{align*}
	where 
	\begin{align}\label{deflambdapro}
	\begin{split}
	\lambda(\ty, \xi):=&\sqrt{ |{\nu}^\varphi(\ty) \wedge\xi|^2  + 1 },\\
	\tilde{\nu}^{\varphi}(\ty) :=& \frac{1 }{\sqrt{1+ | \nabla \chi |^2}} \nu^{\varphi}(\ty),\\
	 \Pi_\pm  ( \ty,\xi) :=& \frac{1}2 \left(  \mathbb{I}_4 \pm    \frac{S \cdot (\nu^\varphi(\ty) \wedge \xi)   - i \beta (\alpha \cdot \nu^\varphi (\ty)) }{\lambda(\ty, \xi)}  \right).
	 \end{split}
	\end{align}	
 In particular, the symbol $L_0(\ty, \xi)$ is elliptic in symbol class ${\sS}^1$ (defined in Section \ref{section5}) and it admits two eigenvalues $\varrho_\pm (\cdot, \cdot) \in {\sS}^1$ of multiplicity $2$ which are given by 
	\begin{align*}\label{}
	\varrho_{\pm}(\ty,\xi)= \frac{i   \nu^\varphi(\ty)\cdot \xi \pm  \lambda(\ty,\xi) }{\sqrt{1+ | \nabla \chi |^2}} ,
	\end{align*}
	and for which there exists $c>0$ such that 
	\begin{equation}\label{ellip}
\pm \Re \varrho_\pm (\ty, \xi) > c \langle \xi \rangle,
\end{equation}
uniformly with respect to $\ty$. Moreover, $\Pi_\pm(\ty , \xi)$ are the projections onto $\mathrm{Kr}(L_0(\ty,\xi) -  \varrho_\pm (\ty, \xi)  \mathbb{I}_4 )$, belong to the symbol class ${\sS}^0$ and satisfy:
\begin{align}\label{PPiP}
\begin{split}
P^\varphi_\pm \,  \Pi_\pm ( \ty,\xi)\,  P^\varphi_\pm = k_{+}^\varphi ( \ty,\xi)  P^\varphi_\pm \quad\text{and }\, P^\varphi_\pm \,  \Pi_\mp ( \ty,\xi)\,  P^\varphi_\mp=\mp\Theta^{\varphi}(\ty, \xi)P^\varphi_\mp,
\end{split}
\end{align}
	with 
\begin{align*} \label{}
\begin{split}
k_\pm^\varphi ( \ty,\xi) = \frac{1}{2}\left(1\pm \frac{1}{\lambda(\ty, \xi)}\right)  ,\quad \Theta^{\varphi}(\ty, \xi)=\frac{1}{2\lambda(\ty, \xi)}\left(S\cdot (\nu^\varphi(\ty)\wedge\xi) \right).
\end{split}
\end{align*}
Now,  using Lemma \ref{properties3.1} and  the properties  \eqref{deflambdapro},  a simple computation shows that
\begin{align*}
	P^{\varphi}_+   \Pi_{\pm} &= k_\pm^\varphi P^{\varphi}_+ \pm\frac{1}{2\lambda}\left(S\cdot (\nu^\varphi(\ty)\wedge\xi) \right) \,  P^{\varphi}_-,\\
	P^{\varphi}_-  \Pi_\pm&= k_\mp^\varphi P^{\varphi}_- \pm\frac{1}{2\lambda}\left(S\cdot (\nu^\varphi(\ty)\wedge\xi)\right) \,  P^{\varphi}_+.
\end{align*}
That is, $k^{\varphi}_+$ is a positive function of ${\sS}^0$,  $(k^{\varphi}_+)^{-1} \in {\sS}^0$  and  $\Theta^{\varphi}\in {\sS}^0$ where $\mathcal{S}^0$ is zero-order symbol class defined in Section \ref{section5}. 
\end{proposition}	
\subsection{Semiclassical parametrix for the boundary problem}\label{ssSclBPb}
In this section,  we construct the approximate solution of the system \eqref{T1bh}. For simplicity of notation, in the sequel  we will use $y$, $\tau$, and $P_{\pm}$ instead of $\ty$, $y_3$, and $P_{\pm}^{\varphi}$, respectively. We  are going to construct a local approximate solution of the following first order system:
\begin{equation*}
	\left\{
	\begin{aligned}\label{}	
&h \partial_{y_3} u^{h}	= L_0(y, h D_{y})  u^{h} + h L_1(y, h D_{y})  u^{h}  \\& +\sum_{k=1}^{+\infty}h^{k}\dfrac{(\alpha\cdot\nu^{\varphi}c_3)^{k}}{(1+|\nabla \chi |^{2})^{k/2}}\left ( L_0(y, h D_{y})  u^{h} + h L_1(y, h D_{y}) \right)  u^{h} 
 , \,  &\text{ in } &\mathbb{R}^{2} \times (\varepsilon, +\infty) ,\\\\
	& P_+ t_{\lbrace \tau=\varepsilon \rbrace}u^{h} = g^\varphi ,  \,    &\text{ on } &\rr^{2} \times \{\varepsilon\}.
	\end{aligned}
	\right.
	\end{equation*}
 This system is equivalent to 
 \begin{equation}
	\left\{
	\begin{aligned}\label{T1bh'}	
&h \partial_{y_3} u^{h}	= L_0(y, h D_{y})  u^{h}   +\sum_{k=1}^{+\infty}h^{k}\dfrac{(\alpha\cdot\nu^{\varphi}c_3)^{k-1}}{(1+|\nabla \chi |^{2})^{\frac{k-1}{2}}} \widetilde{L}_{1}(y,h D_{y})  u^{h} 
 , \,  &\text{ in } &\mathbb{R}^{2} \times (\varepsilon, +\infty) ,\\
	 &P_+ t_{\lbrace \tau  = \varepsilon \rbrace}u^{h} = g^\varphi ,  \,    &\text{ on } &\rr^{2} \times \{\varepsilon\},
	\end{aligned}
	\right.
	\end{equation}
with $\widetilde{L}_{1}(y,\xi)=L_{1}(y,\xi)+(\alpha\cdot\tilde{\nu}^{\varphi}c_{3})L_{0}(y,\xi).$
 \vspace{0.5cm}\\
To be precise, we will look for a solution $u^h$  in the following form:
\begin{equation*}\label{}
u^h(y,\tau)= 
\mathit{Op}^h (A^h(\cdot, \cdot, \tau) ) f = 
 \int_{\rr^2} A^h(y, h \xi , \tau) e^{i y \cdot \xi} \hat{f}(\xi) \text{d} \xi ,
\end{equation*}
with $A^h(\cdot , \cdot , \tau) \in {\sS}^0$ for any $\tau>0$  constructed inductively in the form:
$$A^h(y,  \xi , \tau) \sim \sum_{j \geq 0} h^j A_j(y, \xi , \tau).$$
The action of 
    $h \partial_{y_3} 	- L_0 - \sum_{k=1}^{+\infty}h^{k}\dfrac{(\alpha\cdot\nu^{\varphi}c_3)^{k-1}}{(1+|\nabla \chi |^{2})^{\frac{k-1}{2}}} \widetilde{L}_{1}$ on $A^h(y, h D_y, \tau) f$ is given by $T^h(y, h D_y, \tau) f$, with
\begin{align*}
 T^h(y, \xi , \tau) &=   h (\partial_\tau A^{h}) (y, \xi , \tau) - L_0(y,\xi) A^{h} (y, \xi , \tau)  - h \Big(  \widetilde{L}_1(y,\xi) A^{h} (y, \xi , \tau) - i\partial_\xi L_0(y, \xi) \cdot \partial_y  A^{h} (y, \xi , \tau) \Big) \\& - h^{2}\Big( L_{0}A^{h}+\widetilde{L}_{1}(y,\xi)A^{h} + \partial_{\xi} L_{0}\cdot\partial_{y} A^{h} - i\partial_{\xi}\widetilde{L}_{1}\cdot\partial_{y}A^{h}+ (\alpha\cdot\tilde{\nu}^{\varphi}c_{3})\widetilde{L}_{1}(y,\xi)A^{h} \Big) + ....
\end{align*}
Then, by identifications of the coefficients of $j$, $j\geqslant 0,$ we look for $A_0$ satisfying:
\begin{equation}\label{SystS0}
	\left\{
	\begin{aligned}
	h \partial_{\tau} A_0 (y, \xi , \tau)  & =   L_0(y, \xi) A_0 (y, \xi , \tau) , \\
	P_+ (y) A_0 (y, \xi , \varepsilon)   & =   P_+ (y) ,
	\end{aligned}
	\right.
	\end{equation}
and for $j \geq 1$,
\begin{equation}\label{SystSj}
	\left\{
	\begin{aligned}
 & h \partial_{\tau} A_j (y, \xi , \tau)   =   L_0(y, \xi) A_j (y, \xi , \tau) +  \Big(\widetilde{L}_{1} (y,\xi)-i\partial_{\xi}L_{0}(y,\xi)\cdot\partial_{y}\Big) A_{j-1} (y, \xi , \tau) \\& + \sum_{l\geq2}^{l=j}(\alpha\cdot\tilde{\nu}^{\varphi}c_{3})^{j-l}\Big ((\alpha\cdot\tilde{\nu}^{\varphi}c_{3})\widetilde{L}_{1}(y,\xi)-i\partial_{\xi}\widetilde{L}_{1}(y,\xi)\cdot\partial_{y}\Big)A_{l-2}(y,\xi,\tau), \\\\
	& P_+ (y) A_j (y, \xi , \varepsilon)    =  0 .
	\end{aligned}
	\right.
	\end{equation}
 	Let us introduce a class of parametrized symbols, in which we will construct the family $A_j $:
\begin{equation*}
{\cP}_{h}^{m}:= \{ b(\cdot, \cdot, \tau ) \in {\sS}^{m} ; \; \forall (k,l) \in \nn ^2 , \, \tau^k \partial_\tau^l b(\cdot, \cdot, \tau ) \in h^{k-l} {\sS}^{m-k+l}  \}; \quad m \in \zz.
\end{equation*} 	
 \begin{proposition}\label{SolSystS0}
There exists $A_0 \in {\cP}_{h}^{0}$ solution of \eqref{SystS0} given by:
\begin{align*}
A_0 (y, \xi , \tau) &=   e^{h^{-1} (\tau-\varepsilon) \varrho_-(y,\xi)}\quad \frac{\Pi_-(y,\xi) P_+(y)A_{0}(y,\xi,\varepsilon)}{k^{\varphi}_-(y,\xi)}\\&= e^{h^{-1} (\tau-\varepsilon) \varrho_-(y,\xi)}\dfrac{\Pi_-(y,\xi) P_+(y)}{k^{\varphi}_-(y,\xi)}\\& = e^{h^{-1} (\tau-\varepsilon) \varrho_-(y,\xi)}\left( \mathbb{I}_{4} - \frac{\Theta^{\varphi}}{k_{-}^{\varphi}}\right)P_+.
\end{align*}
\end{proposition}
\textbf{Proof.} 
The proof follows the same argument as \cite[Proposition 5.2]{BBZ}. The solution of the differential system $h \partial_{\tau} A_0  =   L_0  A_0$ is $A_0 (y, \xi , \tau) =   e^{h^{-1} (\tau-\varepsilon) L_0} A_0 (y, \xi , \varepsilon)$. By definition of $\varrho_\pm$ and $\Pi_\pm$, we have:
\begin{equation}\label{expL0}
e^{h^{-1} \tau L_0}  = e^{h^{-1} (\tau-\varepsilon) \varrho_-} \Pi_-(y,\xi) + e^{h^{-1} (\tau-\varepsilon) \varrho_+} \Pi_+(y,\xi).
\end{equation} 
It follows from \eqref{ellip} that $A_0 $ belongs to ${\sS}^0$  for any $\tau >\varepsilon$ if and only if $ \Pi_+(y,\xi) A_0 (y, \xi , \varepsilon) =0$. Moreover, the boundary condition $P_+ A_0  =   P_+ $ implies $P_+ (y) A_0 (y, \xi , \varepsilon)  =   P_+ (y)$.  Thus, we deduce that 
$$A_0 (y, \xi , \varepsilon) =  P_+(y) - \frac{P_- \Pi_+ P_+}{k^{\varphi}_-} (y,\xi) = P_+(y) + \frac{P_- \Pi_- P_+}{k^{\varphi}_-} (y,\xi) = \frac{\Pi_- P_+}{k^{\varphi}_-} (y,\xi).$$
The properties of $\varrho_-$, $\Pi_-$, $P_-$ and $k_+$ given in Proposition \ref{expo de q}, imply that $(k^{\varphi}_+)^{-1}\Pi_- P_- \in \sS^0$ and that $e^{h^{-1} \tau \varrho_-(y,\xi)} \in \cP^0_h$. This concludes the proof of Proposition \ref{SolSystS0}.\qed
\begin{proposition}\label{SolSystSj}
Let $A_0$ be defined by Proposition \ref{SolSystS0}. Then for any $j\geq 1$, there exists $A_j$ solution of \eqref{SystSj} which has the form:
\begin{equation}\label{eqAj}
 A_j (y, \xi , \tau) =  e^{h^{-1} (\tau - \varepsilon) \varrho_-(y,\xi)} \, \sum_{k=0}^{2j} (h^{-1} (\tau - \varepsilon)  \langle \xi \rangle   )^k B_{j,k} (y, \xi),\quad \text{with} \quad B_{j,k} \, \in \,h\, \sS^{0}.
 \end{equation} 	
\end{proposition}
\begin{remark}\label{Remark3.3}
An important difference in the approximation between the solution $A_j$ resulting from this work and the solution presented in the work \cite[Proposition 5.3]{BBZ} lies in the order of the standard symbol class $\mathcal{S}^m$. Indeed, by referring to the form of $A_2$ (see \eqref{A2} from Appendix \ref{Appendix}) one can deduce that the optimal order of the term $\Pi_- a_0 \Big(P_+ - \dfrac{P_+ \Theta^\varphi}{k^\varphi_-} + \Pi_+ a_0 \Big)$ in $B_{2,0}$ is in $h\,\mathcal{S}^0$, and this property is reflected in the construction of $A_j$ for $j \geq 3$. However, in \cite[Proposition 5.3]{BBZ}, it was possible to obtain all $A_j$ in $h^j\mathcal{S}^{-j}$. This discrepancy leads us to deduce the following propositions concerning the solutions $A_j$.
\end{remark}
\begin{remark}
We mention that this difference in the symbol class of terms $B_{j,k}$ with that obtained in \cite{BBZ} is mainly due to the difference discussed in Remark \ref{Remark3.1}, \emph{i.e.}, to the influence of $c\cdot\xi$ as presented in the formula of $L_1$ in system \eqref{T1bh}, and subsequently to that mentioned in Remark \ref{Remark3.3}.
\end{remark}
\textbf{Proof of Proposition \ref{SolSystSj}.} 
For initialization and calculation of $A_1$ and $A_2$, see Appendix \ref{Appendix}. So, for  $A_j$ with $j\geq 1$, it is sufficient to prove the induction step. Thus,  assume that the $A_j$ solution of \eqref{SystSj} satisfies the above property and let us prove that the same holds for $A_{j+1}$. In order to be a solution to the differential system 
\begin{align*}
    h \partial_{\tau} A_{j+1} (y, \xi , \tau)  & =   L_0(y, \xi) A_{j+1} (y, \xi , \tau) +  \Big(\widetilde{L}_{1} (y,\xi) - i\partial_{\xi}L_{0}(y,\xi)\cdot\partial_{y}\Big) A_{j} (y, \xi , \tau) \\& + \sum_{l=2}^{l=j+1}(\alpha\cdot\tilde{\nu}^{\varphi}c_{3})^{j+1-l}\Big ((\alpha\cdot\tilde{\nu}^{\varphi}c_{3})\widetilde{L}_{1}(y,\xi) - i\partial_{\xi}\widetilde{L}_{1}(y,\xi)\cdot\partial_{y}\Big)A_{l-2}(y,\xi,\tau),
\end{align*}
then, for $A_{j+1}$ we have: 
 \begin{equation}\label{Agene}
 \begin{aligned}
 A_{j+1}  &=   e^{h^{-1}L_0 (\tau -\varepsilon)} {A_{j+1}}_{| \tau=\varepsilon}   + e^{h^{-1} \tau L_0} \int_{\varepsilon}^\tau e^{-h^{-1} s L_0} \underbrace{\Big(\widetilde{L}_{1}  - i\partial_{\xi}L_{0}\cdot\partial_{y}\Big) A_{j} (y, \xi , \tau)}_{(a)}\mathrm{d}s \\& + e^{h^{-1} \tau L_0}\int_{\varepsilon}^\tau  e^{-h^{-1} s L_0}\,\,\underbrace{\sum_{l=2}^{l=j+1}(\alpha\cdot\tilde{\nu}^{\varphi}c_{3})^{j+1-l}\Big ((\alpha\cdot\tilde{\nu}^{\varphi}c_{3})\widetilde{L}_{1} -i \partial_{\xi}\widetilde{L}_{1}(y,\xi)\cdot\partial_{y}\Big)A_{l-2}(y,\xi,\tau))}_{(b)}\mathrm{d}s\\&
 := e^{h^{-1}L_0 (\tau -\varepsilon)} {A_{j+1}}_{| \tau=\varepsilon} + e^{h^{-1} \tau L_0}\int_{\varepsilon}^\tau e^{h^{-1}sL_0 }\boldsymbol{\Big ( (a)+(b)\Big)}\text{d}s.
 \end{aligned}
 \end{equation} 	
In order to  know the form of $\boldsymbol{(a)}$ and $\boldsymbol{(b)}$, let us consider the formula \eqref{TL2}. Then for the quantity $\boldsymbol{(a)}$, we have 
 $$  \partial_y  A_{j}  =  e^{h^{-1} (\tau-\varepsilon)  \varrho_-} \, \Big( h^{-1} (\tau - \varepsilon) \partial_y \varrho_- + \partial_y \Big) \sum_{k=0}^{2j} \Big(h^{-1} (\tau - \varepsilon) \langle \xi \rangle   \Big)^k B_{j,k}. $$
 Now, applying $\big(\widetilde{L}_1 - i\partial_{\xi}L_{0}\cdot\partial_{y}\big)$ to $A_j(y,\xi,\tau)$:
 \begin{align*}
\big(\widetilde{L}_1 - i\partial_{\xi}L_{0}\cdot\partial_{y}\big)A_{j}&=a_{0}(y)\big(-z+c\cdot\xi-ic_{3}L_{0} - i\alpha\cdot\partial_{y}\big)e^{h^{-1} (\tau-\varepsilon)  \varrho_-}\sum_{k=0}^{2j} \Big(h^{-1} (\tau  - \varepsilon)  \langle \xi \rangle   \Big)^k B_{j,k}\\&:=\underbrace{e^{h^{-1} (\tau-\varepsilon)  \varrho_-}a_{0}(y)\Big( -z+ c_{3}\alpha\cdot\tilde{\nu}^{\varphi}\beta - i\alpha\cdot\partial_{y}\Big)\sum_{k=0}^{2j} \Big(h^{-1} (\tau  - \varepsilon)  \langle \xi \rangle   \Big)^k B_{j,k}}_{\mathrm{(i)}}\\& + \underbrace{e^{h^{-1} (\tau-\varepsilon)  \varrho_-}a_{0}(y)\Big(c+c_{3}\alpha\cdot\tilde{\nu}^{\varphi}\alpha\Big)\cdot\xi\sum_{k=0}^{2j} \Big(h^{-1} (\tau  - \varepsilon)  \langle \xi \rangle   \Big)^k B_{j,k}}_{\mathrm{(ii)}}\\& + \underbrace{e^{h^{-1} (\tau-\varepsilon)  \varrho_-}a_{0}(y)\Big(-ih^{-1}(\tau - \varepsilon) \alpha\cdot\partial_{y}\varrho_-\Big)\sum_{k=0}^{2j} \Big(h^{-1} (\tau  - \varepsilon)  \langle \xi \rangle   \Big)^k B_{j,k}}_{\mathrm{(iii)}} .
\end{align*} 
Thanks to the properties of $ \varrho_-$ and $B_{j,k}$, (i), (ii) and (iii) have respectively the form:
 \begin{equation}\label{Agene1}
\mathrm{(i)}=\quad  e^{h^{-1} (\tau - \varepsilon) \varrho_-} \, \sum_{k=0}^{2j} (h^{-1} (\tau - \varepsilon) \langle \xi \rangle   )^k B^{'}_{j,k} (y, \xi),\quad \quad  
\end{equation} 
\begin{equation}\label{Agene2}
\mathrm{(ii)}=\quad  e^{h^{-1} (\tau - \varepsilon) \varrho_-} \, \sum_{k=0}^{2j} (h^{-1} (\tau - \varepsilon) \langle \xi \rangle   )^k \,\,\,\langle \xi \rangle \,\,\overline{B}_{j,k} (y, \xi) ,
\end{equation} 
\begin{equation}\label{Agene3}
\mathrm{(iii)}=\quad  e^{h^{-1} (\tau - \varepsilon) \varrho_-} \, \sum_{k=0}^{2j} (h^{-1} (\tau - \varepsilon) \langle \xi \rangle   )^{k+1} B^{''}_{j,k} (y, \xi) ,
\end{equation} 	
with $B^{'}_{j,k}$ and  $B^{''}_{j,k}$ verifying 
the properties of $B_{j,k}$, and $\langle \xi \rangle \,\,\overline{B}_{j,k}\in h\,\mathcal{S}^{1}$. Therefore, toghether \eqref{Agene1}, \eqref{Agene2} and \eqref{Agene3} give that 
 \begin{equation}\label{(a)}
\boldsymbol{\mathrm{(a)}}=\quad  e^{h^{-1} (\tau - \varepsilon) \varrho_-} \, \sum_{k=0}^{2j+1} (h^{-1} (\tau - \varepsilon) \langle \xi \rangle   )^{k} \widetilde{B}_{j,k} (y, \xi) ,
  \end{equation} 	
  where $\widetilde{B}_{j,k}$ verifies
\begin{align*}\label{}
\widetilde{B}_{j,k} \, \in \,h \,\sS^{1} \text{ for } k=0,...,2j, \text{ and } \widetilde{B}_{j,2j+1}\in h\,\mathcal{S}^{0}.
	\end{align*}
Similarly, to calculate $\boldsymbol{(b)}$, applying $\big(-i\partial_{\xi}\widetilde{L}_1\cdot\partial_{y}+(\alpha\cdot\tilde{\nu}^{\varphi}c_{3})\widetilde{L}_{1}\big)$ (see \eqref{TL2}) to the identity \eqref{eqAj} yields that 
\begin{equation}\label{TL12}
\begin{aligned}
   \big(-i\partial_{\xi}\widetilde{L}_1\cdot\partial_{y} & +(\alpha\cdot\tilde{\nu}^{\varphi}c_{3})\widetilde{L}_{1}\big)A_j =\\& e^{h^{-1} (\tau - \varepsilon) \varrho_-}  a_{0}(y)\Bigg (  d + e\cdot \xi -i h^{-1}(y_{3}-\varepsilon) f \cdot\partial_{y}\varrho_- \Bigg) \sum_{k=0}^{2j} \Big(h^{-1} (\tau  - \varepsilon)  \langle \xi \rangle   \Big)^k B_{j,k},
\end{aligned}
\end{equation}
with $d$, $e$ and $f$ defined in \eqref{Tdf}.
Let us decompose $\boldsymbol{(b)}$ as the following
\begin{align*}
    \sum_{l=2}^{l=j+1}(\alpha\cdot\tilde{\nu}^{\varphi}c_{3})^{j+1-l}&\Big (-i\partial_{\xi}\widetilde{L}_{1}\cdot\partial_{y}+(\alpha\cdot\tilde{\nu}^{\varphi}c_{3})\widetilde{L}_{1}\Big)A_{l-2}(y,\xi,\tau)):=\\& \underbrace{(\alpha\cdot\tilde{\nu}^{\varphi}c_{3})^{j-1}\Big (-i\partial_{\xi}\widetilde{L}_{1}\cdot\partial_{y}+(\alpha\cdot\tilde{\nu}^{\varphi}c_{3})\widetilde{L}_{1}\Big)A_{0}(y,\xi))}_{\mathrm{(m1)}}\\& + \underbrace{\sum_{l\geq 3}^{l=j+1}(\alpha\cdot\tilde{\nu}^{\varphi}c_{3})^{j+1-l}\Big (-i\partial_{\xi}\widetilde{L}_{1}\cdot\partial_{y}+(\alpha\cdot\tilde{\nu}^{\varphi}c_{3})\widetilde{L}_{1}\Big)A_{l-2}(y,\xi,\tau))}_{\mathrm{(m2)}}.
\end{align*}
 Since $A_{0}\in \mathcal{S}^{0},$ this gives that 
 \begin{align}\label{m1}
     \mathrm{(m1)}\quad =    \xi \cdot\dot{B}_{0,0} + \widehat{B}_{0,0} + (-ih^{-1} ( \tau - \varepsilon)f\cdot\partial_{y}\varrho_-)B_{0,0},
 \end{align}
 where $\dot{B}_{0,0}$, $\widehat{B}_{0,0}$  $\in\mathcal{S}^{0}$ are respectively the constants obtained by applying $d$ and $e$ to $\dfrac{\Pi_- P_+}{k_-^{\varphi}}$ and $f\cdot\partial_{y}\varrho_- \in \mathcal{S}^{1}$. Thus, $\mathrm{(m1)} \in \mathcal{S}^{1}$, $\forall j\geq 1$. \\
In the other hand, and for all $l\geq 3$ (\emph{i.e.}, $l-2 \geq 1)$, $A_{l-2}$ has the form \begin{equation}\label{eqAl2}
 A_{l-2} (y, \xi , \tau) =  e^{h^{-1} (\tau - \varepsilon) \varrho_-} \, \sum_{k=0}^{2(l-2)} (h^{-1} (\tau - \varepsilon)  \langle \xi \rangle   )^k B_{l-2,k} (y, \xi),
 \end{equation} 	 with $B_{l-2,k}\in h\,\mathcal{S}^{0}$. 
Applying \eqref{TL12} to the identity \eqref{eqAl2} we get 
\begin{align}\label{m2}
\mathrm{(m2)} \quad = e^{h^{-1} ( \tau - \varepsilon) \varrho_-}\sum_{l\geq 3}^{l=j+1} \sum_{k=0}^{2(l-2)+1} (h^{-1} (\tau - \varepsilon) \langle \xi \rangle   )^k \ddot{B}_{j,k} (y, \xi), 
\end{align}
with $\ddot{B}_{j,k}\in h\,\mathcal{S}^{1}$ and $\ddot{B}_{j,2(l-2)+1}\in h\,\mathcal{S}^{0}$. Therefore, for $i =(l-2)\geq 1$ and $j\geq 2,$ toghether \eqref{m1}, \eqref{m2} with \eqref{(a)} give that 
 \begin{equation}\label{(ab)}
 \begin{aligned}
\boldsymbol{\mathrm{(a)+(b)}}&=  e^{h^{-1} (\tau - \varepsilon) \varrho_-} \, \Bigg( \sum_{k=0}^{2j+1} (h^{-1} (\tau - \varepsilon) \langle \xi \rangle   )^{k} \widetilde{B}_{j,k}   + \sum_{l\geq 3}^{l=j+1} \sum_{k=0}^{2(l-2)+1} (h^{-1} (\tau - \varepsilon) \langle \xi \rangle   )^k \ddot{B}_{l-2,k}  + \mathrm{m1}\Bigg)\\&= e^{h^{-1} (\tau - \varepsilon) \varrho_-} \, \Bigg( \sum_{k=0}^{2j+1} (h^{-1} (\tau - \varepsilon) \langle \xi \rangle   )^{k} \widetilde{B}_{j,k}  + \sum_{i\geq 1}^{i=j-1} \sum_{k=0}^{2i+1} (h^{-1} (\tau - \varepsilon) \langle \xi \rangle   )^k \ddot{B}_{i,k}  + \mathrm{m1}\Bigg)\quad \\& = e^{h^{-1} (\tau - \varepsilon) \varrho_-} \, \Bigg( \sum_{k=0}^{2i+1} (h^{-1} (\tau - \varepsilon) \langle \xi \rangle   )^{k} \underbrace{\Bigg(\widetilde{B}_{i,k}  + \sum_{i\geq 1}^{i=j-1}\ddot{B}_{i,k} \Bigg)}_{C_{i,j,k}} + \,\mathrm{m1}\Bigg),
 \end{aligned}
 \end{equation} 	
with  ${C}_{i,j,k} \, \in \,h \,\sS^{1}$, and   ${C}_{i,j,k}\in h \,\sS^{0}$ $\mathrm{for}\quad k=2i+1 .$\\
 So, using the decomposition \eqref{expL0}, for the second term of the r.h.s. of \eqref{Agene} we have:
 \begin{equation}\label{2ndAgene}
 e^{h^{-1} \tau L_0} \int_\varepsilon^\tau e^{-h^{-1} s L_0} \Big(\boldsymbol{(a)} + \boldsymbol{(b)} \Big) \mathrm{d}s
 = e^{h^{-1} \tau \varrho_-}  \Pi_-  \mathbb{I}_-^j (\tau) 
 +
  e^{h^{-1} \tau \varrho_+}  \Pi_+  \mathbb{I}_+^j (\tau) ,
 \end{equation} 	
with  
\begin{align*}
    \mathbb{I}_\pm^j (\tau) = e^{-h^{-1}\varepsilon\varrho_-}\int_\varepsilon^\tau  e^{h^{-1} s( \varrho_- - \varrho_\pm)} \Bigg(\sum_{k=0}^{2i+1} (h^{-1} (s - \varepsilon)  \langle \xi \rangle   )^k C_{i,j,k}\quad +\quad \mathrm{m1}\Bigg)\, \mathrm{d}s.
\end{align*}
For $\mathbb{I}_-^j $, the exponential term is equal to $1$ and by integration of $s^k$, we obtain:
\begin{equation*}\label{preI-}
\begin{aligned}
&\mathbb{I}_-^j (\tau)= \\&   e^{-h^{-1}\varepsilon\varrho_-} \Bigg(\sum_{k=0}^{2i+1} (h^{-1} (\tau - \varepsilon)  \langle \xi \rangle   )^{k+1} \frac{ h  \langle \xi \rangle^{-1} }{k+1}C_{i,j,k} +  \Big( (\tau - \varepsilon)(\xi\cdot \dot{B}_{0,0} + \widehat{B}_{0,0}\Big) - ih^{-1}(\tau - \varepsilon)^{2}\frac{f\cdot\partial_{y}\varrho_-}{2}B_{0,0}\Bigg)\\& = e^{-h^{-1}\varepsilon\varrho_-} \sum_{k=0}^{2i+1} (h^{-1} (\tau - \varepsilon)  \langle \xi \rangle   )^{k+1} \frac{ h  \langle \xi \rangle^{-1} }{k+1}C_{i,j,k} +\\& e^{-h^{-1}\varepsilon\varrho_-}\Bigg( (h^{-1}(\tau - \varepsilon)\langle\xi\rangle)\Big(h\langle\xi\rangle^{-1} \xi\cdot \dot{B}_{0,0} + h \langle\xi\rangle^{-1} \widehat{B}_{0,0}\Big) -i \bigg(h^{-1}(\tau - \varepsilon) \langle\xi\rangle\bigg)^{2}\dfrac{h\langle\xi\rangle^{-1} f\cdot\partial_{y}\varrho_-}{2}B_{0,0}\Bigg),
\end{aligned}
\end{equation*} 	
then $e^{h^{-1} \tau \varrho_-} \Pi_- \mathbb{I}^{j}_-(\tau)$ has the following form:
\begin{equation}\label{I-}
\begin{aligned}
     & e^{h^{-1}\tau \varrho_-}  \Pi_-  \mathbb{I}_-^j (\tau)  = e^{ h^{-1}(\tau - \varepsilon)\varrho_-} Pi_- \,  \sum_{k=0}^{2i+1} \Big ( h^{-1} (\tau - \varepsilon) \langle\xi\rangle\Big)^{k+1} \dfrac{h \langle\xi\rangle^{-1}}{k+1} C_{i,j,k} \quad + \\&  e^{ h^{-1}(\tau - \varepsilon)\varrho_-} \Pi_- \Bigg( (h^{-1}(\tau - \varepsilon)\langle\xi\rangle)\Big(h\langle\xi\rangle^{-1} \xi\cdot \dot{B}_{0,0} + h \langle\xi\rangle^{-1} \widehat{B}_{0,0}\Big) -i \bigg(h^{-1}(\tau - \varepsilon) \langle\xi\rangle\bigg)^{2}\dfrac{h\langle\xi\rangle^{-1} f\cdot\partial_{y}\varrho_-}{2}B_{0,0}\Bigg).
    \end{aligned}
\end{equation} 	
For $\mathbb{I}_+^j $, let us introduce $\mathbb{P}_k$ the polynomial of degree $k$ such that 
$$ \int_\varepsilon^\tau e^{\lambda s} s^k \text{d}s = \frac{1}{\lambda^{k+1}} ( e^{\tau \lambda} \mathbb{P}_k(\tau \lambda) - e^{\varepsilon\lambda}\mathbb{P}_k(0) ), \quad \text{ for any }  \lambda \in \cc^* .$$
Using the above formula, then we obtain:
\begin{equation*}\label{preI+}
\begin{aligned}
\mathbb{I}_{+}^{j}& (\tau) = e^{-h^{-1}\varepsilon\varrho_-}\int_\varepsilon^\tau  e^{h^{-1} s( \varrho_- - \varrho_+)} \Bigg(\sum_{k=0}^{2i+1} (h^{-1} (s-\varepsilon)  \langle \xi \rangle   )^k C_{i,j,k} \quad + \quad \mathrm{m1}\Bigg)\, \mathrm{d}s\\& = e^{-h^{-1}\varepsilon\varrho_-} \sum_{k=0}^{2j+1} \dfrac{h \langle \xi \rangle^k}{(\varrho_- - \varrho_+)^{k+1}}\Big (e^{
h^{-1}\tau(\varrho_- - \varrho_+)} \mathbb{P}_{k}\big(h^{-1}(\tau - \varepsilon)(\varrho_- - \varrho_+) \big) - e^{h^{-1}\varepsilon(\varrho_- - \varrho_+)} \mathbb{P}_{k}(0) \Big)C_{i,j,k} \\& +e^{-h^{-1}\varepsilon\varrho_-} e^{h^{-1} (\varrho_- - \varrho_+)\tau}\Bigg[ h\Big ( \dfrac{\xi\cdot\dot{B}_{0,0} + \widehat{B}_{0,0}}{\varrho_- - \varrho_+}\Big)  +  i \dfrac{(\tau - \varepsilon)}{\varrho_- - \varrho_+}f\cdot\partial_{y}\varrho_- B_{0,0} - \quad  i\dfrac{\varepsilon h }{(\varrho_- - \varrho_+)^{2}}f\cdot\partial_{y}\varrho_- B_{0,0}\Bigg] \\& - e^{-h^{-1}\varepsilon\varrho_-} e^{h^{-1} (\varrho_- - \varrho_+)\varepsilon} \Bigg[ h\Big ( \dfrac{\xi\cdot\dot{B}_{0,0} + \widehat{B}_{0,0}}{\varrho_- - \varrho_+}\Big)  \quad - \quad  i\dfrac{\varepsilon h }{(\varrho_- - \varrho_+)^{2}}f\cdot\partial_{y}\varrho_- B_{0,0}\Bigg]. 
\end{aligned}
\end{equation*} 	
With this notation in hand, we easily see that the term $e^{h^{-1}\tau \varrho_+}  \Pi_+  \mathbb{I}_+^j (\tau)$ has the following form:
\begin{equation}\label{I+}
\begin{aligned}
    & e^{h^{-1}\tau \varrho_+}  \Pi_+  \mathbb{I}_+^j (\tau) = \\&  \Pi_+ \,  \sum_{k=0}^{2i+1} \frac{h \langle \xi \rangle^k}{(\varrho_- - \varrho_+)^{k+1}}  C_{i,j,k}  \Big(   e^{ h^{-1}(\tau - \varepsilon)\varrho_-} \mathbb{P}_k \big(h^{-1} (\tau - \varepsilon) ( \varrho_- - \varrho_+)  \big ) - e^{h^{-1}(\tau - \varepsilon) \varrho_+} \mathbb{P}_k(  0 )   \Big) \\& +  e^{h^{-1} (\tau - \varepsilon) \varrho_-}\Pi_+\Bigg[ h\Big ( \dfrac{\xi\cdot\dot{B}_{0,0} + \widehat{B}_{0,0}}{\varrho_- - \varrho_+}\Big) \quad + \quad i \dfrac{(\tau - \varepsilon)}{\varrho_- - \varrho_+}f\cdot\partial_{y}\varrho_- B_{0,0} - \quad  i\dfrac{\varepsilon h}{(\varrho_- - \varrho_+)^{2}}f\cdot\partial_{y}\varrho_- B_{0,0}\Bigg] \\& - e^{h^{-1} (\tau - \varepsilon) \varrho_+} \Pi_+\Bigg[ h\Big ( \dfrac{\xi\cdot\dot{B}_{0,0} + \widehat{B}_{0,0}}{\varrho_- - \varrho_+}\Big)  \quad - \quad  i\dfrac{\varepsilon h}{(\varrho_- - \varrho_+)^{2}}f\cdot\partial_{y}\varrho_- B_{0,0}\Bigg].
    \end{aligned}
\end{equation} 	
Thus, combining \eqref{I-} and \eqref{I+} with \eqref{Agene}, \eqref{2ndAgene} and \eqref{expL0}, yield that 
 \begin{equation}\label{preAj1}
 \begin{aligned}
 &A_{j+1}  = e^{h^{-1} (\tau - \varepsilon) \varrho_-} \Bigg[ \Pi_- {A_{j+1}}_{| \tau=\varepsilon}  + \Pi_- \,  \sum_{k=0}^{2i+1} \Big ( h^{-1} (\tau - \varepsilon) \langle\xi\rangle\Big)^{k+1} \dfrac{h \langle\xi\rangle^{-1}}{k+1} C_{i,j,k}\\&+\Pi_- \Bigg( \big(h^{-1}(\tau - \varepsilon)\langle\xi\rangle\big)\Big(h\langle\xi\rangle^{-1} \xi\cdot \dot{B}_{0,0} + h \langle\xi\rangle^{-1} \widehat{B}_{0,0}\Big) -i \bigg(h^{-1}(\tau - \varepsilon) \langle\xi\rangle\bigg)^{2}\dfrac{h\langle\xi\rangle^{-1} f\cdot\partial_{y}\varrho_-}{2}B_{0,0}\Bigg)\\&+  \Pi_+ \,  \sum_{k=0}^{2i+1} \frac{h \langle \xi \rangle^k}{(\varrho_- - \varrho_+)^{k+1}}  C_{i,j,k}  \Big( \mathbb{P}_k \big(h^{-1} (\tau - \varepsilon) ( \varrho_- - \varrho_+)  \big ) \Big)\Bigg]  \\& + e^{h^{-1} (\tau - \varepsilon) \varrho_-} \Pi_+\Bigg[  h\Big ( \dfrac{\xi\cdot\dot{B}_{0,0}+\widehat{B}_{0,0}}{\varrho_- - \varrho_+}\Big)  - i\big(h^{-1}(\tau - \varepsilon)\langle\xi\rangle\big) \dfrac{h\langle\xi\rangle^{-1} f\cdot\partial_{y}\varrho_- B_{0,0}}{\varrho_- - \varrho_+} -  i\dfrac{\varepsilon h}{(\varrho_- - \varrho_+)^{2}}f\cdot\partial_{y}\varrho_- B_{0,0}\Bigg] \\&+  e^{h^{-1} (\tau - \varepsilon) \varrho_+} \Bigg[ \Pi_+ {A_{j+1}}_{| \tau=\varepsilon} - \Pi_+ \,  \sum_{k=0}^{2i+1} \frac{h \langle \xi \rangle^k}{(\varrho_- - \varrho_+)^{k+1}}  C_{i,j,k}  \Big(\mathbb{P}_k(  0 )   \Big) \Bigg] \\& + e^{h^{-1} (\tau - \varepsilon) \varrho_+}\Pi_+\Bigg( h\Big ( \dfrac{\xi\cdot\dot{B}_{0,0} + \widehat{B}_{0,0}}{\varrho_- - \varrho_+}\Big)  \quad - \quad  i\dfrac{\varepsilon h}{(\varrho_- - \varrho_+)^{2}}f\cdot\partial_{y}\varrho_- B_{0,0}\Bigg).
 \end{aligned}
  \end{equation} 	
We set 
\begin{align}\label{TB+}
    \widetilde{{B}^+_{j+1}}  := \Pi_+ \,  \sum_{k=0}^{2i+1} \frac{h \langle \xi \rangle^k}{(\varrho_- - \varrho_+)^{k+1}}  C_{i,j,k}  \mathbb{P}_k(  0 )  - \Pi_+\Bigg( h\Big ( \dfrac{\xi\cdot\dot{B}_{0,0} + \widehat{B}_{0,0}}{\varrho_- - \varrho_+}\Big)  \quad - \quad  i\dfrac{\varepsilon h}{(\varrho_- - \varrho_+)^{2}}f\cdot\partial_{y}\varrho_- B_{0,0}\Bigg)
\end{align} belongs to $h\,\mathcal{S}^{0}$ 
as a linear combination of products of $\Pi_+ \in \sS^0$,  $h \langle \xi \rangle^k (\varrho_--\varrho_+)^{-k-1}$ $\in h \,\sS^{-1}$,  and of $C_{i,j,k}$ which verify the properties  as in \eqref{(ab)}. \\\\
Now, in order to have $A_{j+1} \in \sS^0$,  we let the contribution of the exponentially growing term vanish by choosing 
 \begin{equation*}
 \Pi_+ {A_{j+1}}(y, \xi, \varepsilon) = \widetilde{{B}^+_{j+1,k}} (y, \xi) . 
   \end{equation*} 
Then, we obtain 
 \begin{equation}\label{AJ1}
 \begin{aligned}
&A_{j+1}  = e^{h^{-1} (\tau - \varepsilon) \varrho_-} \Bigg[ \Pi_- {A_{j+1}}_{| \tau=\varepsilon}  + \Pi_- \,  \sum_{k=0}^{2i+1} \Big ( h^{-1} (\tau - \varepsilon) \langle\xi\rangle\Big)^{k+1} \dfrac{h \langle\xi\rangle^{-1}}{k+1} C_{i,j,k}\\&+\Pi_- \Bigg( \big(h^{-1}(\tau - \varepsilon)\langle\xi\rangle\big)\Big(h\langle\xi\rangle^{-1} \xi\cdot \dot{B}_{0,0} + h \langle\xi\rangle^{-1} \widehat{B}_{0,0}\Big) -i \bigg(h^{-1}(\tau - \varepsilon) \langle\xi\rangle\bigg)^{2}\dfrac{h\langle\xi\rangle^{-1} f\cdot\partial_{y}\varrho_-}{2}B_{0,0}\Bigg)\\&+  \Pi_+ \,  \sum_{k=0}^{2i+1} \frac{h \langle \xi \rangle^k}{(\varrho_- - \varrho_+)^{k+1}}  C_{i,j,k}  \Big( \mathbb{P}_k \big(h^{-1} (\tau - \varepsilon) ( \varrho_- - \varrho_+)  \big ) \Big)\Bigg]  \\& + e^{h^{-1} (\tau - \varepsilon) \varrho_-} \Pi_+\Bigg[  h\Big ( \dfrac{\xi\cdot\dot{B}_{0,0}+\widehat{B}_{0,0}}{\varrho_- - \varrho_+}\Big)  - i\big(h^{-1}(\tau - \varepsilon)\langle\xi\rangle\big) \dfrac{h\langle\xi\rangle^{-1} f\cdot\partial_{y}\varrho_- B_{0,0}}{\varrho_- - \varrho_+} - i\dfrac{\varepsilon h}{(\varrho_- - \varrho_+)^{2}}f\cdot\partial_{y}\varrho_- B_{0,0}\Bigg],
 \end{aligned}
  \end{equation} 	
since the boundary condition $P_+(y) {A_{j+1}}(y, \xi, \varepsilon)=0$, gives  
 \begin{equation*}
\Pi_-{A_{j+1}}(y, \xi, \varepsilon) = \Pi_-(P_+ + P_-){A_{j+1}}(y, \xi, \varepsilon) =  \Pi_- P_-{A_{j+1}}(y, \xi, \varepsilon) ,
   \end{equation*} 
   using the formula of $A_{j+1}(y,\xi,\tau)$ above, we get that 
   \begin{align*}
    P_-A_{j+1}(y,\xi,\varepsilon)=\dfrac{P_-\Pi_+}{k_-^{\varphi}}\widetilde{B^{+}_{j+1,k}},
   \end{align*}
   therefore 
   \begin{align}\label{CB+}
        \Pi_-A_{j+1}(y,\xi,\varepsilon)= \dfrac{\Pi_- P_-\Pi_+}{k_-^{\varphi}}\widetilde{B^{+}_{j+1,k}} . 
   \end{align}
In the other hand, regarding the following two series mentioned  in \eqref{preAj1}
   \begin{align}\label{sum}
 \Pi_- \,  \sum_{k=0}^{2i+1} \Big ( h^{-1} (\tau - \varepsilon) \langle\xi\rangle\Big)^{k+1} \dfrac{h \langle\xi\rangle^{-1}}{k+1} C_{i,j,k} +  \Pi_+ \,  \sum_{k=0}^{2i+1} \frac{h \langle \xi \rangle^k}{(\varrho_- - \varrho_+)^{k+1}}  C_{i,j,k}  \Big( \mathbb{P}_k \big(h^{-1} (\tau - \varepsilon) ( \varrho_- - \varrho_+)  \big ) \Big), 
   \end{align}
by calculation, it is easy to verify that for all $j \geq 2$ (\emph{i.e.}, $i \geq 1$), this quantity can be written as follows
\begin{align}\label{TB-}
    \sum_{k=0}^{2(j+1)} \big ( h^{-1}(\tau - \varepsilon)\langle\xi\rangle\big)^{k} \widetilde{B^{-}_{j+1,k}},
\end{align}
such that $\widetilde{B^{-}_{j+1,k}}$, as a linear combination, belong to $h^{2}\mathcal{S}^{0}$ for $k=0,...,2j+1$ and  $\widetilde{B^{-}_{j+1,2(j+1)}} \in h^{2}\mathcal{S}^{-1}$ .\\

Finally, the fact that we have the other terms (first and last) of the equality \eqref{AJ1} of order $h\,\mathcal{S}^{0}$ and admit the same structure as that of the terms in \eqref{sum}, then thanks to \eqref{CB+}, and \eqref{TB+}, \eqref{TB-},  together with  \eqref{AJ1} give that  
 $${A_{j+1}}(y, \xi, \tau ) = e^{h^{-1} (\tau - \varepsilon) \varrho_-(y, \xi)} \Bigg(  \dfrac{\Pi_- P_- \Pi_+ }{k^{\varphi}_-}  \widetilde{{B}^+_{j+1}} (y, \xi)  + \sum_{k=0}^{2(j+1)} (h^{-1} (\tau - \varepsilon)  \langle \xi \rangle   )^{k}  \widetilde{{B}^-_{j+1,k}}(y, \xi) \Bigg),$$
 where $\widetilde{{B}^+_{j+1}} (y, \xi),\,\widetilde{{B}^-_{j+1,k}}(y, \xi) \text{ belong to } h\,\mathcal{S}^{0}$, and Proposition \ref{SolSystSj} is proven with 
$$ {B}_{j+1,0} = \frac{\Pi_- P_+ \Pi_+ }{k^{\varphi}_+}  \widetilde{{B}^+_{j+1}} + \widetilde{{B}^-_{j+1,0}}, \text{ and for } k \geq 1, {{B}_{j+1,k}} = \widetilde{{B}^-_{j+1,k}}.\hspace{5cm}\qed$$
  \begin{proposition}\label{SobolevMap}
Let $A_j$, $j\geq 0$, be of the form \eqref{eqAj}. Then, for any $s \geq -\frac{1}{2}$, the operator ${\mathcal A}_j$ defined by 
\begin{align*}
	{\mathcal A}_j : f \longmapsto  ({\mathcal A}_j f) (y,y_3) = \frac{1}{(2 \pi)^2} \int_{\rr^2} A_j(y, h \xi , y_3) e^{i y \cdot \xi} \hat{f}(\xi) \mathrm{d} \xi
\end{align*}
gives rise to a bounded operator from $\mathit{H}^s(\rr^2)$ into $\mathit{H}^{s+\frac12}(\rr^2\times (\varepsilon,+ \infty))$. Moreover,  for any $l \in [0,\frac12]$ we have:
\begin{equation}\label{SobolevEsti}
	\| {\mathcal A}_j\|_{\mathit{H}^s \to \mathit{H}^{s+\frac12-l}}  =  O(h^{l-|s|+1}).
\end{equation} 		
\end{proposition}
\textbf{Proof.} The proof of this proposition follows exactly the arguments of \cite[Proposition 5.4]{BBZ}. However, this difference obtained at the rate level on $h$ is because of the presence of a parameter $h$ in the terms $B_{j,k}$ of the solution $A_j$.
\begin{proposition}\label{ApproxSystSC}
Let $f \in H^s(\rr^2)$ and $A_j$, $j\geq 0$, be as in Propositions \ref{SolSystS0}, \ref{SolSystSj}. Then for any $N\in \mathbb{N},$ the function  
$u_N^h = \sum_{j=0}^N h^j {\mathcal A}_j f$ satisfies: 
\begin{equation}\label{4.53}
\left\{
\begin{aligned}
h \partial_{\tau} u_N^h & - L_0(y, h D_y)  u_N^h - h \sum_{k=1}^{+\infty}h^{k-1}\dfrac{(\alpha\cdot\nu^{\varphi}c_3)^{k-1}}{(1+|\nabla \chi |^{2})^{\frac{k-1}{2}}} \widetilde{L}_{1}(y,h D_{y}) u_N^h =  h^{N+1} {\mathcal R}_N^h f ,  &  \quad \text{ in }  &  \rr^2 \times (\varepsilon, + \infty),\\
P_+ u_N^h  & = f,  &  \text{ on } &  \rr^2 \times \{\varepsilon\},
\end{aligned}
\right.
\end{equation}
with 
\begin{align*}
&\mathcal R_N^h f =
  \frac{-1}{(2 \pi)^2} \int_{\rr^2} \Bigg( \sum_{k=1}^{+\infty}h^{k}\dfrac{(\alpha\cdot\nu^{\varphi}c_3)^{k-1}}{(1+|\nabla \chi |^{2})^{\frac{k-1}{2}}}\Big( h^{-1} \widetilde{L}_{1} A_N - i \partial_{\xi}\widetilde{L}_{1}\cdot\partial_{y}A_N\Big) - i \partial_\xi L_0  \cdot \partial_y  A_N \Bigg)  e^{i y \cdot \xi} \hat{f}(\xi) \mathrm{d} \xi,
\end{align*}
a bounded operator from $H^s(\rr^2)$ into $H^{s+\frac12}(\rr^2\times (\varepsilon,+ \infty))$ satisfying  for any $l \in [0, \frac12]$:
\begin{equation}\label{EstiRN}
\| {\mathcal R}^h_N\|_{\mathit{H}^s \to \mathit{H}^{s+\frac12-l}}  =  O(h^{l-|s|+1}).
\end{equation} 	
\end{proposition}
\textbf{Proof.} 
By construction of the sequence $(A_j)_{j \in \{0,\cdots,N-1\} }$ as in \eqref{T1bh'}, we have the system \eqref{4.53} with  ${\mathcal R}^h_N = \mathit{Op}^h (r_N^h(\cdot, \cdot, \tau) )$, such that
$$ r_N^h(y, \xi , \tau) = - \Bigg(\sum_{k=1}^{+\infty}h^{k}\dfrac{(\alpha\cdot\nu^{\varphi}c_3)^{k-1}}{(1+|\nabla \chi |^{2})^{\frac{k-1}{2}}}\Big( h^{-1} \widetilde{L}_{1} A_N - i  \partial_{\xi}\widetilde{L}_{1}\cdot\partial_{y} A_N\Big) - i \partial_\xi L_0  \cdot \partial_y  A_N \Bigg) (y, h \xi , \tau) .$$
As in the proof of  Proposition \ref{SolSystSj}, $\widetilde{L}_{1}A_N$ has the form \eqref{Agene2}, and $\partial_{\xi}\widetilde{L}_{1}\cdot\partial_{y}A_N$ and  $\partial_{\xi}L_{0}\cdot\partial_{y}A_N$ have the form \eqref{Agene3}. Then, $r_N^h$ has the form \eqref{(a)} (with $j=N$).  Therefore,  as in the proof of Proposition \ref{SobolevMap}, we obtain the estimate \eqref{EstiRN}.\qed
\begin{proposition}\label{prop4.6}
Let us consider the Poincaré-Steklov operator $\mathcal{A}^{h}$ introduced at the beginning of Section \ref{section4}. For $h=\varepsilon\in (0,1]$ and for all $N\in \mathbb{N}$, there is a $h$-pseudodifferential operator of order $0$, $\mathcal{A}^{h}_{N}$ such that for $h$ sufficiently small, we have the following estimate:
\begin{align}\label{est}
|| \mathcal{A}^{h} - \mathcal{A}^{h}_N||_{H^{1/2}(\Sigma^{\varepsilon})\rightarrow H^{\frac{3}{2} - l} (\S^{\varepsilon})} = O(h^{2l + \frac{1}{2}}), \quad \text{ for any } l\in [0,\frac{1}{2}].
\end{align}
\end{proposition}
\textbf{Proof.} The proof of this proposition follows the same argument of \cite[Theorem 5.1]{BBZ}. That is a consequence of the above Proposition \ref{SobolevMap} and \ref{ApproxSystSC}, combined with the regularity estimates from Theorem \ref{MITLipshictz}-(iii). More precisely, let $(U^{\varepsilon}_\varphi, V^{\varepsilon}_{\varphi},\varphi^{\varepsilon})$ a chart of an atlas $\mathbb{A}^{\varepsilon}$ of $\Sigma^{\varepsilon}$, and $\psi_{1},\psi_{2}\in C^{\infty}_{0}(U^{\varepsilon}_{\varphi})$. Let also $h^{\varepsilon}\in P_-H^{1/2}(\Sigma^{\varepsilon})$ be such that $f^{\varepsilon}:=(\varphi_{\varepsilon}^{-1})^* [\psi_2 h^{\varepsilon}]\in H^{1/2}(V^{\varepsilon}_\varphi)^4$,  which can be extended by 0 to a function of $H^{1/2}(\mathbb{R}^2)^4$. Then, for $\varepsilon=h=\kappa^{-1}$ and $N\in\mathbb{N},$ the previous construction provides a function $u^{h}_{N}\in H^{1}(\mathbb{R}^{2}\times (\varepsilon,+\infty))^{4}$ which verifies the following system
\begin{equation*}
	\left\{
\begin{aligned}\label{}	
&(\widetilde{D}^\varphi_\kappa  -z )u^{h}_{N}	=h^{N+1}\mathcal{R}_N^{h}f^{\varepsilon},  \quad  &\text{ in } & \mathbb{R}^2 \times (\varepsilon, +\infty) ,\\
	&P^\varphi_- t_{\Sigma^{\varepsilon}} u^{h}_{N}= f^{\varepsilon},  \quad    &\text{ on }& \mathbb{R}^2 \times \{\varepsilon\}, 
	\end{aligned}
	\right.
	\end{equation*}
 where $u^{h}_{N}$, $\mathcal{R}_N^{h}$ are defined in Proposition \ref{ApproxSystSC}. Moreover, from the latter, we know that $\mathcal{R}_N^{h}\in H^{N+1}(\mathbb{R}^2 \times (\varepsilon, +\infty))$ with norm in $\mathit{H}^{1-l},$ $l\in[0,\frac{1}{2}],$ bounded by $O(h^{l+\frac{1}{2}}).$ Consequently, $v^{h}_{N}:=\phi_{\varepsilon}^* u^{h}_{N},$ defined on $\mathcal{V}_{\varphi,\eta}^{\varepsilon}$ satisfies
 \begin{equation*}
	\left\{
\begin{aligned}\label{}	
&(D_\kappa  -z )v^{h}_{N}	=h^{N+1}(\phi_{\varepsilon}^{-1})^{*}\big(\mathcal{R}_N^{h}f^{\varepsilon}\big),  \quad  &\text{ in } & \mathcal{V}^{\varepsilon}_{\varphi,\eta} ,\\
	&P_- t_{\Sigma^{\varepsilon}} v^{h}_{N}= \psi_{2}h^{\varepsilon},  \quad    &\text{ on }& U^{\varepsilon}_{\varphi}.
	\end{aligned}
	\right.
	\end{equation*}
Recall the definition of the lifting operator $\mathcal{E}^{\varepsilon}_\kappa$, given in Definition \ref{def}. We have for $h^{\varepsilon}\in P_- H^{1/2}(\Sigma^{\varepsilon})^{4}$, $\mathcal{E}^{\varepsilon}_\kappa[\psi_2 h^{\varepsilon}] \in H^{1}(\mathcal{U}^{\varepsilon})^4$. Since $P_- t_{\Sigma^{\varepsilon}} v^{h}_{N} = P_- t_{\Sigma^{\varepsilon}} \mathcal{E}^{\varepsilon}_\kappa[\psi_2 h^{\varepsilon}]=\psi_{2}h^{\varepsilon}$, it follows that 
$$v^{h}_{N} -  \mathcal{E}^{\varepsilon}_\kappa[\psi_2 h^{\varepsilon}] = h^{N+1} (D^{\varepsilon}_{\mathrm{MIT}}(\kappa) - 1)^{-1} (\phi_{\varepsilon}^{-1})^{*}\big(\mathcal{R}_N^{h} (\varphi_{\varepsilon}^{-1})^* [\psi_2 h^\varepsilon]\big).$$ Thanks to the estimation of \cite[Theorem 3.2-(i)]{BBZ}, and also by continuing the steps of the proof of Theorem 5.1 in \cite{BBZ}, we obtain that $\mathcal{A}_{N}^{h} \in h \, Op^{h}\mathcal{S}^{0}(\Sigma^{\varepsilon})$ and the estimate \eqref{est} holds for any $l\in[0,\frac{1}{2}]$.\qed\\

At the end of this section, let's give some pseudodifferential properties of the Poincaré-Steklov operators, $\mathscr{A}_m$ and $\mathscr{A}_{m}^\varepsilon$, introduced in Definition \ref{def}, in order to use it in Section \ref{app}.
\begin{remark}\label{remark5.2}
 We mention that the fixed Poincaré-Steklov operator $\mathscr{A}_m$ have been introduced and studied in details in the paper \cite[Theorem 4.1]{BBZ}. Moreover, it is a pseudodifferential operator of order $0$, which can be considered as a $h$-pseudodifferential operator, and whose semiclassical principal symbol (in local coordinate) is given by
      \begin{align*}
     \mathscr{P}_{h,\mathscr{A}_m}(x_{\Sigma},\xi)= \dfrac{S\cdot(\xi\wedge\nu(x_{\Sigma}))}{|\xi\wedge\nu(x_\Sigma)|}P_-,\quad \text{for any }  x_\Sigma\in\Sigma.
\end{align*}
  \end{remark}
For $\mathscr{A}_m^\varepsilon,$ we have the following results:
 \begin{theorem}\label{pseudoAme}
 Let $z\in\rho(D_m)$ and $x_\Sigma \in \Sigma$ and recall the definition of $\mathcal{T}_\varepsilon$ from  Definition \ref{Transformationoperator}. We define the Cauchy operator $\mathscr{C}_m^{z,\varepsilon}: L^{2}(\Sigma^\varepsilon)^4 \longrightarrow  L^{2}(\Sigma^\varepsilon)^4$ as the singular integral operator acting as 
      \begin{align*}\label{}
\mathscr{C}^{z,\varepsilon}_{m}[g](x)&:=  \lim\limits_{\rho\searrow 0}\int_{|x-y|>\rho}\phi^{z}_{m}(x-y)g(y)\mathrm{d}\sigma(y), \quad \text{for } \mathrm{d}\sigma\emph{-a.e.,}\,x=x_\Sigma + \varepsilon\nu(x_\Sigma)\in\Sigma^\varepsilon,\,g\in\mathit{L}^2(\Sigma^\varepsilon)^4.
\end{align*} 
Also, we consider the Poincaré-Steklov operator $\mathscr{A}^{\varepsilon}_m$ given in Definition \ref{def}. Then, $\mathcal{T}_\varepsilon^{-1} \mathscr{C}_m^{z,\varepsilon} \mathcal{T}_\varepsilon$ and $\mathcal{T}_\varepsilon^{-1} \mathscr{A}_m^{\varepsilon} \mathcal{T}_\varepsilon$ are homogeneous pseudodifferential operators of order $0$, and we have 
\begin{align*}
    &\mathcal{T}_\varepsilon^{-1} \mathscr{C}_m^{z,\varepsilon} \mathcal{T}_\varepsilon = \mathrm{det}(1-\varepsilon W(x_\Sigma)) \Big[\dfrac{1}{2}\alpha\cdot\dfrac{\nabla_{\Sigma}}{\sqrt{-\Delta_\Sigma}} + \varepsilon\, Op(b_0(x_\Sigma, \xi)) + Op(b_{-1}(x_\Sigma,\varepsilon\xi))\Big],\\
    &\mathcal{T}_\varepsilon^{-1} \mathscr{A}_m^{\varepsilon} \mathcal{T}_\varepsilon =\mathrm{det}(1-\varepsilon W(x_\Sigma))\Big[ S\cdot\dfrac{(\nabla_{\Sigma}\wedge\nu)}{\sqrt{-\Delta_\Sigma}} P^\varepsilon_- + \varepsilon\, Op(b^{'}_0(x_\Sigma, \xi)) + Op(b^{'}_{-1}(x_\Sigma,\varepsilon\xi))\Big],
\end{align*}
where $\nabla_\Sigma = \nabla - \nu(\nu\cdot\nabla)$ is the surface gradient along $\Sigma$, and $-\Delta_\Sigma$ is the Laplace-Beltrami operator, with $b_0$, $b_0^{'},$ resp. $b_{-1}$, $b_{-1}^{'}$ the symbols of order $0$, resp. $-1$. 
  \end{theorem}
  \textbf{Proof.} The proof follows similar arguments as in \cite[Theorem 4.1]{BBZ}. Let $f\in L^2(\Sigma)^4$ and consider the operator $\mathcal{T}_\varepsilon^{-1} \mathscr{C}_m^{z,\varepsilon} \mathcal{T}_\varepsilon f.$ Using the explicit formula of $\mathscr{A}_m^\varepsilon$, we have the following connection  
  $$ L^2(\Sigma)^4 \ni \mathcal{T}_\varepsilon^{-1} \mathscr{A}_m^{\varepsilon} \mathcal{T}_\varepsilon f= - P^\varepsilon_+ \beta \Big(\beta/2 + \mathcal{T}_\varepsilon^{-1} \mathscr{C}_m^{z,\varepsilon} \mathcal{T}_\varepsilon f\Big)^{-1} P^\varepsilon_- .$$
  Now, fix a local chart $(U, V, \varphi)$ of $\Sigma$ and let $\psi_k : \Sigma \longrightarrow \mathbb{R},$ $k=1,2,$ be a $C^\infty$-smooth function with $\mathrm{supp(\psi_1)}\cap \mathrm{supp}(\psi_2)=\emptyset$.
 For $x_\Sigma \in \Sigma,$ 
\begin{equation}\label{TCT}
\begin{aligned}
     \Big(\mathcal{T}_\varepsilon^{-1} \mathscr{C}_m^{z,\varepsilon} \mathcal{T}_\varepsilon f\Big)(x_\Sigma) &= \mathrm{det}(1-\varepsilon W(x_\Sigma))\,\mathrm{p.v.}\int_{|x_\Sigma + \varepsilon\nu(x_\Sigma) - y|>\rho}\phi^z_m(x_\Sigma + \varepsilon\nu(x_\Sigma) - y)\mathcal{T}_\varepsilon f (y) \mathrm{d}\sigma(y)\\& = \mathrm{det}(1-\varepsilon W(x_\Sigma))\,\mathrm{p.v.}\int_{|x_\Sigma - y_\Sigma|>\rho^{'}} \phi^z_m(x_\Sigma + \varepsilon\nu(x_\Sigma) - y_\Sigma - \varepsilon \nu(y_\Sigma))f (y_\Sigma) \mathrm{d}\sigma(y_\Sigma)\\& = \mathrm{det}(1-\varepsilon W(x_\Sigma))\int_{V}\phi^z_m(x_\Sigma - y_\Sigma + \varepsilon \big(\nu(x_\Sigma) - \nu(y_\Sigma))\big)f (y_\Sigma) \mathrm{d}\sigma(y_\Sigma).
\end{aligned}
\end{equation}
Now, recall the definition of $\phi^z_m$ from \eqref{defsolo}, and observe that 
$$ \phi_m^z(x-y) =  k (x-y) + a (x-y),$$
where 
		\begin{align*}\label{}
			k^{z}(x-y)&=\frac{e^{i\sqrt{z^2-m^2}|x-y|}}{4\pi|x-y|}\left(z +m\beta+ \sqrt{z^2-m^2}\alpha\cdot\frac{x-y}{|x-y|}\right) +i\frac{e^{i\sqrt{z^2-m^2}|x-y|}-1}{4\pi|x-y|^3}\alpha\cdot(x-y),\\
			a(x-y)&=\frac{i}{4\pi|x-y|^3}\alpha\cdot(x-y).
		\end{align*}  
		Using this, it follows that 
		\begin{align*}\label{}
			\begin{split}
				\mathscr{C}^{z,\varepsilon}_{m}[g](x)=&   \lim\limits_{\rho\searrow 0}\int_{|x-y|>\rho}a(x-y)g(y)\mathrm{d\sigma}(y)+\int_{\Sigma^\varepsilon}k^{z}(x-y)g(y)\mathrm{d\sigma}(y)\\
				=& A[g](x)+ K[g](x).
			\end{split}
		\end{align*}
As $|k^{z}(x-y)|=\mathcal{O}(|x-y|^{-1})$ when $|x-y|\rightarrow0$,  using the standard layer potential techniques 
		(see, \emph{e.g.} \cite[Chap. 3, Sec. 4]{T2} and \cite[Chap. 7, Sec. 11]{T1}) it is not hard to prove that the integral operator $\mathcal{T}_\varepsilon^{-1} K \mathcal{T}_\varepsilon$ gives rise to a pseudodifferential operator of order $-1$, \emph{i.e.}, $\mathcal{T}_\varepsilon^{-1} K \mathcal{T}_\varepsilon\in  \mathit{Op} \sS^{-1}(\S)$. Thus,  we can (formally) write
		\begin{align}\label{rela}
			\mathcal{T}_\varepsilon^{-1} \mathscr{C}_{m}^{z,\varepsilon} \mathcal{T}_\varepsilon = \mathcal{T}_\varepsilon^{-1}  A \mathcal{T}_\varepsilon \quad \mathrm{mod} \, \mathit{Op} \sS^{-1}(\S),
		\end{align}
		which means that the operator $A$ encodes the main contribution in the pseudodifferential character of $\mathcal{T}_\varepsilon^{-1} \mathscr{C}^{z,\varepsilon}_{m}\mathcal{T}_\varepsilon$. \\
  
For $\Sigma^\varepsilon \ni x = x_\Sigma  + \varepsilon \nu (x_\Sigma)$, $y = y_\Sigma  + \varepsilon \nu (y_\Sigma)$,
$$ a\big( x_\Sigma - y_\Sigma + \varepsilon (\nu(x_\Sigma) - \nu(y_\Sigma))\big)= i \alpha\cdot\dfrac{\big( x_\Sigma - y_\Sigma + \varepsilon (\nu(x_\Sigma) - \nu(y_\Sigma))\big)}{\big| x_\Sigma - y_\Sigma + \varepsilon (\nu(x_\Sigma) - \nu(y_\Sigma))\big|^3}.$$
Set $X= x_\Sigma - y_\Sigma.$ Then, $| x_\Sigma - y_\Sigma + \varepsilon (\nu(x_\Sigma) - \nu(y_\Sigma))| = |X+\varepsilon\nu X|.$ And $|X + \varepsilon\nu X|^{-3}$ yields
$$ |X + \varepsilon\nu X|^{-3} = (1+\varepsilon^2)^{-3/2} |X|^{-3} \Big(  1+ 2\varepsilon (1+\varepsilon^2)^{-1} \dfrac{\langle X,\nu X \rangle}{|X|^2}\Big)^{-3/2}.$$
By a series expansion (first order), we get 
$$|X + \varepsilon\nu X|^{-3}= |X|^{-3} + \varepsilon \Big( -3 |X|^{-3} \dfrac{\langle X,\nu X \rangle}{|X|^2} \Big). $$
For any $X\in U$ we have  $X=(\widetilde{X}, \chi(\widetilde{X}))$ with $X\in V$ and where the graph of $\chi : V \longrightarrow \mathbb{R}$ coincides with $U.$ With the same argument in \cite[Theorem 4.1]{BBZ} we get that, uniformly with respect to $\varepsilon\in(0,\varepsilon_0),$ with $\varepsilon_0$ sufficiently small
\begin{align*}
    |X + \varepsilon \nu X|^{-3} &= \dfrac{1}{\langle\widetilde{X}, G(\widetilde{x}_\Sigma)\widetilde{X}\rangle^{3/2}} + k_1(\widetilde{X}), \quad &\text{ with } |k_1(\widetilde{X})| = \mathcal{O}(|\widetilde{X}|^{-2})\text{ when } |\widetilde{X}|\longrightarrow 0,\\
        |X + \varepsilon \nu X|^{-5}\langle X,\nu X\rangle &= \dfrac{\langle \widetilde{X},\nu \widetilde{X}\rangle}{\langle\widetilde{X}, G(\widetilde{x}_\Sigma)\widetilde{X}\rangle^{5/2}} + \langle \widetilde{X},\nu \widetilde{X}\rangle k_2(\widetilde{X}), \quad &\text{ with } |k_2(\widetilde{X})| = \mathcal{O}(|\widetilde{X}|^{-4})\text{ when } |\widetilde{X}|\longrightarrow 0,
\end{align*}
where $G(\widetilde{x}_\Sigma)$ is the metric tensor.
We deduce that 
\begin{align}\label{rela1}
    \psi_2( \mathcal{T}_\varepsilon^{-1} A \mathcal{T}_\varepsilon [\psi_1 f] )(x_\Sigma)= \psi_2 Op(a_0(x_{\Sigma}, \xi))\psi_1 f (x_\Sigma) + \varepsilon\, \psi_2 Op(b_0(x_{\Sigma}, \xi))\psi_1 f (x_\Sigma) + \psi_2 L \psi_1, 
\end{align} 
with $L$ a pseudodifferential operator of order $-1.$ Thus, $\mathcal{T}_\varepsilon^{-1} A \mathcal{T}_\varepsilon $ is a zero-order pseudodifferential operator. Furthermore, thanks to  \eqref{rela} and \eqref{rela1} we get that $\mathcal{T}_\varepsilon^{-1} \mathscr{C}^{z,\varepsilon}_m \mathcal{T}_\varepsilon $ is a homogeneous pseudodifferential operator of order $0$, with principal symbol given by 
\begin{align*}
    \mathcal{T}_\varepsilon^{-1} \mathscr{C}_m^{z,\varepsilon} \mathcal{T}_\varepsilon = \mathrm{det}(1-\varepsilon W(x_\Sigma)) \Big[\dfrac{1}{2}\alpha\cdot\dfrac{\nabla_{\Sigma}}{\sqrt{-\Delta_\Sigma}} + \varepsilon\, Op(b_0(x_\Sigma, \xi)) + Op(b_{-1}(x_\Sigma,\varepsilon\xi))\Big].
\end{align*}
Consequently, thanks to the relation between $\mathscr{C}^{z,\varepsilon}_m$ and $\mathscr{A}_m^\varepsilon$, we have that $\mathcal{T}_\varepsilon^{-1} \mathscr{A}_m^{\varepsilon} \mathcal{T}_\varepsilon$ is a homogeneous pseudodifferential operators of order 0
\begin{align*}
\mathcal{T}_\varepsilon^{-1} \mathscr{A}_m^{\varepsilon} \mathcal{T}_\varepsilon =\mathrm{det}(1-\varepsilon W(x_\Sigma))\Big[ S\cdot\dfrac{(\nabla_{\Sigma}\wedge\nu(x_\Sigma))}{\sqrt{-\Delta_\Sigma}} P^\varepsilon_- + \varepsilon\, Op(b^{'}_0(x_\Sigma, \xi)) + Op(b^{'}_{-1}(x_\Sigma,\varepsilon\xi))\Big].
\end{align*}\qed
\begin{corollary}\label{corollsry5.1}
    The Poincaré-Steklov operator $\mathscr{A}_m^{\varepsilon}$ is a homogeneous pseudodifferential operator of order $0$, and we have that
\begin{align*}
\mathscr{A}_m^{\varepsilon} &=  S\cdot\dfrac{(\nabla_{\Sigma^\varepsilon}\wedge N^\varepsilon(p(x_\Sigma)))}{\sqrt{-\Delta_{\Sigma^\varepsilon}}} P^\varepsilon_- + \varepsilon\, Op(b^{p}_0(x_\Sigma, \xi)) + Op(b^{p}_{-1}(x_\Sigma,\varepsilon\xi))\\&= - S\cdot\dfrac{(\nabla_{\Sigma^\varepsilon}\wedge \nu(x_\Sigma))}{\sqrt{-\Delta_{\Sigma^\varepsilon}}} P^\varepsilon_- + \varepsilon\, Op(b^{p}_0(x_\Sigma, \xi)) + Op(b^{p}_{-1}(x_\Sigma,\varepsilon\xi)),\quad\text{with }\varepsilon\in(0,\varepsilon_0),
\end{align*}
where $\nabla_{\Sigma^\varepsilon}$ is the surface gradient along $\Sigma^\varepsilon$, $-\Delta_{\Sigma^\varepsilon}$ is the Laplace-Beltrami operator, and $b_j^{p}(x_\Sigma, \xi)$ has the following form
\begin{align*}
    b_j^{p}(x_\Sigma, \xi)= b_j \bigg( p(x_\Sigma), \Big({\nabla {p(x_\Sigma)}^{-1}}\Big)^{t}\xi\bigg),\quad \text{for } j\in \lbrace -1,0\rbrace,
\end{align*}
with $p(x_\Sigma) = x_\Sigma + \varepsilon \nu (x_{\Sigma})$ the diffeomorphism from Definition \ref{Transformationoperator}, and $\Big({\nabla {p(x_\Sigma)}^{-1}}\Big)^{t} = \Big({\big(1-\varepsilon W(x_\Sigma)\big)^{-1}}\Big)^{t} = \big(1-\varepsilon W(x_\Sigma)\big)^{-1},$ where $W(x_\Sigma)$ is the Weingarten matrix, symmetric, given in Definition \ref{Weingarten}.
\end{corollary}
\textbf{Proof.} The proof of this corollary is a consequence of Theorem \ref{pseudoAme} and the arguments of \cite[Theorem 9.3]{MZ}.
\qed
\section{Reduction to a MIT bag problem.}\label{app}
Throughout the section, we denote $\Omega_+$,  $\O^{\varepsilon}_{-}$ and  $\mathcal{U}^{\varepsilon}$ the domains as in Figure \ref{figure1} such that $\S=\partial\Omega_+$, $\S^{\varepsilon}:=\partial\O^{\varepsilon}_-$ and  $\partial\mathcal{U}^{\varepsilon}=\S\cup\S^{\varepsilon},$ respectively,  and we let $N^{\varepsilon}$ be the outward pointing unit normal to $\O^{\varepsilon}_{-}.$ We set $\nu$ the outward unit normal to the fixed domain $\O_+\subset\mathbb{R}^3$. Fix $m>0$ and let $M>0$. Remember our perturbed Dirac operator 
\begin{align*}   
\mathfrak{D}^{\varepsilon}_{M}\varphi=(D_m+M\beta \mathbbb{1}_{\mathcal{U}^{\varepsilon}})\varphi,\quad\forall\varphi\in\mathrm{dom}(\mathfrak{D}_M^{\varepsilon}):=H^1(\rr^{3})^{4},
\end{align*}
where $\mathbbb{1}_{\mathcal{U}^{\varepsilon}}$ is the characteristic function of $\mathcal{U}^{\varepsilon}$.\\\\
Let us now recall the definition of the MIT bag operator from Section \ref{DefMIT} by $D^{\Omega_+}_{\mathrm{MIT}}$, $D^{\Omega^\varepsilon_-}_{\mathrm{MIT}}$, and $D^{\mathcal{U}^\varepsilon}_{\mathrm{MIT}}$, which act in $L^{2}(\O_+)^{4}, \,L^{2}(\O^{\varepsilon}_-)^{4}$, and $L^{2}(\mathcal{U}^{\varepsilon})^{4}$ repsectively. The aim of this section is to use the properties of the Poincaré-Steklov operators carried out in the previous sections to study the resolvent of $\mathfrak{D}_{M}^{\varepsilon}$ when $M$ is large enough. Namely, we give a Krein-type resolvent formula of the Dirac operator $\mathfrak{D}^{\varepsilon}_{M}$ in terms of the resolvent of the MIT bag operator $D^{\Omega_+}_{\mathrm{MIT}}\oplus D_{\mathrm{MIT}}^{\Omega^\varepsilon_-}$, and we show that the convergence of $\mathfrak{D}_{M}^{\varepsilon}$ toward $\mathscr{D}_L$ holds in the norm resolvent sense when $M$ and $\varepsilon$ converge to $\infty$ and $0^{+}$, respectively. To set up Krein's formula between the resolvent of $\mathfrak{D}^{\varepsilon}_{M}$ and the resolvent of $D^{\Omega_+}_{\mathrm{MIT}}\oplus D_{\mathrm{MIT}}^{\Omega^\varepsilon_-}$, we will fix $\nu$ the only normal acting in our domain. Throughout this section, the projections associated with the surface $\Sigma^\varepsilon$ (\emph{i.e.,} $P_{\pm}^{\varepsilon}(x)$, for $x\in\Sigma^{\varepsilon}$)  verify the properties of Remark \ref{remark3.1}.
\subsection{Notations}
Let $z\in\rho(D^{\Omega^\varepsilon_-}_{\mathrm{MIT}})\cap\rho(\mathfrak{D}^{\varepsilon}_M)$. We recall $\Omega_{+-}^{\varepsilon}:=\Omega_+\cup\Omega^{\varepsilon}_-$. We define the resolvents associated with the operators   $\mathfrak{D}_{M}^{\varepsilon}$, $D^{\mathcal{U}^\varepsilon}_{\mathrm{MIT}}$, and $D^{\Omega^\varepsilon_{+-}}_{\mathrm{MIT}}:=D^{\Omega_+}_{\mathrm{MIT}}\oplus D^{\Omega^\varepsilon_-}_{\mathrm{MIT}}$, respectively, by
\begin{itemize}
    \item $R_{M}^{\varepsilon}(z):=(\mathfrak{D}^{\varepsilon}_{M}-z)^{-1}:L^{2}(\mathbb{R}^{3})^{4}\rightarrow H^{1}(\mathbb{R}^{3})^{4}$.
\item $R^{\mathcal{U}^\varepsilon}_{\mathrm{MIT}}(z):=(D^{\mathcal{U}^\varepsilon}_{\mathrm{MIT}}-z)^{-1}:L^{2}(\mathcal{U}^{\varepsilon})^{4}\rightarrow$ $\mathrm{dom}$($D^{\mathcal{U}^\varepsilon}_{\mathrm{MIT}}$).
 \item $R_{\mathrm{MIT}}^{\Omega^{\varepsilon}_{+-}}(z):=(D^{\Omega^\varepsilon_{+-}}_{\mathrm{MIT}}-z)^{-1}:L^{2}(\Omega_+)^4 \oplus L^2(\Omega^\varepsilon_-)^{4}\rightarrow$ $\mathrm{dom}(D^{\Omega_+}_{\mathrm{MIT}})\oplus\mathrm{dom}$($D^{\Omega_-^\varepsilon}_{\mathrm{MIT}})\subset L^{2}(\Omega_+)^4 \oplus L^2(\Omega^\varepsilon_-)^{4}$
\end{itemize}
can be read as the following matrix:
\begin{align}\label{MRMIT}
    R_{\mathrm{MIT}}^{\Omega^{\varepsilon}_{+-}} = \begin{pmatrix}
	R^{\Omega_+}_{\mathrm{MIT}}r_{\Omega_+} & 0\\
	0 & R^{\Omega^\varepsilon_-}_{\mathrm{MIT}}r_{\Omega^\varepsilon_-}
	\end{pmatrix}\equiv r_{\Omega^\varepsilon_{+-}}e_{\Omega_+}R^{\Omega_+}_{\mathrm{MIT}}r_{\Omega_+} + r_{\Omega^\varepsilon_{+-}}e_{\Omega^\varepsilon_-}R^{\Omega^\varepsilon_-}_{\mathrm{MIT}}r_{\Omega_-^\varepsilon} =\big(R^{\Omega_+}_{\mathrm{MIT}}r_{\O_+},R_{\mathrm{MIT}}^{\Omega_-^\varepsilon}r_{\O^\varepsilon_-}\big),
\end{align}
where $R^{\Omega_+}_{\mathrm{MIT}}(z)$, $R^{\Omega^\varepsilon_-}_{\mathrm{MIT}}(z)$ are the resolvents of $D^{\Omega_+}_{\mathrm{MIT}}$, $D_{\mathrm{MIT}}^{\Omega^\varepsilon_-},$ respectively, and $r_{\Omega^\varepsilon_{+-}}$, $e_{\Omega^\varepsilon_{+-}}$ are defined below.\\\\
We define $r_{\Omega^\varepsilon
_{+-}}$ and $e_{\Omega^\varepsilon_{+-}}$ as the restriction operator in $\Omega^\varepsilon_{+-}$ and its adjoint operator, \emph{i.e.}, the extension by $0$ outside of $\Omega^\varepsilon_{+-}$, respectively, by
\begin{equation}\label{re}
\begin{aligned}
\begin{array}{rcl}
r_{\Omega^\varepsilon_{+-}}: L^{2}(\rr^3)^4 &\to &  L^{2}(\O_+)^4\oplus L^{2}(\O^{\varepsilon}_-)^4\\
w &\mapsto & r_{\Omega^\varepsilon_{+-}}w:=( r_{\O_{+}}w \oplus r_{\O^{\varepsilon}_-} w)\equiv (r_{\O_+},r_{\O^{\varepsilon}_-})w,\\\vspace{0.1cm}
\end{array}\\
\begin{array}{rcl}
e_{\Omega^{\varepsilon}_{+-}}: L^{2}(\O_-^{\varepsilon})^4\oplus L^{2}(\O_{+})^4 &\to & L^{2}(\rr^3)^4\\
v=(v^{\varepsilon},v_{+}) &\mapsto & e_{\Omega^\varepsilon_{+-}}(v^{\varepsilon},v_+):=e_{\O_-^{\varepsilon}}v^{\varepsilon}+ e_{\Omega_+}v_+.
\end{array}
\end{aligned}
\end{equation}
   Let us recall for $z\in\rho(D_m)$, the lifting operators associated with boundary value problems \eqref{L1b}, \eqref{L1b'} and \eqref{L2b} are defined respectively, by 
   \begin{align*}
\begin{array}{rcl}
       E_{m}(z): P_- H^{1/2}(\S)^{4}&\to & H^{1}(\O_+)^{4}\\ 
         g_{+} &\mapsto & E_{m}(z)g_+:=\Phi^{z}_{m}(\Lambda^{z}_{+,m})^{-1}P_-,
       \end{array} \end{align*}
       \begin{align*}
       \begin{array}{rcl}
       E^{\varepsilon}_{m}(z): P_+ H^{1/2}(\S^{\varepsilon})^{4}&\to & H^{1}(\O_-^{\varepsilon})^{4}\\ 
         g^{\varepsilon} &\mapsto & E^{\varepsilon}_{m}(z)g^{\varepsilon}:=\Phi^{z,\varepsilon}_{m}(\Lambda^{z,\varepsilon}_{+,m})^{-1}P_+,
       \end{array}       \end{align*}
       $$\mathcal{E}^{\varepsilon}_{m+M}(z): P_+H^{1/2}(\S)^{4}\oplus P_-H^{1/2}(\S^{\varepsilon})^{4}\to  H^{1}(\mathcal{U}^{\varepsilon})^{4},$$
    with $\mathcal{E}^{\varepsilon}_{m+M}(z)(h_+,h^{\varepsilon}):=\Phi^{z}_{m+M}(\Lambda^{z}_{+,m+M})^{-1}P_+ h_+ + \Phi^{z,\varepsilon}_{m+M}(\Lambda^{z,\varepsilon}_{+,m+M})^{-1}P_- h^{\varepsilon}$.\\\\
In addition, we also recall the Poincaré-Steklov operators from Definition \ref{def}
  \begin{equation*}
   \begin{aligned}
\begin{array}{rcl}
       \mathscr{A}_{m}(z): P_- H^{1/2}(\S)^{4}&\to & P_+ H^{1/2}(\S)^{4}\\ 
         g_{+} &\mapsto & \mathscr{A}_{m}(z)g_{+}:=-P_+\beta(\Lambda^{z}_{+,m})^{-1}P_-g_{+},
       \end{array}\\
\begin{array}{rcl}
       \mathscr{A}^{\varepsilon}_{m}(z): P_+ H^{1/2}(\S^{\varepsilon})^{4}&\to & P_- H^{1/2}(\S^{\varepsilon})^{4}\\ 
         g^{\varepsilon}_- &\mapsto & \mathscr{A}^{\varepsilon}_{m}(z)g^{\varepsilon}:=-P_-\beta(\Lambda^{z,\varepsilon}_{+,m})^{-1}P_+g^{\varepsilon},
       \end{array}
       \end{aligned}
        \end{equation*}
      \begin{align*}\label{}
          \mathcal{A}^{\varepsilon}_{m+M}(z):  P_+H^{1/2}(\S)^{4}\oplus P_-H^{1/2}(\S^{\varepsilon})^{4}\to  P_-H^{1/2}(\S)^{4}\oplus P_+H^{1/2}(\S^{\varepsilon})^{4},\quad \text{with }
      \end{align*}
      $\quad\quad\quad\hspace{2cm} \mathcal{A}^{\varepsilon}_{m+M}(h_+,h^{\varepsilon}):=\big (-P_-\beta(\Lambda^{z}_{+,m+M})^{-1}P_+h_{+},-P_+\beta (\Lambda^{z,\varepsilon}_{+,m+M})^{-1}P_-h^{\varepsilon}\big)$.
\subsection{The Krein resolvent formula of $R^\varepsilon_M$}\label{kreinresolform}
Let $f\in\mathit{L}^{2}(\rr^3)^4$ and set 
	\begin{align*}
	 u^{\varepsilon}= r_{\mathcal{U}^{\varepsilon}}R^{\varepsilon}_{M}(z)f\quad \text{ and} \quad v=r_{\Omega^\varepsilon_{+-}}R^{\varepsilon}_{M}(z)f:=(v^{\varepsilon} \oplus v_+).
	\end{align*}
	Then $u^{\varepsilon}$ and $v$ satisfy the following system 
	\begin{equation*}\label{} 
	\left\{
	\begin{aligned}
(D_m- z) v_+ =f \quad &\text{ in }\Omega_+,\\
(D_m- z) v^{\varepsilon} =f \quad &\text{ in }\Omega_-^{\varepsilon},\\
(D_{m+M}-z)u^{\varepsilon} =f \quad &\text{ in }\mathcal{U}^{\varepsilon}, \\
 P_\pm t_{\S}v_+= P_{\pm} t_{\S}u^{\varepsilon}\quad& \text{ on } \S,\\
 P^\varepsilon_\mp t_{\S^{\varepsilon}}v^{\varepsilon}= P^\varepsilon_{\mp} t_{\S^{\varepsilon}}u^{\varepsilon}\quad & \text{ on } \S^{\varepsilon}.
	\end{aligned}
	\right.
	\end{equation*}  
Using  Lemma \ref{prop of C}, it is straightforward to check that the following resolvent formulas hold: 
\begin{align}\label{RMIT}
	R_{\mathrm{MIT}}^{\Omega^\varepsilon_-}(z)&= r_{\O_-^\varepsilon}(D_m -z)^{-1}e_{\O_-^\varepsilon} - \Phi_{m,-}^{z,\varepsilon}(\Lambda_{+,m}^{z,\varepsilon})^{-1}t_{\S^{\varepsilon}}(D_m -z)^{-1}e_{\O_-^\varepsilon},
		\end{align} 
	\begin{align*}\label{RMIT'}
	R_{\mathrm{MIT}}^{\Omega^\varepsilon_{+-}}(z)&= r_{\Omega^\varepsilon_{+-}}(D_m -z)^{-1}e_{\Omega^\varepsilon_{+-}} -r_{\Omega^\varepsilon_{+-}}e_{\Omega_+}\Phi_{m,+}^{z}(\Lambda_{+,m}^{z})^{-1}t_{\S}(D_m -z)^{-1}r_{\Omega_+}e_{\Omega^\varepsilon_{+-}} \\& -r_{\Omega^\varepsilon_{+-}}e_{\Omega_-^\varepsilon} \Phi_{m,-}^{z,\varepsilon}(\Lambda_{+,m}^{z,\varepsilon})^{-1}t_{\S^{\varepsilon}}(D_m -z)^{-1}r_{\Omega_-^\varepsilon}e_{\Omega^\varepsilon_{+-}},
		\end{align*} 
			\begin{align*}
	R_{\mathrm{MIT}}^{\mathcal{U}^\varepsilon}(z)&= r_{\mathcal{U}^{\varepsilon}}(D_{m}+M\beta-z)^{-1}e_{\mathcal{U}^{\varepsilon}} - \Phi_{m+M}^{z,\varepsilon}(\Lambda_{+,m+M}^{z,\varepsilon})^{-1}t_{\partial\mathcal{U}^{\varepsilon}}(D_{m}+M\beta -z)^{-1}e_{\mathcal{U}^{\varepsilon}}.
			\end{align*}  
In the whole following sections, and for simplicity, we'll use the following notation: 
\begin{align*}
    (\bullet,\bullet) := \mathrm{diag} (\bullet,\bullet)=\begin{pmatrix}
        \bullet &0\\0&\bullet
    \end{pmatrix}.
\end{align*}
Now, we set $\G_{\pm}:=P_{\pm}t_{\S} $ and $\G^{\varepsilon}_{\pm}:=P_{\pm}t_{\S^{\varepsilon}}.$ Since $E_{m}(z)$, $E^{\varepsilon}_{m}(z)$ and $\mathcal{E}^{\varepsilon}_{m+M}(z)$ gives the unique solution to the boundary value problem \eqref{L1b}, \eqref{L1b'} and  \eqref{L2b}, respectively, and the fact
   \begin{equation*} 
	\left\{
	\begin{aligned}
	\G_{-}R^{\Omega_+}_{\mathrm{MIT}}(z)r_{\O_{+}}f &=0,\quad	\G_+R^{\mathcal{U}^\varepsilon}_{M\mathrm{
MIT}}(z)r_{\mathcal{U}^{\varepsilon}}f &=0,\\
 \G^{\varepsilon}_{+}R_{\mathrm{MIT}}^{\Omega^\varepsilon_-}(z)r_{\O_-^{\varepsilon}}f &=0,\quad 	\G^{\varepsilon}_{-}R^{\mathcal{U}^\varepsilon}_{\mathrm{MIT}}(z)r_{\mathcal{U}^{\varepsilon}}f &=0.
	\end{aligned}
	\right.
	\end{equation*}
Then, if we let 
 \begin{equation*} 
	\left\{
	\begin{aligned}
	&\varphi =\G_{+}r_{\O_{+}}R_{M}^{\varepsilon}(z),\quad 	\varphi^{\varepsilon} =\G^{\varepsilon}_{-}r_{\O_-^{\varepsilon}}R_{M}^{\varepsilon}(z),\\&
 \psi =\G_{-}r_{\mathcal{U}^{\varepsilon}}R_{M}^{\varepsilon}(z),\quad	\psi^{\varepsilon} =\G^{\varepsilon}_{+}r_{\mathcal{U}^{\varepsilon}}R_{M}^{\varepsilon}(z),
	\end{aligned}
	\right.
	\end{equation*}
	it is easy to check that
	\begin{equation}\label{S2} 
	\left\{
	\begin{aligned}
	v_+ &=R^{\Omega_+}_{\mathrm{MIT}}(z)r_{\O_+}f +E_{m}(z)\psi,\\
 v^{\varepsilon} &=R_{\mathrm{MIT}}^{\Omega^\varepsilon_-}(z)r_{\O_-^{\varepsilon}}f +E^{{\varepsilon}}_{m}(z)\psi^\varepsilon,\\
	u^{\varepsilon} &=R^{\mathcal{U}^\varepsilon}_{\mathrm{MIT}}(z)r_{\mathcal{U}^{\varepsilon}}f +\mathcal{E}^{\varepsilon}_{m+M}(z)(\varphi,\varphi^{\varepsilon}).
	\end{aligned}
	\right.
	\end{equation}
Hence, to get an explicit formula for $R^\varepsilon_M (z)$ it remains to find the unknowns $(\varphi,\varphi^\varepsilon, \psi,\psi_\varepsilon)$. To do this, from \eqref{S2} we get
\begin{equation}\label{PPPP} 
	\left\{
	\begin{aligned}
	\varphi &= \G_+ v_+ = \G_+ R^{\Omega_+}_{\mathrm{MIT}}r_{\Omega_+}f + \mathscr{A}_{m}(z)\psi,\\
	\varphi^\varepsilon &= \G^\varepsilon_- v^\varepsilon = \G^\varepsilon_- R^{\Omega^\varepsilon_-}_{\mathrm{MIT}}r_{\Omega_\varepsilon}f + \mathscr{A}^\varepsilon_{m}(z)\psi^\varepsilon,\\
 \psi &= \G_- u^\varepsilon = \G_- R^{\mathcal{U}^\varepsilon}_{\mathrm{MIT}}(z)r_{\mathcal{U}^{\varepsilon}}f + \G_- \mathcal{E}^{\varepsilon}_{m+M}(z)(\varphi,\varphi^{\varepsilon}),\\
	 \psi^\varepsilon &= \G_+^\varepsilon u^\varepsilon = \G^\varepsilon_+ R^{\mathcal{U}^\varepsilon}_{\mathrm{MIT}}(z)r_{\mathcal{U}^{\varepsilon}}f + \G_+ \mathcal{E}^{\varepsilon}_{m+M}(z)(\varphi,\varphi^{\varepsilon}).
	\end{aligned}
	\right.
	\end{equation}
Using the restriction map $r_{\bullet}$ and the extension map $e_{\bullet}$ given in \eqref{re}, we get 
\begin{equation*}\label{S3} 
	\left\{
	\begin{aligned}
	v &=e_{\Omega^{\varepsilon}_{+-}}\big(R^{\Omega_+}_{\mathrm{MIT}}(z),R_{\mathrm{MIT}}^{\Omega_-^\varepsilon}(z)\big)r_{\Omega^\varepsilon_{+-}}f + e_{\Omega^\varepsilon_{+-}}\Big(E_{m}(z)P_-,E^{{\varepsilon}}_{m}(z)P_+\Big)(\G_-,\G^{\varepsilon}_+)r_{\mathcal{U}^{\varepsilon}}R_{M}^{\varepsilon}(z)f,\\
	u^{\varepsilon} &=R^{\mathcal{U}^\varepsilon}_{\mathrm{MIT}}(z)r_{\mathcal{U}^{\varepsilon}}f +\mathcal{E}^{\varepsilon}_{m+M}(z)(P_+,P_-)\big(\G_+r_{\O_+},\G^{\varepsilon}_-r_{\O_-^\varepsilon}\big)R_{M}^{\varepsilon}(z)f.
	\end{aligned}
	\right.
	\end{equation*}
	Thus, we obtain 
 \begin{align}\label{preformula}
	\begin{split}
	R^{\varepsilon}_M  (z) & =e_{\Omega_+}R^{\Omega_+}_{\mathrm{MIT}}r_{\Omega_+} + e_{\Omega_-^\varepsilon}R^{\Omega^\varepsilon_-}_{\mathrm{MIT}}r_{\Omega_-^\varepsilon} + e_{\mathcal{U}^{\varepsilon}}R^{\mathcal{U}^\varepsilon}_{\mathrm{MIT}}(z)r_{\mathcal{U}^{\varepsilon}}  \\
&\hspace{1cm} + \bigg(e_{\Omega^\varepsilon_{+-}}\Big(E_{m}(z)P_-,E^{{\varepsilon}}_{m}(z)P_+\Big)(\G_-,\G^{\varepsilon}_+)r_{\mathcal{U}^{\varepsilon}}+  e_{\mathcal{U}^{\varepsilon}} \mathcal{E}^{\varepsilon}_{m+M}(z)\big(\G_+ r_{\Omega_+},\G^{\varepsilon}_- r_{\O_-^{\varepsilon}}\big)\bigg) R^{\varepsilon}_M(z)  \\
	&= e_{\Omega^\varepsilon_{+-}} R^{\Omega^{\varepsilon}_{+-}}_{\mathrm{MIT}}(z) r_{\Omega^\varepsilon_{+-}}+ e_{\mathcal{U}^{\varepsilon}}R^{\mathcal{U}^\varepsilon}_{\mathrm{MIT}}(z)r_{\mathcal{U}^{\varepsilon}}+E^{\varepsilon}_M(z)\G^{\varepsilon} R^{\varepsilon}_{M}(z),
	\end{split}
	\end{align} 
	with
	$R_{\mathrm{MIT}}^{\Omega^\varepsilon_{+-}}(z)$ as in \eqref{MRMIT}.\\ \\  
 Here $\G^{\varepsilon}$ and $E^{\varepsilon}_{M}(z)$ are defined as follows:
\begin{align*}
    \begin{array}{rcl}
       \G^{\varepsilon}: H^{1}(\Omega_+)^4\oplus H^{1}(\Omega_-^\varepsilon)^4\oplus H^{1}(\mathcal{U}^\varepsilon)^4 & \to & H^{1/2}(\Sigma)^4\oplus H^{1/2}(\Sigma^\varepsilon)^4\oplus H^{1/2}(\Sigma)^4\oplus H^{1/2}(\Sigma^\varepsilon)^4\\ 
         (r_{\Omega_+},\, r_{\Omega_-^\varepsilon}, \, r_{\mathcal{U}^\varepsilon}) &\mapsto & (\G_+r_{\O_+}\,\, \G^{\varepsilon}_-r_{\O_-^{\varepsilon}}\,\, \G_-r_{\mathcal{U}^{\varepsilon}}\, \,\G^{\varepsilon}_+r_{\mathcal{U}^{\varepsilon}})^{t},
       \end{array}
\end{align*}
and  $E^{\varepsilon}_{M}(z)=e_{\Omega^\varepsilon_{+-}} E^{\Omega^\varepsilon_{+-}}_m(z)+  e_{\mathcal{U}^{\varepsilon}} \mathcal{E}^{\varepsilon}_{m+M}(z)(P_+,P_-)$, with $E_m^{\Omega^\varepsilon_{+-}}(z) = r_{\Omega^\varepsilon_{+-}}e_{\Omega_+} E_m (z) P_- + r_{\Omega^\varepsilon_{+-}} e_{\Omega_-^\varepsilon} E_m^{\varepsilon}(z)P_+$ can be read as the following matrix:
\begin{align}\label{liftom+-}
    \begin{array}{rcl}
     E_m^{\Omega^\varepsilon_{+-}}:P_- H^{1/2}(\Sigma)^4\oplus P_+ H^{1/2}(\S^\varepsilon)^4& \to & H^{1}(\Omega_+)^4\oplus H^{1}(\Omega_-^\varepsilon)^4\\ 
         (\psi,\, \psi_{\varepsilon}) &\mapsto & (E_m P_- \psi,\, E_m^{\varepsilon}P_+ \psi_\varepsilon)\equiv \begin{pmatrix}
	E_m P_-  & 0\\
	0 & E_m^\varepsilon P_+
	\end{pmatrix}\begin{pmatrix}
	\psi\\\psi_\varepsilon
	\end{pmatrix}.
       \end{array}
\end{align}
Now, applying $\G^{\varepsilon}$ to the identity \eqref{preformula}, it yields   
\begin{align}\label{LI}
    \G^{\varepsilon} R^{\varepsilon}_{\mathrm{MIT}}(z)=\Big(\mathbb{I}- \big(\mathscr{A}_{m}(z)P_-,\mathscr{A}^{\varepsilon}_{m}(z)P_+\big)- \mathcal{A}^{\varepsilon}_{m+M}(z)(P_+,P_-)\Big) \G^{\varepsilon}R_{M}^{\varepsilon}(z):=\Upsilon_{M}^{\varepsilon}(z)\G^{\varepsilon}R_{M}^{\varepsilon}(z),
\end{align}
with $R^{\varepsilon}_{\mathrm{MIT}}(z):=e_{\Omega^\varepsilon_{+-}}R_{\mathrm{MIT}}^{\Omega^\varepsilon_{+-}}(z) + e_{\mathcal{U}^{\varepsilon}}R^{\mathcal{U}^\varepsilon}_{\mathrm{MIT}}(z).$ Similarly, we mention that $\big[\mathscr{A}_{m}(z),\mathscr{A}^{\varepsilon}_{m}(z)\big]$ means the sum of both terms $\mathscr{A}_{m}$, $\mathscr{A}^{\varepsilon}_{m}$ and can be read as the following matrix 
\begin{align}\label{MAm}
    \begin{array}{rcl}
   \mathscr{A}_m^{\Omega^\varepsilon_{+-}}:=\big(\mathscr{A}_m, \mathscr{A}_m^{\varepsilon} \big):P_- H^{1/2}(\Sigma)^4\oplus P_+ H^{1/2}(\S^\varepsilon)^4& \to & P_+ H^{1/2}(\Sigma)^4\oplus H^{1/2}(\Sigma^\varepsilon)^4\\ 
         (\psi,\, \psi_{\varepsilon}) &\mapsto & \big(\mathscr{A}_m, \mathscr{A}_m^{\varepsilon} \big)(\psi,\, \psi_{\varepsilon}) = \begin{pmatrix}
	\mathscr{A}_m P_- & 0\\
	0 & \mathscr{A}_m^\varepsilon P_+ 
	\end{pmatrix}\begin{pmatrix}
	\psi\\\psi_\varepsilon
	\end{pmatrix}.
       \end{array}
\end{align}
Using the formula of $\mathcal{A}_{m+M}^\varepsilon$,  the term $(\G_-,\G_+^\varepsilon)\mathcal{E}^{\varepsilon}_{m+M}(z)$ is identified with $(P_-\mathcal{A}_{m+M}^\varepsilon, P_+\mathcal{A}_{m+M}^\varepsilon) = \\ (P_-, P_+)\mathcal{A}_{m+M}^\varepsilon(z)\equiv (P_-, 0) \mathcal{A}_{m+M}^\varepsilon(z) +(0, P_+) \mathcal{A}_{m+M}^\varepsilon(z) .$
\\ \\
Now, applying also $\Big(\mathbb{I}+ \mathscr{A}_m^{\Omega^\varepsilon_{+-}}(z) + (P_-,P_+)\mathcal{A}^{\varepsilon}_{m+M}(z)\Big)$ to the identity \eqref{LI} we get
\begin{align*}
\G^{\varepsilon}R_{M}^{\varepsilon}(z)=\Xi^{\varepsilon}_{M}(z)\Big(\mathbb{I}+\mathscr{A}_m^{\Omega^\varepsilon_{+-}}(z) + (P_-, P_+) \mathcal{A}_{m+M}^\varepsilon(z)\Big)\G^{\varepsilon} R^{\varepsilon}_{\mathrm{MIT}}(z),
\end{align*}
with $\Xi_{M}^{\varepsilon}(z): H^{1/2}(\Sigma)^4 \oplus H^{1/2}(\Sigma^{\varepsilon})^4 \rightarrow  H^{1/2}(\Sigma)^4 \oplus H^{1/2}(\Sigma^{\varepsilon})^4$ the following quantity
\begin{align}\label{XIME}
    \Xi_{M}^{\varepsilon}(z):=\bigg(\mathbb{I}_8-\mathscr{A}_m^{\Omega^\varepsilon_{+-}}(z) (P_-, P_+)\mathcal{A}^{\varepsilon}_{m+M}(z)- \mathcal{A}^{\varepsilon}_{m+M}(z)(P_+, P_-)\mathscr{A}_m^{\Omega^\varepsilon_{+-}}(z) \bigg)^{-1}.
\end{align}
From which it follows that,
\begin{align}\label{kfullreso}
    R^{\varepsilon}_{M}(z)=R_{\mathrm{MIT}}^{\varepsilon}(z)+E^{\varepsilon}_{M}(z)[\Upsilon_{M}^{\varepsilon}(z)]^{-1} \G^{\varepsilon}R_{\mathrm{MIT}}^{\varepsilon}(z),
    \end{align}
    with \begin{align*}
  [\Upsilon^{\varepsilon}_M]^{-1}(z)= \Xi_{M}^{\varepsilon}(z)\Big(\mathbb{I}+\mathscr{A}_m^{\Omega^\varepsilon_{+-}}(z) + (P_-, P_+) \mathcal{A}_{m+M}^\varepsilon(z)\Big). 
  \end{align*} 
  \begin{remark}\label{remark6.1}
      The identity \eqref{LI} has the following matrix form
      \begin{align*}
   \begin{pmatrix} \G_+ r_{\Omega_+} R^\varepsilon_M \\ \G_- r_{\Omega_-^\varepsilon} R^\varepsilon_M \\\G_- r_{\mathcal{U}^\varepsilon} R^\varepsilon_M \\ \G^{\varepsilon}_+ r_{\mathcal{U}^\varepsilon}R^\varepsilon_M
	\end{pmatrix}=   \begin{pmatrix}
	 \G_+ R^{\Omega_+}_{\mathrm{MIT}}r_{\Omega_+}\\ \G_-^\varepsilon R^{\Omega^\varepsilon_-}_{\mathrm{MIT}}r_{\Omega_-^\varepsilon} \\\G_- R^{\mathcal{U}^\varepsilon}_{\mathrm{MIT}}r_{\mathcal{U}^\varepsilon} \\ \G^{\varepsilon}_+ R^{\mathcal{U}^\varepsilon}_{\mathrm{MIT}}r_{\mathcal{U}^\varepsilon}
	\end{pmatrix} + \begin{pmatrix} 0 & 0 & \mathscr{A}_m P_- & 0 \\ 0 & 0  & 0 & \mathscr{A}_m^\varepsilon P_+  \\ \mathcal{A}_{m+M}^\varepsilon(P_+, P_-) & \mathcal{A}_{m+M}^\varepsilon(P_+, P_-) & 0 & 0 \\ \mathcal{A}_{m+M}^\varepsilon(P_+, P_-) & \mathcal{A}_{m+M}^\varepsilon(P_+, P_-)  & 0 &  0
	\end{pmatrix}\begin{pmatrix} \G_+ r_{\Omega_+} R^\varepsilon_M \\ \G^\varepsilon_- r_{\Omega_-^\varepsilon} R^\varepsilon_M \\\G_- r_{\mathcal{U}^\varepsilon} R^\varepsilon_M \\ \G^{\varepsilon}_+ r_{\mathcal{U}^\varepsilon}R^\varepsilon_M
	\end{pmatrix}.
  \end{align*}
  Moreover, if we note by $\G^\varepsilon_{+-} = (\G_+r_{\O_+}\,\, \G^{\varepsilon}_-r_{\O_-^{\varepsilon}})^{t}$ and $\G^\varepsilon_{-+}= ( \G_-r_{\mathcal{U}^{\varepsilon}}\, \,\G^{\varepsilon}_+r_{\mathcal{U}^{\varepsilon}})^{t}$. Then, using the quantities of \eqref{PPPP}, we remark that the Krein resolvent formula \ref{kfullreso} can be also written in the following matrix
    \begin{align*}
   \begin{pmatrix} r_{\Omega^\varepsilon
_{+-}} R^\varepsilon_M \\  r_{\mathcal{U}^\varepsilon}R^\varepsilon_M
	\end{pmatrix}=   \begin{pmatrix}
	 R^{\Omega^\varepsilon_{+-}}_{\mathrm{MIT}}r_{\Omega^\varepsilon_{+-}} \\ R^{\mathcal{U}^\varepsilon}_{\mathrm{MIT}}r_{\mathcal{U}^\varepsilon} 
	\end{pmatrix} + \begin{pmatrix} E^{\Omega^\varepsilon_{+-}}_m \Xi_M^{\varepsilon,-+}& 0  \\ 0 & E_{m+M}^\varepsilon \Xi_M^{\varepsilon,+-}
	\end{pmatrix}\begin{pmatrix}\mathcal{A}_{m+M}^\varepsilon  & \mathbb{I}_4  \\  \mathbb{I}_4 &\mathscr{A}_m^{\Omega^\varepsilon_{+-}}\end{pmatrix} \begin{pmatrix}
	 \G^{\varepsilon}_{+-}R^{\Omega^\varepsilon_{+-}}_{\mathrm{MIT}}r_{\Omega^\varepsilon_{+-}}\\ \G^{\varepsilon}_{-+}R^{\mathcal{U}^\varepsilon}_{\mathrm{MIT}}r_{\mathcal{U}^\varepsilon} 
\end{pmatrix},
\end{align*}
  where $\mathscr{A}_m^{\Omega^\varepsilon_{+-}} $ is the matrix in \eqref{MAm} and $\Xi_M^{\varepsilon,\pm\mp}$ are given in the following corollary.
  \end{remark}
\begin{corollary}\label{une autre est} Consider the operator $\Xi_{M}^{\varepsilon}(z)$  given in \eqref{XIME}. Then, there is $M_0>0$ such that for every  $M>M_0$, $h\equiv\varepsilon=1/M$ and for all $ z\in\rho(D_{\mathrm{MIT}}^{\Omega^\varepsilon_{+-}})\cap\rho(\mathfrak{D}^{\varepsilon}_M)$,  the operator $\Xi^{\varepsilon}_M(z)$ is everywhere defined and uniformly bounded with respect to $M$. Moreover, the operators $\Xi^{\varepsilon,+-}_{M}(z)$ and $\Xi^{\varepsilon,-+}_{M}(z)$ defined by 
\begin{align*}
\begin{array}{rcl}
       \Xi^{\varepsilon,+-}_{M}(z):  P_+H^{s}(\S)^4\oplus P_-H^{s}(\S^{\varepsilon})^{4}&\to & P_+H^{s}(\S)^{4}\oplus P_- H^{s}(\S^{\varepsilon})^{4},
       \end{array}
   \end{align*}
   \begin{align*}
       \begin{array}{rcl}
       \Xi^{\varepsilon,-+}_{M}(z):  P_-H^{s}(\S)^{4}\oplus P_+ H^{s}(\S^{\varepsilon})^{4} & \to & P_-H^{s}(\S)^4\oplus P_+H^{s}(\S^{\varepsilon})^{4},
       \end{array}
   \end{align*}
 which have  the following formula
   \begin{align*}
\Xi^{\varepsilon,+-}_{M}(z)&=\Big(\mathbb{I} -\mathscr{A}_m^{\Omega^\varepsilon_{+-}}(z)(z)(P_-, P_+)\mathcal{A}^{\varepsilon}_{m+M}(z)(P_+, P_-)\Big)^{-1},\\
\Xi^{\varepsilon,-+}_{M}(z) &=\Big(\mathbb{I}-(P_-, P_+)\mathcal{A}^{\varepsilon}_{m+M}(z)(P_+, P_-)\mathscr{A}_m^{\Omega^\varepsilon_{+-}}(z)\Big)^{-1}
 \end{align*}
 are bounded for any $s\in\mathbb{R}$, and it holds that
 \begin{align}\label{XIM+-1}
     ||\Xi^{\varepsilon,\pm\mp}_{M}(z) ||_{P_\pm H^{-1/2}(\S)^4\oplus P_\mp H^{-1/2}(\S^{\varepsilon})^{4}\rightarrow P_\pm H^{-1/2}(\S)^4\oplus P_\mp H^{-1/2}(\S^{\varepsilon})^{4}} \lesssim 1,
 \end{align}
 uniformly with respect to $M>M_0$.

 Moreover, the Poincaré-Steklov $\mathcal{A}_{m+M}^\varepsilon$, satisfies the following estimate
 \begin{align}\label{ESTAME}
     ||\mathcal{A}_{m+M}^\varepsilon||_{P_{+} H^{1/2}(\Sigma)^4\oplus P_{-} H^{1/2}(\Sigma^\varepsilon)^4 \rightarrow P_{-} H^{-1/2}(\Sigma)^4\oplus P_{+} H^{-1/2}(\Sigma^\varepsilon)^4} \lesssim M^{-1}.
 \end{align}
\end{corollary}
\textbf{Proof.} 
Set $\kappa=m+M$ and $h=\kappa^{-1}$. The proof of this corollary follows a similar argument as in \cite[Proposition 6.1]{BBZ}. It is based on the pseudodifferential properties of the Poincaré-Steklov operators $\mathscr{A}^{\varepsilon}_m$ and $\mathcal{A}_{\kappa}^{\varepsilon}$. Since $\mathscr{A}_m$ (resp. $\mathscr{A}^{\varepsilon}_m$) are a pseudodifferential operators of order $0$, see Remark \ref{remark5.2} (resp. Corollary \ref{corollsry5.1}), we can consider it as an $h$-pseudodifferential operator of order $0$ whose principal symbol is given by:
\begin{align*}
\mathscr{P}_{h,\mathscr{A}_m}(x_\Sigma,\xi)&= \dfrac{S\cdot(\xi\wedge\nu(x_\Sigma))}{\big|\xi\wedge\nu(x_\Sigma)\big|}P_{-},\quad x_\Sigma\in\Sigma,\\
\mathscr{P}_{h,\mathscr{A}_m^{\varepsilon}}(x,\xi)&= - \big(1-\varepsilon W(x_\Sigma)\big)^{-1}\dfrac{S\cdot(\xi\wedge\nu(p^{-1}(x)))}{\big|\big(1-\varepsilon W(x_\Sigma)\big)^{-1}\big|\Big|\xi\wedge\nu(p^{-1}(x))\Big|}P_{+},\, \Sigma^\varepsilon\ni x=p(x_\Sigma)=x_\Sigma + \varepsilon \nu(x_\Sigma),
\end{align*}
where $S$ is the spin angular momentum given in Lemma \ref{properties3.1}, $\xi\in\mathbb{R}^{2}$ can be identify with $\bar{\xi}=(\xi_1,\xi_2,0)^{t}\in \mathbb{R}^{3}$, $p$ is the diffeomorphism from Remark \ref{remark3.1}, and for $x = \varphi(\tilde{x})$ stands for $\nu^{\varphi}(\tilde{x})$. On the other hand, Proposition \ref{prop4.6} follows that $\mathcal{A}_{\kappa}^{\varepsilon}$ is $h$-pseudodifferential operator of order $0$ has the following principal symbol 
\begin{align*}
     \mathscr{P}_{h,\mathcal{A}_{\kappa}^{\varepsilon}}(x,\xi)= \big(1-\varepsilon W(x_\Sigma)\big)^{-1}\dfrac{S\cdot(\xi\wedge\nu(p^{-1}(x)))}{\sqrt{\Big(\big(1-\varepsilon W(x_\Sigma)\big)^{-1}\xi\wedge\nu(p^{-1}(x))\Big)^{2}+1} + 1}\begin{pmatrix}
        -P_+&0\\0& P_-
     \end{pmatrix}.
\end{align*}
Consequently, the symbol calculus yields for all $h<h_0$ that $$ \bigg(\mathbb{I}_8-\mathscr{A}_m^{\Omega^\varepsilon_{+-}}(z)(P_-, P_+) \mathcal{A}^{\varepsilon}_{\kappa}(z)- \mathcal{A}^{\varepsilon}_{\kappa}(z)(P_+, P_-)\mathscr{A}_m^{\Omega^\varepsilon_{+-}}(z)\bigg)$$
is a $\kappa^{-1}$-pseudodifferential operator of order $0$.\\\\
Now, using the principal symbols of $\mathscr{A}_m$, $\mathscr{A}_m^\varepsilon$, the principal symbol of $\mathscr{A}_m^{\Omega^\varepsilon_{+-}}$ can be written as the following:
\begin{align*}
    \mathscr{P}_{h,\,\mathscr{A}_m^{\Omega^\varepsilon_{+-}}}(x_{\Sigma},\xi)&=\begin{pmatrix}
	\mathscr{P}_{h,\mathscr{A}_m}(x_\Sigma,\xi) & 0\\
	0 & \mathscr{P}_{h,\mathscr{A}_m^{\varepsilon}}(p(x_\Sigma),\xi) 
	\end{pmatrix}\\&=\dfrac{S\cdot(\xi\wedge\nu(x_\Sigma))}{\big|\xi\wedge\nu(x_\Sigma)\big|}\begin{pmatrix}
	P_- & 0\\
	0 & -\dfrac{\big(1-\varepsilon W(x_\Sigma)\big)^{-1}}{\big| \big(1-\varepsilon W(x_\Sigma)\big)^{-1}\big|}P_+ 
	\end{pmatrix}.
\end{align*}
Using Lemma \ref{properties3.1}, we obtain 
\begin{align*}
\mathscr{P}_{h,\,\mathscr{A}_m^{\Omega^\varepsilon_{+-}}}(x_{\Sigma},\xi) \mathscr{P}_{h,\mathcal{A}_{\kappa}^{\varepsilon}}(x,\xi)=&\\&\hspace{-3cm} -\dfrac{\big(1-\varepsilon W(x_\Sigma)\big)^{-1}\big|\xi\wedge\nu(x_\Sigma)\big|}{\sqrt{\Big(\big(1-\varepsilon W(x_\Sigma)\big)^{-1}\xi\wedge\nu(p^{-1}(x))\Big)^{2}+1} + 1} \begin{pmatrix}
    P_+ &0\\0& \dfrac{\big(1-\varepsilon W(x_\Sigma)\big)^{-1}}{\big| \big(1-\varepsilon W(x_\Sigma)\big)^{-1}\big|}P_- 
\end{pmatrix}.
\end{align*} 
Then, it yields
\begin{align*}
    \mathbb{I}_8 -\mathscr{P}_{h,\,\mathscr{A}_m^{\Omega^\varepsilon_{+-}}}(x_{\Sigma},\xi) \mathscr{P}_{h,\mathcal{A}_{\kappa}^{\varepsilon}}(x,\xi) - \mathscr{P}_{h,\mathcal{A}_{\kappa}^{\varepsilon}}(x,\xi)\mathscr{P}_{h,\,\mathscr{A}_m^{\Omega^\varepsilon_{+-}}}(x_{\Sigma},\xi)=&\\&\hspace{-7cm}\mathbb{I}_8  +  \dfrac{\big(1-\varepsilon W(x_\Sigma)\big)^{-1}\big|\xi\wedge\nu(x_\Sigma)\big|}{\sqrt{\Big(\big(1-\varepsilon W(x_\Sigma)\big)^{-1}\xi\wedge\nu(p^{-1}(x))\Big)^{2}+1} + 1}\begin{pmatrix}
    \mathbb{I}_4&0 \\0&\dfrac{\big(1-\varepsilon W(x_\Sigma)\big)^{-1}}{\big| \big(1-\varepsilon W(x_\Sigma)\big)^{-1}\big|}\mathbb{I}_4
\end{pmatrix}\gtrsim 1.
\end{align*}
Thus, $\Xi_{M}^\varepsilon$ is a zero-order pseudoddiferential operator. \\

Thanks to the following relationship: $\Xi^{\varepsilon,\pm\mp}_{M}(z) = (P_\pm, P_\mp) \Xi_{M}^{\varepsilon}(z)(P_\pm,P_\mp),$ it yields the same properties for $\Xi^{\varepsilon,\pm\mp}_{M}(z)$ and therefore \eqref{XIM+-1} is established.\\

Regarding the estimate of $\mathcal{A}_{\kappa}^{\varepsilon}$, exploits also the Calderón-Vaillancourt theorem which shows that for any operator in $h\,Op^h \mathcal{S}^0(\partial\mathcal{U}^\varepsilon)$ is uniformly bounded by $O(h)$, with respect to $h = \kappa^{-1} \in (0,1)$, from $H^{1/2}(\partial\mathcal{U}^\varepsilon)^4$ into $H^{1/2}(\partial\mathcal{U}^\varepsilon)^4 \hookrightarrow H^{-1/2}(\partial\mathcal{U}^\varepsilon)^4$, see \eqref{Calderon-VaillancourtHS}. Thus,
 $$ \Big|\Big|\mathcal{A}_{\kappa}^{\varepsilon} -  \dfrac{S\cdot(\nabla_{\partial\mathcal{U}^\varepsilon}\wedge\nu(p^{-1}(x)))}{\sqrt{-\kappa^{-2}\Delta_{\partial\mathcal{U}^\varepsilon
 }+\mathbb{I}} + \mathbb{I}}(P_+, P_-)\Big|\Big|_{H^{1/2}(\partial\mathcal{U}^\varepsilon)^4 \rightarrow H^{-1/2}(\partial\mathcal{U}^\varepsilon)^4}\lesssim \kappa^{-1},$$
 uniformly with respect to $\kappa$ big enough and $\varepsilon \in (0, \varepsilon_0)$.  Then we conclude the proof of the estimate by using that  $\Big(\sqrt{-\kappa^{-2}\Delta_{\partial\mathcal{U}^\varepsilon
 }+\mathbb{I}} + \mathbb{I}\Big)^{-1}$ is uniformly bounded from  $H^{1/2}(\partial\mathcal{U}^\varepsilon)^4$ into itself and $(\nabla_{\partial\mathcal{U}^\varepsilon}\wedge \nu(p^{-1}(x)))$ is  uniformly bounded from  $H^{1/2}(\partial\mathcal{U}^\varepsilon)^4$ into $H^{-1/2}(\partial\mathcal{U}^\varepsilon)^4$.
    \begin{remark}
Let $E^{\Omega^\varepsilon_{+-}}_{m}$ from \eqref{liftom+-}.	Thanks to \cite[Proposition 4.1 $(\mathrm{ii})$]{BBZ}, we have that 
	\begin{align*}
	\left(  E^{\Omega^\varepsilon_{+-}}_{m}(z)\right)^{\ast}=-\beta 
(\G_{+-}^\varepsilon)^t R^{\Omega^\varepsilon_{+-}}_\mathrm{MIT} (\overline{z}) \quad\mathrm{and}\quad  \left(  \mathcal{E}^{\varepsilon}_{m+M}(z)\right)^{\ast}=-\beta  (\G_{-+}^\varepsilon)^t\begin{pmatrix}
	   R^{\mathcal{U}^\varepsilon}_{\mathrm{MIT}} (\overline{z}) \\ R^{\mathcal{U}^\varepsilon}_{\mathrm{MIT}} (\overline{z})
	\end{pmatrix},
	\end{align*}
for any $ z\in\rho(D^{\Omega^\varepsilon_{+-}}_{\mathrm{MIT}})\cap\rho(\mathfrak{D}_{M}^{\varepsilon})$. 
	\end{remark}
\section{Resolvent convergence to the Dirac operator with Lorentz scalar.}\label{convergences}
In this section, we gather the necessary elements to prove the main result of this work. The components of the proof for the main theorem (\emph{i.e.}, Theorem \ref{maintheo}) are dedicated to examining the convergence of the terms present in the resolvent formula \eqref{preformula}. It is important to note that this resolvent formula includes certain terms independent of $M$ and $\varepsilon$, namely $E_m, \mathscr{A}_m,$ and $R^{\Omega_+}_{\mathrm{MIT}}r_{\Omega_+}$, which remain fixed and act within $\Omega_+$. Consequently, our focus shifts to examining the convergence of terms dependent on $\varepsilon$ but independent of $M$, namely $R^{\Omega^\varepsilon_-}_{\mathrm{MIT}}r_{\Omega^{\varepsilon_-}}$ and $E_{m}^{\varepsilon}$ (see, Proposition \ref{coro1.1}). Subsequently, we will proceed to estimate the remaining terms in relation to $M$ and $\varepsilon$ (see, Proposition \ref{Prop1.1}).  
\begin{proposition}\label{RMITCONV} Let $\varepsilon_0 >0$ be small enough, and let $z\in\mathbb{C}\backslash\mathbb{R}$. We set $\Omega_-:=\mathbb{R}^3\setminus\overline{\Omega_+}$ the exterior fixed domain and by $\Sigma=\partial\Omega_- = \partial\Omega_+$ its boundary. We denote by $R_{\mathrm{MIT}}^{\Omega_-}$ the resolvent of the fixed MIT bag operator, which we denote by $D_{\mathrm{MIT}}^{\Omega_-}$, acts in $\Omega_-$. Then, for any  $\varepsilon\in (0, \varepsilon_0)$ the following holds:
\begin{align}\label{approresol}
\Big|\Big| e_{\Omega^{\varepsilon}_-}R^{\Omega^\varepsilon_-}_{\mathrm{MIT}}(z)r_{\Omega^{\varepsilon}_-} - e_{\Omega_-}R_{\mathrm{MIT}}^{\Omega_-}(z)r_{\Omega_-}\Big|\Big|_{\mathit{L}^2(\rr^3)^4\rightarrow \mathit{L}^2(\rr^3)^4}& = \mathcal{O}(\varepsilon).
\end{align} 
\end{proposition}
\textbf{Proof.} The Krein formula for the resolvent $R^{\Omega^\varepsilon_-}_{\mathrm{MIT}}$ (from equality \eqref{RMIT} )
\begin{align*}
    e_{\Omega^{\varepsilon}_-}R^{\Omega^\varepsilon_-}_{\mathrm{MIT}}(z)r_{\Omega^{\varepsilon}_-} &= (D_m - z)^{-1} - e_{\Omega^{\varepsilon}_-} \Phi^{z,\varepsilon}_{m,-} (\Lambda^{z,\varepsilon}_{+,m})^{-1} t_{\Sigma^\varepsilon} (D_m - z)^{-1},\\
    e_{\Omega_-}R^{\Omega_-}_{\mathrm{MIT}}(z)r_{\Omega_-} &= (D_m - z)^{-1} - e_{\Omega_-} \Phi^{z}_{m,-} (\Lambda^{z}_{m,+})^{-1} t_{\Sigma} (D_m - z)^{-1}
\end{align*}
yield that 
\begin{equation}\label{estimresol}
\begin{aligned}
&\Big|\Big| e_{\Omega^{\varepsilon}_-}R^{\Omega^\varepsilon_-}_{\mathrm{MIT}}(z)r_{\Omega^{\varepsilon}_-} - e_{\Omega_-}R_{\mathrm{MIT}}^{\Omega_-}(z)r_{\Omega_-}\Big|\Big|_{\mathit{L}^2(\rr^3)^4\rightarrow \mathit{L}^2(\rr^3)^4} \\& = \Big|\Big|e_{\Omega_-} \Phi^{z}_{m,-} (\Lambda^{z}_{+,m})^{-1} t_{\Sigma} (D_m - z)^{-1} - e_{\Omega^{\varepsilon}_-} \Phi^{z,\varepsilon}_{m,-} (\Lambda^{z,\varepsilon}_{+,m})^{-1} t_{\Sigma^\varepsilon} (D_m - z)^{-1} \Big|\Big|_{\mathit{L}^2(\rr^3)^4\rightarrow \mathit{L}^2(\rr^3)^4}\\& \leq \Big|\Big|e_{\Omega_-} \Phi^{z}_{m,-} (\Lambda^{z}_{+,m})^{-1} t_{\Sigma} - e_{\Omega^{\varepsilon}_-} \Phi^{z,\varepsilon}_{m,-} (\Lambda^{z,\varepsilon}_{+,m})^{-1} t_{\Sigma^\varepsilon}  \Big|\Big|_{\mathit{H}^1(\rr^3)^4\rightarrow \mathit{L}^2(\rr^3)^4} \Big|\Big| (D_m - z)^{-1} \Big|\Big|_{\mathit{L}^2(\rr^3)^4\rightarrow \mathit{H}^1(\rr^3)^4} \\& \lesssim \Big|\Big|e_{\Omega_-} \Phi^{z}_{m,-} (\Lambda^{z}_{+,m})^{-1} t_{\Sigma} - e_{\Omega^{\varepsilon}_-} \Phi^{z,\varepsilon}_{m,-} (\Lambda^{z,\varepsilon}_{+,m})^{-1} t_{\Sigma^\varepsilon}  \Big|\Big|_{\mathit{H}^1(\rr^3)^4\rightarrow \mathit{L}^2(\rr^3)^4}
\end{aligned} 
\end{equation}
since $(D_m - z)^{-1}$ is bounded from $\mathit{L}^2(\rr^3)^4$ into $\mathit{H}^1(\rr^3)^4$.\\\\
To obtain a rigorous estimate of the right-hand side of \eqref{estimresol}, we'll use the unitary transformation $\mathcal{T}_{\varepsilon}$ from Definition \ref{Transformationoperator} and the explicit formula for $\Lambda^{z}_{+,m}$ (resp. $\Lambda^{z,\varepsilon}_{+,m}$). Let $f,g\in \mathit{L}^2(\rr^3)^4.$ Since $t_{\S}(D_m-z)^{-1}=(\Phi^{\bar{z}}_{m})^{\ast}$  \big(resp. $t_{\S^\varepsilon}(D_m-z)^{-1}=(\Phi^{\bar{z},\varepsilon}_{m})^{\ast}$ \big) by duality and interpolation arguments, we get that
\begin{align*}
  & \Big| \Big\langle \big[e_{\Omega_-} \Phi^{z}_{m,-} \Big(\beta/2 + \mathscr{C}^{z}_m\Big )^{-1} t_{\Sigma} - e_{\Omega^{\varepsilon}_-} \Phi^{z,\varepsilon}_{m,-} \Big(\beta/2 + \mathscr{C}^{z, \varepsilon}_m\Big )^{-1} t_{\Sigma^\varepsilon} \big]f, g \Big\rangle_{\mathit{L}^2(\rr^3)^4, \mathit{L}^2(\rr^3)^4}  \Big|\\& = \Big| \Big\langle e_{\Omega_-} \Phi^{z}_{m,-} \Big(\beta/2 + \mathscr{C}^{z}_m\Big )^{-1} t_{\Sigma}f, g \Big\rangle_{\mathit{L}^2(\rr^3)^4,\mathit{L}^2(\rr^3)^4} -  \Big\langle  e_{\Omega^{\varepsilon}_-} \Phi^{z,\varepsilon}_{m,-} \Big(\beta/2 + \mathscr{C}^{z, \varepsilon}_m\Big )^{-1} t_{\Sigma^\varepsilon}f, g \Big\rangle_{ \mathit{L}^2(\rr^3)^4,\mathit{L}^2(\rr^3)^4}  \Big| \\& = \Big| \Big\langle \Big(\beta/2 + \mathscr{C}^{z}_m\Big )^{-1} t_{\Sigma}f, t_{\S}(D_m-z)^{-1} r_{\Omega_-} g \Big\rangle_{\mathit{L}^2(\Sigma)^4, \mathit{L}^2(\Sigma)^4} - \\& \hspace{5cm}\Big\langle \Big(\beta/2 + \mathscr{C}^{z, \varepsilon}_m\Big )^{-1} \mathcal{T}_\varepsilon \mathcal{T}_\varepsilon^{-1}t_{\Sigma^\varepsilon}f, t_{\S^\varepsilon}(D_m-z)^{-1} r_{\Omega^{\varepsilon}_-} g \Big\rangle_{ \mathit{L}^2(\Sigma^\varepsilon)^4, \mathit{L}^2(\Sigma^\varepsilon)^4}  \Big| \\& = \Big| \Big\langle \Big(\beta/2 + \mathscr{C}^{z}_m\Big )^{-1} t_{\Sigma}f, t_{\S}(D_m-z)^{-1} r_{\Omega_-} g \Big\rangle_{\mathit{L}^2(\Sigma)^4,\mathit{L}^2(\Sigma)^4} - \\& \hspace{3cm}\Big\langle \Big(\beta/2 + \mathcal{T}_\varepsilon \mathcal{T}_\varepsilon^{-1} \mathscr{C}^{z, \varepsilon}_m \mathcal{T}_\varepsilon \mathcal{T}_\varepsilon^{-1} \Big )^{-1} \mathcal{T}_\varepsilon \mathcal{T}_\varepsilon^{-1}t_{\Sigma^\varepsilon}f, t_{\S^\varepsilon}(D_m-z)^{-1} r_{\Omega^{\varepsilon}_-} g \Big\rangle_{ \mathit{L}^2(\Sigma^\varepsilon)^4,  \mathit{L}^2(\Sigma^\varepsilon)^4} \Big| 
  \\& = \Big| \Big\langle \Big(\beta/2 + \mathscr{C}^{z}_m\Big )^{-1} t_{\Sigma}f, t_{\S}(D_m-z)^{-1} r_{\Omega_-} g \Big\rangle_{\mathit{L}^2(\Sigma)^4, \mathit{L}^2(\Sigma)^4} - \\& \hspace{4cm}\Big\langle \Big(\beta/2 + \mathcal{T}_\varepsilon^{-1} \mathscr{C}^{z, \varepsilon}_m \mathcal{T}_\varepsilon  \Big )^{-1} \mathcal{T}_\varepsilon^{-1}t_{\Sigma^\varepsilon}f, \mathcal{T}_\varepsilon^{-1} t_{\S^\varepsilon}(D_m-z)^{-1} r_{\Omega^{\varepsilon}_-} g \Big\rangle_{ \mathit{L}^2(\Sigma)^4, \mathit{L}^2(\Sigma)^4}  \Big| .
\end{align*}
By adding and subtracting the term $\Big\langle \Big(\beta/2 + \mathcal{T}_\varepsilon^{-1}\mathscr{C}^{z,\varepsilon}_m \mathcal{T}_\varepsilon\Big )^{-1} \mathcal{T}_\varepsilon^{-1} t_{\Sigma^\varepsilon}f, t_{\S}(D_m-z)^{-1} r_{\Omega_-} g \Big\rangle_{\mathit{L}^2(\Sigma)^4, \mathit{L}^2(\Sigma)^4}$ in the last quantity, we obtain that
\begin{align*}
    \Bigg| \Big\langle& \big[e_{\Omega_-} \Phi^{z}_{m,-} \Big(\beta/2 + \mathscr{C}^{z}_m\Big )^{-1} t_{\Sigma} - e_{\Omega^{\varepsilon}_-} \Phi^{z,\varepsilon}_{m,-} \Big(\beta/2 + \mathscr{C}^{z, \varepsilon}_m\Big )^{-1} t_{\Sigma^\varepsilon} \big]f, g \Big\rangle_{\mathit{L}^2(\rr^3)^4, \mathit{L}^2(\rr^3)^4}  \Bigg| \\& \leq \bigg| \bigg|\Big[\Big(\beta/2 + \mathscr{C}^{z}_m\Big )^{-1} t_{\Sigma} - \Big(\beta/2 + \mathcal{T}_\varepsilon^{-1} \mathscr{C}^{z, \varepsilon}_m \mathcal{T}_\varepsilon  \Big )^{-1} \mathcal{T}_\varepsilon^{-1}t_{\Sigma^\varepsilon}\Big]f
    \bigg| \bigg|_{\mathit{L}^2(\Sigma)^4} \bigg| \bigg| t_{\S}(D_m-z)^{-1} r_{\Omega_-} g \bigg| \bigg|_{\mathit{L}^2(\Sigma)^4} \\& + \bigg| \bigg|\Big(\beta/2 + \mathcal{T}_\varepsilon^{-1} \mathscr{C}^{z, \varepsilon}_m \mathcal{T}_\varepsilon  \Big )^{-1} \mathcal{T}_\varepsilon^{-1}t_{\Sigma^\varepsilon}f\bigg| \bigg|_{\mathit{L}^2(\Sigma)^4}  \bigg| \bigg| \Big[t_{\S}(D_m-z)^{-1} r_{\Omega_-} - \mathcal{T}_\varepsilon^{-1} t_{\S^\varepsilon}(D_m-z)^{-1} r_{\Omega^{\varepsilon}_-} \Big]g \bigg| \bigg|_{\mathit{L}^2(\Sigma)^4} \\& =: r_1 + r_2.
      \end{align*}

Now, let $\mathscr{C}_m^z$ and $\mathcal{T}_\varepsilon^{-1} \mathscr{C}^{z, \varepsilon}_m \mathcal{T}_\varepsilon$ from \eqref{CauchyOpe} and \eqref{TCT} respectively. Then, for a fixed $\rho, \rho^{'}> 0$ such that $\rho^{''} = \mathrm{min}\lbrace \rho ,\rho^{'}\rbrace$, the regularity of $\S$ and $\phi_{m}^{z}$, and  a combination of the mean value theorem give
\begin{align*}
\big| \phi^z_m \big(x_\Sigma - y_\Sigma + \varepsilon (\nu(x_\Sigma) - \nu(y_\Sigma))\big) - \phi_m^z(x_\Sigma - y_\Sigma )\big|\leq \varepsilon\,|\partial \phi_m^z| \leq \varepsilon\, C, \quad\text{with } C \text{ only depending on }z.
\end{align*}
We set $f_{\varepsilon}(y_\Sigma):= \mathrm{det}(1-\varepsilon\nu(x_{\Sigma}))f(y_\Sigma).$ On one hand, using the Cauchy-Schwarz inequality, we obtain that
\begin{align*}
      \Big | \mathscr{C}_m^z f (x_\Sigma)  & - \big(\mathcal{T}_\varepsilon^{-1} \mathscr{C}^{z, \varepsilon}_m \mathcal{T}_\varepsilon f\big)(x_\Sigma)  \Big | \\&\leq \int_{|x_{\Sigma} - y_\Sigma| >\rho^{''}} \big| \phi^z_m \big(x_\Sigma - y_\Sigma + \varepsilon (\nu(x_\Sigma) - \nu(y_\Sigma))\big)f(y_\Sigma) - \phi_m^z(x_\Sigma - y_\Sigma ) f_\varepsilon (y_\Sigma)\big|\mathrm{d}\sigma(y_\Sigma) \\& \leq \int_{\Sigma} \Big| \Big(\phi^z_m \big(x_\Sigma - y_\Sigma + \varepsilon (\nu(x_\Sigma) - \nu(y_\Sigma))\big) - \phi_m^z(x_\Sigma - y_\Sigma) \Big)f(y_\Sigma)\Big| \mathrm{d}\sigma (y_\Sigma) \\& \hspace{5.4cm }+  \int_{\Sigma} \big|  \phi_m^z(x_\Sigma - y_\Sigma) \big( f_\varepsilon (y_\Sigma) - f(y_\Sigma)\big)\big| \mathrm{d}\sigma (y_\Sigma).
\end{align*}
On the other hand, Proposition \ref{PWM} gives us 
$$\mathrm{det}\big(1-\varepsilon W(x_\Sigma)\big) = 1 - \varepsilon \lambda_1(x_\Sigma) - \varepsilon \lambda_2(x_\Sigma) + \varepsilon^2 \lambda_1(x_\Sigma)\lambda_2(x_\Sigma), $$
where $ \lambda_1(x_\Sigma),\,\lambda_2(x_\Sigma)$ are the eigenvalues of the Weingarten map $W(x_\Sigma)$. Then, we get \\
$$ | f_\varepsilon(y_\Sigma) - f(y_\Sigma)| = |\mathrm{det}(1 - \varepsilon W (y_\Sigma)) - 1 ||f(y_\Sigma)|\lesssim \varepsilon ||f||_{L^2(\Sigma)^4}.$$
We conclude that 
\begin{align}\label{diffcauchy}
    \big|\big| \left(\mathscr{C}_{m}^{z}- \mathcal{T}_\varepsilon^{-1} \mathscr{C}^{z, \varepsilon}_m \mathcal{T}_\varepsilon \right) \big|\big|_{L^{2}(\S)^{4}\rightarrow L^{2}(\S)^{4}} = \mathcal{O}(\varepsilon).
\end{align} 
Now, we are going to establish the estimate $r_1$. First, we have that $t_{\S}(D_m-z)^{-1} r_{\Omega_-} $ is bounded from $ L^{2}(\mathbb{R}^3)^4$ into $L^{2}(\Sigma)^4$. On the other hand, using triangular inequality, we get that 
\begin{align*}
    \bigg| \bigg|\Big[\Big(\dfrac{\beta}{2} &+ \mathscr{C}^{z}_m\Big )^{-1} t_{\Sigma} - \Big(\dfrac{\beta}{2} + \mathcal{T}_\varepsilon^{-1} \mathscr{C}^{z, \varepsilon}_m \mathcal{T}_\varepsilon  \Big )^{-1} \mathcal{T}_\varepsilon^{-1}t_{\Sigma^\varepsilon}\Big]f
    \bigg| \bigg|_{\mathit{L}^2(\Sigma)^4} \\&\leq \bigg| \bigg|\Big[\Big(\dfrac{\beta}{2} + \mathscr{C}^{z}_m\Big )^{-1} - \Big(\dfrac{\beta}{2} + \mathcal{T}_\varepsilon^{-1} \mathscr{C}^{z, \varepsilon}_m \mathcal{T}_\varepsilon  \Big )^{-1}\Big]t_{\Sigma}f
    \bigg| \bigg|_{\mathit{L}^2(\Sigma)^4} \\& \hspace{2cm}+ \bigg| \bigg|\Big(\dfrac{\beta}{2} + \mathcal{T}_\varepsilon^{-1} \mathscr{C}^{z, \varepsilon}_m \mathcal{T}_\varepsilon  \Big )^{-1} \Big[\mathcal{T}_\varepsilon^{-1}t_{\Sigma^\varepsilon} - t_\Sigma\Big]f
    \bigg| \bigg|_{\mathit{L}^2(\Sigma)^4} \leq q_1 + q_2.
\end{align*}
To prove the estimate $q_1$, we let $f\in L^{2}(\S)^{4}$ and we set $h=\left(\dfrac{\beta}{2}+\mathscr{C}_{m}^{z}\right)^{-1} t_{\Sigma} f$ bounded from $L^{2}(\S)^{4}$ into itself. Then, the Cauchy-Schwarz inequality and the following statement  
\begin{align}\label{Diff}
    \left(\dfrac{\beta}{2}+ \mathcal{T}_\varepsilon^{-1} \mathscr{C}^{z, \varepsilon}_m \mathcal{T}_\varepsilon\right)^{-1}-\left(\dfrac{\beta}{2}+ \mathscr{C}^{z}_m \right)^{-1}=\left(\dfrac{\beta}{2}+\mathscr{C}_{m}^{z,\varepsilon}\right)^{-1}\left(\mathscr{C}_{m}^{z}-  \mathcal{T}_\varepsilon^{-1} \mathscr{C}^{z, \varepsilon}_m \mathcal{T}_\varepsilon\right)\left(\dfrac{\beta}{2}+\mathscr{C}_{m}^{z}\right)^{-1}
\end{align}
yields that 
\begin{align*}
    q_1 & = \left|\left|\left(\dfrac{\beta}{2}+\mathcal{T}_\varepsilon^{-1} \mathscr{C}^{z, \varepsilon}_m \mathcal{T}_\varepsilon \right)^{-1} \left(\mathscr{C}_{m}^{z}- \mathcal{T}_\varepsilon^{-1} \mathscr{C}^{z, \varepsilon}_m \mathcal{T}_\varepsilon  \right)h\right|\right|_{L^{2}(\S)^{4}}
    \\&\leq \left|\left|\left(\dfrac{\beta}{2}+\mathcal{T}_\varepsilon^{-1} \mathscr{C}^{z, \varepsilon}_m \mathcal{T}_\varepsilon \right)^{-1}\right|\right|_{L^{2}(\S)^{4}\rightarrow L^{2}(\S)^{4}}\big|\big| \left(\mathscr{C}_{m}^{z}- \mathcal{T}_\varepsilon^{-1} \mathscr{C}^{z, \varepsilon}_m \mathcal{T}_\varepsilon \right)h \big|\big|_{L^{2}(\S)^{4}}\\&\leq\big|\big| \left(\mathscr{C}_{m}^{z}- \mathcal{T}_\varepsilon^{-1} \mathscr{C}^{z, \varepsilon}_m \mathcal{T}_\varepsilon \right)h \big|\big|_{L^{2}(\S)^{4}}\\& \lesssim \big|\big| \left(\mathscr{C}_{m}^{z}- \mathcal{T}_\varepsilon^{-1} \mathscr{C}^{z, \varepsilon}_m \mathcal{T}_\varepsilon \right) \big|\big|_{L^{2}(\S)^{4}\rightarrow L^{2}(\S)^{4}} \big|\big| h \big|\big|_{L^{2}(\S)^{4}}\lesssim \big|\big| \left(\mathscr{C}_{m}^{z}- \mathcal{T}_\varepsilon^{-1} \mathscr{C}^{z, \varepsilon}_m \mathcal{T}_\varepsilon \right) \big|\big|_{L^{2}(\S)^{4}\rightarrow L^{2}(\S)^{4}}
\end{align*}
since $\mathscr{C}_m^z$ and $\mathcal{T}_\varepsilon^{-1} \mathscr{C}^{z, \varepsilon}_m \mathcal{T}_\varepsilon$ are bounded from $L^{2}(\S)^{4}$ into itself. Thanks to the estimate \eqref{diffcauchy}, we get that $q_1 = \mathcal{O}(\varepsilon).$\\\\
To prove the estimate $q_2$, we have for $x\in \S^{\varepsilon}$, the following estimate holds  in $L^{2}(\Sigma)^{4}$
\begin{align}\label{TRACECONV}
\bigg| \bigg| t_{\S}(D_m-z)^{-1} r_{\Omega_-} - \mathcal{T}_\varepsilon^{-1} t_{\S^\varepsilon}(D_m-z)^{-1} r_{\Omega^{\varepsilon}_-}  \bigg| \bigg|_{\mathit{L}^2(\mathbb{R}^3)^4 \rightarrow \mathit{L}^2(\Sigma)^4} = \mathcal{O}(\varepsilon).
\end{align}
Next, based on \eqref{TCT}, we immediately get that $\mathcal{T}_\varepsilon^{-1} \mathscr{C}^{z, \varepsilon}_m \mathcal{T}_\varepsilon $ is uniformly bounded from $\mathit{L}^2(\Sigma)^4$ into itself.\\
Thus, together with \eqref{diffcauchy}, \eqref{TRACECONV}, we deduce that $r_2$ has a convergence rate of $\mathcal{O}(\varepsilon).$\\\\

Now, for the same reasons as those used to prove the estimate $q_2$, subsequently, the fact that we have  we immediately deduce that $\Big(\beta/2 + \mathcal{T}_\varepsilon^{-1} \mathscr{C}^{z, \varepsilon}_m \mathcal{T}_\varepsilon  \Big )^{-1}= \Big(\beta/2 + \mathscr{C}^{z}_m \Big )^{-1}+ \mathcal{O}(\varepsilon)$ (see the estimate $q_1$ for more details), we obtain the estimate $r_2$. \\

Thus, we conclude that the statement \eqref{approresol} is valid in $L^2(\mathbb{R}^3)^4$. The proof of Proposition \ref{RMITCONV} is complete.\qed
\begin{lemma}\label{Lemma7}
If the Lorentz scalar is $\mu= 2$ (confinement case). We can identify the domain \eqref{DOMDL} by the following form
\begin{align*}
 \mathrm{dom}(\mathscr{D}_{L}):=\lbrace (\varphi_{+},\varphi_-)\,\in H^{1}(\Omega_+)^{4}\oplus H^{1}(\Omega_-)^{4},\,g\in H^{1/2}(\S)^{4},\, P_+\varphi_- =P_-\varphi_+=0 \text{ on }\S \rbrace,  
\end{align*}
and then, $\mathscr{D}_L = D_{\mathrm{MIT}}^{\Omega_+}\oplus D_{\mathrm{MIT}}^{\Omega_-},$ where $D_{\mathrm{MIT}}^{\Omega_+}$ resp. $D_{\mathrm{MIT}}^{\Omega_-}$ is introduced in Section \ref{DefMIT} resp. Proposition \ref{RMITCONV}.
\end{lemma}
\textbf{Proof.} Using Plemelj-Sokhotski jump
formula from Lemma \ref{prop of C}-(i), and that $ \varphi_{\pm}=t_{\S}u + \textit{C}^{z}_{\pm,m}[g]$, then we get $P_+\varphi_- = -\beta P_-P_+=0$ and $P_-\varphi_+ = -\beta P_+P_-=0$. Moreover, as $P_+\varphi_- + P_-\varphi_+=t_{\Sigma}u+\Lambda^{z}_{+,m}[g]$, we have that  $t_{\S} u=-\Lambda^{z}_{+,m}[g]$.\qed\\
\textbf{Proof of Proposition \ref{coro1.1}.} For $z\in \rho(\mathscr{D}_L)$, we have the following estimate
\begin{equation*}
    \begin{aligned}
   & \left|\left| e_{\Omega_-^{\varepsilon}}R^{\Omega^\varepsilon_-}_{\mathrm{MIT}}(z)r_{\Omega_-^{\varepsilon}} +e_{\Omega_{+}}R^{\Omega_+}_{\mathrm{MIT}}(z)r_{\Omega_{+}} - R_{L}(z)\right|\right|_{\mathit{L}^2(\rr^3)^4} \\& \leq \left|\left| e_{\Omega_-^{\varepsilon}}R^{\Omega^\varepsilon_-}_{\mathrm{MIT}}(z)r_{\Omega_-^{\varepsilon}} +e_{\Omega_{+}}R^{\Omega_+}_{\mathrm{MIT}}(z)r_{\Omega_{+}} +e_{\Omega_{-}}R^{\Omega_-}_{\mathrm{MIT}}(z)r_{\Omega_{-}}-e_{\Omega_{-}}R^{\Omega_-}_{\mathrm{MIT}}(z)r_{\Omega_{-}} - R_{L}(z)\right|\right|_{\mathit{L}^2(\rr^3)^4} \\& \leq  \left|\left| e_{\Omega_-^{\varepsilon}}R^{\Omega^\varepsilon_-}_{\mathrm{MIT}}(z)r_{\Omega_-^{\varepsilon}} -e_{\Omega_{-}}R^{\Omega_-}_{\mathrm{MIT}}(z)r_{\Omega_{-}} \right|\right|_{\mathit{L}^2(\rr^3)^4} +  \Big|\Big|e_{\Omega_{+}}R^{\Omega_+}_{\mathrm{MIT}}(z)r_{\Omega_{+}} +e_{\Omega_{-}}R^{\Omega_-}_{\mathrm{MIT}}(z)r_{\Omega_{-}}- R_{L}(z)\Big|\Big|_{\mathit{L}^2(\rr^3)^4}.
\end{aligned}
\end{equation*}
Then, Proposition \ref{RMITCONV} and Lemma \ref{Lemma7} yield the statement \eqref{stat}.\qed
\begin{remark}
For all $f\in L^{2}(\rr^{3})^{4},\,g\in P_+ L^{2}(\S)^{4}$ the following convergence holds
\begin{align}\label{Liftingconv}
   \big|\big|  e_{\O_-^{\varepsilon}}E_{m}^{{\varepsilon}}(z)[\mathcal{T}_{\varepsilon}] - e_{\O_-}E^-_{m}(z)\big|\big|_{L^{2}(\Sigma)^{4}\rightarrow L^{2}(\mathbb{R}^3)^{4}}=\mathcal{O}(\varepsilon),
\end{align}
where $E_m^-$ is the lifting operator associated with the boundary value problem $(D_m - z)U= 0$ in $\Omega_-$ with $P_+U=0$ on $\Sigma.$
\end{remark}
\textbf{Proof.}
Now, let me show la convergence considered in \eqref{Liftingconv}. To this end, let $\widetilde{g}:=\mathcal{T}_{\varepsilon}g\in P_+ L^{2}(\S^{\varepsilon})^{4}$, then we have
    \begin{align*}
        &\Big|  \langle  e_{\O^{\varepsilon}_-}E_{m}^{\varepsilon}(z)[\mathcal{T}_{\varepsilon}g],f\rangle_{L^{2}(\rr^{3})^{4}} -\langle e_{\O_-}E_{m}^-(z)g,f\rangle_{L^{2}(\rr^{3})^{4}}\Big|\\&=
        \Big| \langle \beta g,\left(\mathcal{T}^{-1}_{\varepsilon}\G_{+}^{\varepsilon}R_{\mathrm{MIT}}^{\Omega^\varepsilon_-}(\bar{z})r_{\O_-^{\varepsilon}}-\G_+ R^-_{\mathrm{MIT}}(\bar{z})r_{\O_-}\right)f\rangle_{L^{2}(\S)^{4}}\Big|\\&
        \leq ||g||_{L^{2}(\S)^{4}} \Big|\Big|\big(\mathcal{T}^{-1}_{\varepsilon}\G_{+}^{\varepsilon}r_{\O_-^{\varepsilon}}e_{\O_-^{\varepsilon}}R_{\mathrm{MIT}}^{\varepsilon}(\bar{z})r_{\O_-^{\varepsilon}}-\G_+ r_{\Omega_-} e_{\Omega_-} R^{\Omega_-}_{\mathrm{MIT}}(\bar{z})r_{\Omega_-}\big)f\Big|\Big|_{L^{2}(\S)^{4}}\\&\lesssim \Big|\Big|\big(\mathcal{T}^{-1}_{\varepsilon}\G_{+}^{\varepsilon}r_{\O_-^{\varepsilon}}e_{\O_-^{\varepsilon}}R_{\mathrm{MIT}}^{\Omega_-^\varepsilon}(\bar{z})r_{\O_-^{\varepsilon}}-\mathcal{T}_{\varepsilon}^{-1}\G_{+}^{\varepsilon}r_{\O_-^{\varepsilon}}e_{\Omega_-}R^{\Omega_-}_{\mathrm{MIT}}(\bar{z})r_{\Omega_-}+\mathcal{T}_{\varepsilon}^{-1}\G_{+}^{\varepsilon}r_{\O_-^{\varepsilon}}e_{\Omega_-}R^{\Omega_-}_{\mathrm{MIT}}(\bar{z})r_{\Omega_-}\\&\hspace{10cm} -\G_+ r_{\Omega_-} e_{\Omega_-} R^{\Omega_-}_{\mathrm{MIT}}(\bar{z})r_{\Omega_-}\big)f\Big|\Big|_{L^{2}(\S)^{4}}\\& \lesssim  \big|\big|\mathcal{T}^{-1}_{\varepsilon}\G_{+}^{\varepsilon}r_{\O_-^{\varepsilon}}\big|\big|_{L^{2}(\S)^{4}} \Big|\Big| e_{\O_-^{\varepsilon}}R_{\mathrm{MIT}}^{\Omega_-^\varepsilon}(\bar{z})r_{\O_-^{\varepsilon}}f - e_{\Omega_-}R^{\Omega_-}_{\mathrm{MIT}}(\bar{z})r_{\Omega_-}f\Big|\Big|_{L^{2}(\rr^3)^{4}}    \\& \hspace{6cm}+\Big|\Big|\mathcal{T}^{-1}_{\varepsilon}\G_{+}^{\varepsilon}r_{\O_-^{\varepsilon}} - \G_+ r_{\O_-}\Big|\Big|_{L^{2}(\S)^{4}} \big|\big|e_{\O_-}R^{\Omega_-}_{\mathrm{MIT}}r_{\Omega_-} (\bar{z})f\big|\big|_{L^{2}(\rr^3)^{4}}. 
    \end{align*}
    Since $\G^{\varepsilon}_{+}$ is bounded form $L^{2}(\O_-^{\varepsilon})^{4}$ to $L^{2}(\S^{\varepsilon})^{4}$ for $\varepsilon$ small enough, then $\mathcal{T}^{-1}_{\varepsilon}\G_{+}^{\varepsilon}r_{\O_-^{\varepsilon}}$  is bounded in ${L^{2}(\S)^{4}}$. Thus, together with the boundedness of $e_{\Omega_-}R^{\Omega_-}_{\mathrm{MIT}}$ in ${L^{2}(\rr^3)^{4}}$ and the convergence established in Proposition \ref{coro1.1}, we get 
    \begin{align*}
        \left|  \langle  e_{\O_-^{\varepsilon}}E_{m}^{{\varepsilon}}(z)[\mathcal{T}_{\varepsilon}g],f\rangle_{L^{2}(\rr^{3})^{4}} -\langle e_{\O_-}E^-_{m}(z)g,f\rangle_{L^{2}(\rr^{3})^{4}}\right|\lesssim \varepsilon, \quad \text{for all } \,f\in L^{2}(\mathbb{R}^3)^{4}.
    \end{align*}
Since this is true for all $g\in L^{2}(\Sigma)^4$, by duality arguments it follows that
\begin{align*}
    \big|\big|  e_{\O_-^{\varepsilon}}E_{m}^{{\varepsilon}}(z)[\mathcal{T}_{\varepsilon}] - e_{\O_-}E^-_{m}(z)\big|\big|_{L^{2}(\Sigma)^{4}\rightarrow L^{2}(\mathbb{R}^3)^{4}}=\mathcal{O}(\varepsilon).
\end{align*}

\begin{lemma}\label{Estimations des op} Let $K\subset\mathbb{C}$ be a compact set. Then, there exists  $M_0>0$ such that for all  $M>M_0$, for $\varepsilon\in(0,\varepsilon_0)$, $K\subset\rho(D_{\mathrm{MIT}}^{\mathcal{U}^\varepsilon}(m+M))$, and for $z\in K$ the following estimates hold:
\begin{align*}
 \left|\left| e_{\mathcal{U}^{\varepsilon}}R^{\mathcal{U}^\varepsilon}_{\mathrm{MIT}}(z)r_{\mathcal{U}^{\varepsilon}}f\right|\right|_{\mathit{L}^2(\rr^{3})^4}&\lesssim \frac{1}{M} \left|\left| f \right|\right|_{\mathit{L}^2(\rr^{3})^4},\quad &\forall& \,\,f\in \mathit{L}^2(\rr^{3})^4,\\
 \left|\left| \G_{-+}^{\varepsilon}R^{\mathcal{U}^\varepsilon}_{\mathrm{MIT}}(z) r_{\mathcal{U}^{\varepsilon}}f\right|\right|_{\mathit{L}^2(\partial\mathcal{U}^{\varepsilon})^4}&\lesssim \frac{1}{\sqrt{M}} \left|\left| f \right|\right|_{\mathit{L}^2(\rr^{3})^4}, \quad &\forall& \,\,f\in \mathit{L}^2(\rr^{3})^4,\\
 \left|\left| \G_{-+}^{\varepsilon}R^{\mathcal{U}^\varepsilon}_{\mathrm{MIT}}(z)r_{\mathcal{U}^{\varepsilon}}f\right|\right|_{H^{-1/2}(\partial\mathcal{U}^{\varepsilon})^4} & \lesssim \frac{1}{M}\left|\left| f\right|\right|_{\mathit{L}^2(\rr^{3})^4},\quad &\forall& \,\,f\in \mathit{L}^2(\rr^{3})^4,\\
 \left|\left|e_{\mathcal{U}^\varepsilon}\mathcal{E}^{\varepsilon}_{m+M}(z)(\psi,\mathcal{T}_{\varepsilon}\varphi) \right|\right|_{\mathit{L}^2(\mathbb{R}^3)^4}&\lesssim \frac{1}{\sqrt{M}}\left|\left| \psi \right|\right|_{\mathit{L}^{2}(\S)^4}\left|\left| \varphi \right|\right|_{ \mathit{L}^{2}(\S)^4},\quad &\forall&\,\, (\psi,\mathcal{T}_{\varepsilon}\varphi)\in P_+L^2(\Sigma)^4 \oplus P_- \mathit{L}^{2}(\S^{\varepsilon})^4,
 \end{align*}
 \begin{align*}
 \left|\left|e_{\mathcal{U}^\varepsilon}\mathcal{E}^{\varepsilon}_{m+M}(z)(\psi,\mathcal{T}_{\varepsilon}\varphi) \right|\right|_{\mathit{L}^2(\mathbb{R}^3)^4}&\lesssim \frac{1}{M}\left|\left| \psi \right|\right|_{H^{1/2}(\S)^4}\left|\left| \varphi \right|\right|_{ H^{1/2}(\S)^4},\\&\hspace{4cm}\forall\,\, (\psi,\mathcal{T}_{\varepsilon}\varphi)\in P_+H^{1/2}(\Sigma)^4 \oplus P_- H^{1/2}(\S^{\varepsilon})^4.
 \end{align*}
 \end{lemma}	
\textbf{Proof.}
Using the same arguments as in the proof of \cite[Lemma 6.1]{BBZ}, we can show the above estimates with respect to $M$. First, I want to show the claimed estimates for $e_{\mathcal{U}^{\varepsilon}}R^{\mathcal{U}^\varepsilon}_{\mathrm{MIT}}(z)r_{\mathcal{U}^{\varepsilon}}$ and $\G^{\varepsilon}_-R^{\mathcal{U}^\varepsilon}_{\mathrm{MIT}}(z)r_{\mathcal{U}^{\varepsilon}}$. For this,
fix a compact set $K\subset\cc$,  and note that for $z\in K$ and $M_1> \sup_{z\in K}\{  |\mathrm{Re}(z)| -m\} $ it holds that $K\subset\rho(D_{m+M_1})$, and hence $K\subset\rho(D^{\mathcal{U}^\varepsilon}_{\mathrm{MIT}})$  for all $M>M_1$. Let $f\in L^{2}(\rr^{3})^{4}.$ We have that 
\begin{align*}
    ||e_{\mathcal{U}^{\varepsilon}} R^{\mathcal{U}^\varepsilon}_{\mathrm{MIT}}(z)r_{\mathcal{U}^{\varepsilon}}f||_{L^{2}(\rr^{3})^{4}}=||R^{\mathcal{U}^\varepsilon}_{\mathrm{MIT}}(z)r_{\mathcal{U}^{\varepsilon}}f||_{L^{2}(\mathcal{U}^{\varepsilon})^{4}}.
\end{align*}
Now, for $r_{\mathcal{U}^{\varepsilon}}f\in L^{2}(\mathcal{U}^{\varepsilon})^{4}$ and $\varphi\in \mathrm{dom}(D_{\mathrm{MIT}}^{\mathcal{U}^\varepsilon}),$ then a straightforward application of the Green's formula yields that
\begin{align*}
\begin{split}
 \| D^{\mathcal{U}^\varepsilon}_{\mathrm{MIT}}\var \|^2_{\mathit{L}^2(\mathcal{U}^{\varepsilon})^4}=&  \| (\alpha\cdot\nabla)\var \|^2_{\mathit{L}^2(\mathcal{U}^{\varepsilon})^4}+(m+M)^2  \left|\left| \var \right|\right|^2_{\mathit{L}^2(\mathcal{U}^{\varepsilon})^4}+(m+M)  \left|\left| P^{\mathcal{U}^{\varepsilon}}_- t_{\partial\mathcal{U}^{\varepsilon}}\var \right|\right|^2_{\mathit{L}^2(\partial\mathcal{U}^{\varepsilon})^4},
\end{split}
\end{align*}
with $P^{\mathcal{U}^{\varepsilon}}_- t_{\partial\mathcal{U}^{\varepsilon}} = P_- t_{\Sigma}+ P_+ t_{\Sigma^{\varepsilon}}$.
Using this and the Cauchy-Schwarz inequality we obtain that 
\begin{align*}
 \|  (D^{\mathcal{U}^\varepsilon}_{\mathrm{MIT}}-z)\var \|^2_{\mathit{L}^2(\mathcal{U}^{\varepsilon})^4}= & \| D^{\mathcal{U}^\varepsilon}_{\mathrm{MIT}}\var \|^2_{\mathit{L}^2(\mathcal{U}^{\varepsilon})^4} +|z|^2 \| \var \|^2_{\mathit{L}^2(\mathcal{U}^{\varepsilon})^4}  -2\mathrm{Re}(z) \langle  D^{\mathcal{U}^\varepsilon}_{\mathrm{MIT}}\var, \var \rangle_{\mathit{L}^2(\mathcal{U}^{\varepsilon})^4}\\
 \geqslant&  \|  D_{\mathrm{MIT}}^{\mathcal{U}^\varepsilon}\var \|^2_{\mathit{L}^2(\mathcal{U}^{\varepsilon})^4} +|z|^2 \| \var \|^2_{\mathit{L}^2(\mathcal{U}^{\varepsilon})^4} -\frac{1}{2} \|  D^{\mathcal{U}^\varepsilon}_{\mathrm{MIT}}\var \|^2_{\mathit{L}^2(\mathcal{U}^{\varepsilon})^4} -2|\mathrm{Re}(z)|^2 \| \var \|^2_{\mathit{L}^2(\mathcal{U}^{\varepsilon})^4}\\
 \geqslant& \left(\frac{(m+M)^2}{2}+|\mathrm{Im}(z)|^2 -|\mathrm{Re}(z)|^2\right)   \left|\left| \var \right|\right|^2_{\mathit{L}^2(\mathcal{U}^{\varepsilon})^4}  +\frac{M}{2}  \left|\left| P^{\mathcal{U}^{\varepsilon}}_- t_{\partial\mathcal{U}^{\varepsilon}}\var \right|\right|^2_{\mathit{L}^2(\partial\mathcal{U}^{\varepsilon})^4}.
\end{align*}
Therefore,   taking $R^{\mathcal{U}^\varepsilon}_{\mathrm{MIT}}(z)r_{\mathcal{U}^{\varepsilon}}f=\var$ and  $M\geqslant M_2\geqslant \sup_{z\in K}\{  \sqrt{|\mathrm{Re}(z)|^2-|\mathrm{Im}(z)|^2} -m\}$  we obtain the inequality 
\begin{align*}
 \left|\left| R^{\mathcal{U}^\varepsilon}_{\mathrm{MIT}}(z)r_{\mathcal{U}^{\varepsilon}}f\right|\right|_{\mathit{L}^2(\mathcal{U}^{\varepsilon})^4} +\quad\frac{1}{\sqrt{M}}\left|\left|\G_{-+}^\varepsilon R^{\mathcal{U}^\varepsilon}_{\mathrm{MIT}}(z)r_{\mathcal{U}^{\varepsilon}}f\right|\right|_{\mathit{L}^2(\partial\mathcal{U}^{\varepsilon})^4} \lesssim \frac{1}{M}\left|\left| f\right|\right|_{\mathit{L}^2(\rr^{3})^4},\,\text{with}\quad\partial\mathcal{U}^{\varepsilon}=\S\cup\S^{\varepsilon}.
  \end{align*} 
Thus
\begin{align*}
 \left|\left| e_{\mathcal{U}^{\varepsilon}}R^{\mathcal{U}^\varepsilon}_{\mathrm{MIT}}(z) r_{\mathcal{U}^{\varepsilon}}f\right|\right|_{\mathit{L}^2(\rr^{3})^4} \lesssim \frac{1}{M} \left|\left| f \right|\right|_{\mathit{L}^2(\rr^{3})^4}, \quad\text{ and }\quad
      \left|\left| \G_{-+}^{\varepsilon}R^{\mathcal{U}^\varepsilon}_{\mathrm{MIT}}(z)r_{\mathcal{U}^{\varepsilon}}f\right|\right|_{\mathit{L}^2(\partial\mathcal{U}^{\varepsilon})^4} \lesssim \frac{1}{\sqrt{M}}\left|\left| f\right|\right|_{\mathit{L}^2(\rr^{3})^4}.
  \end{align*}
 Since $\G_{-+}^{\varepsilon}:=(\G_-, \G^{\varepsilon}_+)$ is bounded from $\mathit{L}^2(\mathcal{U}^\varepsilon)^4$ into $\mathit{H}^{-1/2}(\partial\mathcal{U}^{\varepsilon})^4$ for $\varepsilon\in (0,\varepsilon_0)$ with $\varepsilon_0$ sufficiently small, it follows from the above inequality that 
\begin{align*}
 \left|\left| \G_{-+}^{\varepsilon}R^{\mathcal{U}^\varepsilon}_{\mathrm{MIT}}(z)r_{\mathcal{U}^{\varepsilon}}f\right|\right|_{\mathit{H}^{-1/2}(\partial\mathcal{U}^{\varepsilon})^4}&\lesssim  \left|\left| \G_{-+}^{\varepsilon}\right|\right|_{\mathit{L}^2(\mathcal{U}^{\varepsilon})^4 \rightarrow\mathit{H}^{-1/2}(\partial\mathcal{U}^{\varepsilon})^4} \left|\left|R^{\mathcal{U}^\varepsilon}_{\mathrm{MIT}}(z) r_{\mathcal{U}^{\varepsilon} }f\right|\right|_{\mathit{L}^{2}(\mathcal{U}^{\varepsilon})^4}\\&\lesssim   \frac{1}{M} \left|\left| f \right|\right|_{\mathit{L}^2(\rr^3)^4},
 \end{align*}
 for any $ f\in \mathit{L}^2(\mathbb{R}^{3})^4$, which gives the last inequality. \\\\
 Let us now turn to the proof of the claimed estimates for $e_{\mathcal{U}^{\varepsilon}}\mathcal{E}^{\varepsilon}_{m+M}(z).$
Let $f,\,\psi$ belong to $L^{2}(\rr^{3})^{4}$ and $L^{2}(\S)^{4}$, respectively, and consider the transformation operator $\mathcal{T}_{\varepsilon}$ defined in \eqref{Transformationoperator}. For $\varphi\in L^{2}(\Sigma)^4$, we set $\varphi_{\varepsilon}= \mathcal{T}_{\varepsilon}\varphi \in L^2(\Sigma^\varepsilon)$.
We mention that $\beta(\G_-,\G^{\varepsilon}_+)R^{\mathcal{U}^\varepsilon}_{\mathrm{MIT}}(\overline{z})$ is the adjoint of the operator $\mathcal{E}^{\varepsilon}_{m+M}(z): P_+ L^2 (\Sigma)^4\oplus P_- L^{2}(\S^{\varepsilon})^4\longrightarrow \mathit{L}^2(\mathcal{U}^{\varepsilon})^4$. 
Using this and the estimate fulfilled by $(\G_-,\G^{\varepsilon}_+)R^{\mathcal{U}^\varepsilon}_{\mathrm{MIT}}(\overline{z})r_{\mathcal{U}^{\varepsilon}}$  we obtain that 
 \begin{align*}
\left| \langle f,e_{\mathcal{U}^{\varepsilon}}\mathcal{E}^{\varepsilon}_{m+M}(z)(\psi,\varphi_\varepsilon) \rangle_{\mathit{L}^{2}(\rr^{3})^4}\right| &=\left| \langle (\G_-, \mathcal{T}^{-1}_\varepsilon\G^{\varepsilon}_+)R^{\mathcal{U}^\varepsilon}_{\mathrm{MIT}}(\overline{z})r_{\mathcal{U}_{\varepsilon}}f,\beta(\psi,\varphi) \rangle_{\mathit{L}^{2}(\S)^4}\right|\\
&\leqslant  \left|\left| \big(\G_-, \mathcal{T}_{\varepsilon}^{-1}\G^{\varepsilon}_+\big)R^{\mathcal{U}^\varepsilon}_{\mathrm{MIT}}(z) r_{\mathcal{U}^{\varepsilon}}f\right|\right|_{\mathit{L}^{2}(\S)^4}\left|\left| \psi \right|\right|_{\mathit{L}^{2}(\S)^4} \left|\left| \varphi\right|\right|_{\mathit{L}^{2}(\S)^4}\\
&\leqslant  \left|\left| \psi \right|\right|_{\mathit{L}^{2}(\S)^4}\left|\left| \varphi \right|\right|_{\mathit{L}^{2}(\S)^4} \left|\left| \mathcal{T}_{\varepsilon}^{-1}\right|\right|_{L^{2}(\partial\mathcal{U}^{\varepsilon})^{4}\rightarrow L^{2}(\S)^{4}}
\left|\left|\G_{-+}^{\varepsilon}R^{\mathcal{U}^\varepsilon}_{\mathrm{MIT}}(z)r_{\mathcal{U}^{\varepsilon}}f\right|\right|_{\mathit{L}^{2}(\partial\mathcal{U}^{\varepsilon})^4}\\
&\lesssim \frac{1}{\sqrt{M}} \left|\left|  f\right|\right|_{L^{2}(\rr^{3})^4}\left|\left| \psi \right|\right|_{\mathit{L}^{2}(\S)^4} \left|\left| \varphi \right|\right|_{\mathit{L}^{2}(\S)^4}.
\end{align*}
So, we get 
\begin{align*}
    \left|\left|e_{\mathcal{U}^{\varepsilon}}\mathcal{E}^{\varepsilon}_{m+M}(z)(\psi,\mathcal{T}_{\varepsilon}\varphi) \right|\right|_{\mathit{L}^2(\rr^{3})^4}&\lesssim \frac{1}{\sqrt{M}}\left|\left| \psi \right|\right|_{\mathit{L}^{2}(\S)^4} \left|\left| \varphi \right|\right|_{\mathit{L}^{2}(\S)^4}.
\end{align*}
Similarly, we established the last inequality of the lemma andthis finishes the proof of the lemma.\qed\\

The last ingredient to prove Theorem \ref{maintheo} is to show that the second term in the ride hand side of the resolvent formula \eqref{kfullreso} converges to zero when $M$ converges to $\infty$, (\emph{i.e.}, $h= \varepsilon = M^{-1}\rightarrow 0).$ 
\\

\textbf{Proof of Proposition \ref{Prop1.1}.}
Recall the following notations: $D^{\Omega^{\varepsilon}_{+-}}_{\mathrm{MIT}} = D^{\Omega_{+}}_{\mathrm{MIT}} \oplus D^{\Omega^{\varepsilon}_{-}}_{\mathrm{MIT}}$ and $R^{\Omega^{\varepsilon}_{+-}}_{\mathrm{MIT}} = R^{\Omega_{+}}_{\mathrm{MIT}}\oplus R^{\Omega^{\varepsilon}_{-}}_{\mathrm{MIT}}$, with $\Omega^{\varepsilon}_{+-}= \Omega_+ \cup \Omega^\varepsilon_-$. Let $z\in \rho(\mathfrak{D}^\varepsilon_M)\cap\rho(D^{\Omega^{\varepsilon}_{+-}}_{\mathrm{MIT}})$ and $f\in L^2(\mathbb{R}^3)^4$. From the resolvent formula \eqref{kfullreso} and Remark \ref{remark6.1}, together give us the following
\begin{equation*}
    \begin{aligned}
         \big|\big| R_M^\varepsilon(z) - e_{\Omega^\varepsilon_{+-}}R^{\Omega^{\varepsilon}_{+-}}_{\mathrm{MIT}}(z)r_{\Omega^\varepsilon_{+-}}\big|\big|_{L^2(\mathbb{R}^3)^4 \rightarrow L^2(\mathbb{R}^3)^4} & \leq \big|\big| e_{\mathcal{U}^\varepsilon}R^{\mathcal{U}^\varepsilon}_{\mathrm{MIT}}(z)r_{\mathcal{U}^\varepsilon}f\big|\big|_{L^2(\mathbb{R}^3)^4} \\& + \big|\big| E^{\Omega^{\varepsilon}_{+-}}_m(z)\Xi_M^{\varepsilon,-+}(z)\mathscr{A}^\varepsilon_{m+M}\G_{+-}^\varepsilon R^{\Omega^{\varepsilon}_{+-}}_{\mathrm{MIT}}(z)r_{\Omega^{\varepsilon}_{+-}}f \big|\big|_{L^2(\Omega^{\varepsilon}_{+-})^4} \\& +  \big|\big| E^{\Omega^{\varepsilon}_{+-}}_m(z)\Xi_M^{\varepsilon,-+}(z)\G_{-+}^\varepsilon  R^{\mathcal{U}^\varepsilon}_{\mathrm{MIT}}(z)r_{\mathcal{U}^\varepsilon}f \big|\big|_{L^2(\Omega^{\varepsilon}_{+-})^4} \\& + \big|\big| \mathcal{E}^{\varepsilon}_{m+M}(z)\Xi_M^{\varepsilon,+-}(z)\G_{+-}^\varepsilon R^{\Omega^{\varepsilon}_{+-}}_{\mathrm{MIT}}(z)r_{\Omega^{\varepsilon}_{+-}}f \big|\big|_{L^2(\mathcal{U}^\varepsilon)^4} \\&+ \big|\big| \mathcal{E}^{\varepsilon}_{m+M}(z)\Xi_M^{\varepsilon,+-}(z)\mathscr{A}_m^{\Omega^{\varepsilon}_{+-}}\G_{-+}^\varepsilon R^{\mathcal{U}^\varepsilon}_{\mathrm{MIT}}(z)r_{\mathcal{U}^\varepsilon}f \big|\big|_{L^2(\mathcal{U}^\varepsilon)^4} \\& =: J_1 +J_2 +J_3 + J_4 + J_5.
    \end{aligned}
\end{equation*}
We start with $J_1$. From the second item of Lemma \ref{Estimations des op}, we get that $J_1 \lesssim M^{-1}.$ Now, thanks to the uniform bound (with respect to $M$) of $\Xi_M^{\varepsilon,\pm\mp}$, see Corollary \ref{une autre est}, $J_2$, $J_3$, $J_4$, $J_5$ become as follows
\begin{align*}
    J_2 &\lesssim \big|\big| E^{\Omega^{\varepsilon}_{+-}}_m(z)\big|\big|_{L^{2}(\Omega^{\varepsilon}_{+-})^4} \big|\big|\mathscr{A}^\varepsilon_{m+M}\big|\big|_{H^{-1/2}(\Sigma)^4\oplus H^{-1/2}(\Sigma^\varepsilon)^4}\big|\big|\G_{+-}^\varepsilon R^{\Omega^{\varepsilon}_{+-}}_{\mathrm{MIT}}(z)r_{\Omega}f \big|\big|_{H^{1/2}(\Sigma)^4\oplus H^{1/2}(\Sigma^\varepsilon)^4},\\
     J_3 &\lesssim \big|\big| E^{\Omega^{\varepsilon}_{+-}}_m(z)\big|\big|_{H^{-1/2}(\Sigma)^4\oplus H^{-1/2}(\Sigma^\varepsilon)^4\rightarrow L^{2}(\Omega^{\varepsilon}_{+-})^4} \big|\big|\G_{-+}^\varepsilon R^{\mathcal{U}^\varepsilon}_{\mathrm{MIT}}(z) r_{\mathcal{U}^\varepsilon}f \big|\big|_{H^{-1/2}(\Sigma)^4\oplus H^{-1/2}(\Sigma^\varepsilon)^4},\\
      J_4 &\lesssim \big|\big| \mathcal{E}^{\varepsilon}_{m+M}(z)\big|\big|_{H^{1/2}(\Sigma)^4\oplus H^{1/2}(\Sigma^\varepsilon)^4\rightarrow L^{2}(\mathcal{U}^\varepsilon)^4} \big|\big|\G_{+-}^\varepsilon R^{\Omega^{\varepsilon}_{+-}}_{\mathrm{MIT}}(z)r_{\Omega}f \big|\big|_{H^{1/2}(\Sigma)^4\oplus H^{1/2}(\Sigma^\varepsilon)^4},\\ 
      J_5 &\lesssim \big|\big| \mathcal{E}^{\varepsilon}_{m+M}(z)\big|\big|_{L^{2}(\mathcal{U}^\varepsilon)^4}\big|\big|\mathscr{A}_m^{\Omega^{\varepsilon}_{+-}}\big|\big|_{L^{2}(\Sigma)^4\oplus L^{2}(\Sigma^\varepsilon)^4}\big|\big|\G_{-+}^\varepsilon R^{\mathcal{U}^\varepsilon}_{\mathrm{MIT}}(z)r_{\mathcal{U}^\varepsilon}f \big|\big|_{L^{2}(\Sigma)^4\oplus L^{2}(\Sigma^\varepsilon)^4}.
\end{align*}
Notice that the terms $E^{\Omega^{\varepsilon}_{+-}}_m, \, \mathscr{A}_m^{\Omega^{\varepsilon}_{+-}},$ and $\G_{+-}^\varepsilon  R^{\Omega^{\varepsilon}_{+-}}_{\mathrm{MIT}}(z)$ are bounded operators for all $\varepsilon\in(0,\varepsilon_0)$, everywhere defined and do not depend on $M.$ Now, thanks to Lemma \ref{Estimations des op},  $\G_{-+}^\varepsilon R^{\mathcal{U}^\varepsilon}_{\mathrm{MIT}}(z)r_{\mathcal{U}^\varepsilon}$ and $e_{\mathcal{U}^\varepsilon}\mathcal{E}^{\varepsilon}_{m+M}(z)$ hold the following estimate
\begin{align*}
     \left|\left| \G_{-+}^{\varepsilon}R^{\mathcal{U}^\varepsilon}_{\mathrm{MIT}}(z) r_{\mathcal{U}^{\varepsilon}}f\right|\right|_{\mathit{L}^2(\partial\mathcal{U}^{\varepsilon})^4}\lesssim \frac{1}{\sqrt{M}} \left|\left| f \right|\right|_{\mathit{L}^2(\rr^{3})^4}&\text{ and }
 \left|\left| \G_{-+}^{\varepsilon}R^{\mathcal{U}^\varepsilon}_{\mathrm{MIT}}(z)r_{\mathcal{U}^{\varepsilon}}f\right|\right|_{H^{-1/2}(\partial\mathcal{U}^{\varepsilon})^4}  \lesssim \frac{1}{M}\left|\left| f\right|\right|_{\mathit{L}^2(\rr^{3})^4},\\
 \left|\left|e_{\mathcal{U}^\varepsilon}\mathcal{E}^{\varepsilon}_{m+M}(z)(\psi,\mathcal{T}_{\varepsilon}\varphi) \right|\right|_{\mathit{L}^2(\mathbb{R}^3)^4}&\lesssim \frac{1}{\sqrt{M}}\left|\left| \psi \right|\right|_{\mathit{L}^{2}(\S)^4}\left|\left| \varphi \right|\right|_{ \mathit{L}^{2}(\S)^4},\\
 \left|\left|e_{\mathcal{U}^\varepsilon}\mathcal{E}^{\varepsilon}_{m+M}(z)(\psi,\mathcal{T}_{\varepsilon}\varphi) \right|\right|_{\mathit{L}^2(\mathbb{R}^3)^4}&\lesssim \frac{1}{M}\left|\left| \psi \right|\right|_{H^{1/2}(\S)^4}\left|\left| \varphi \right|\right|_{ H^{1/2}(\S)^4}.
\end{align*}
Thus, from the above estimates, we deduce that 
\begin{align*}
    J_k \lesssim M^{-1} ||f||_{L^2(\mathbb{R}^3)^4}, \quad \forall\,k\in\lbrace 3,4,5\rbrace.
\end{align*}
Moreover, the following lower bound  of $\mathcal{A}_{m+M}^\varepsilon$, see Corollary \eqref{ESTAME},  
$$||\mathcal{A}_{m+M}^\varepsilon||_{H^{1/2}(\Sigma)^4\oplus H^{1/2}(\Sigma^\varepsilon)^4\rightarrow H^{-1/2}(\Sigma)^4\oplus H^{-1/2}(\Sigma^\varepsilon)^4  }\lesssim M^{-1},$$
yields that $J_2\lesssim M^{-1} ||f||_{L^2(\mathbb{R}^3)^4}.$
Thus, we obtain the estimate 
$$ \big|\big| R_M^\varepsilon(z) - e_{\Omega^{\varepsilon}_{+-}} R_{\mathrm{MIT}}^{\Omega^{\varepsilon}_{+-}}(z)r_{\Omega^{\varepsilon}_{+-}}\big|\big|_{L^2(\mathbb{R}^3)^4 \rightarrow L^2(\mathbb{R}^3)^4}\lesssim M^{-1} ||f||_{L^2(\mathbb{R}^3)^4}.$$
And this achieves the proof of the proposition.\qed\\

Thus, Theorem \ref{maintheo} is then obtained by a simple combination of Propositions \ref{coro1.1}, \ref{Prop1.1}.
\section*{Acknowledgement}
This work was supported by LTC Transmath, BERC.2022-2025 program and BCAM Severo Ochoa research project. I wish to express my gratitude to my thesis advisor Vincent Bruneau for suggesting the problem and for many stimulating conversations, for his patient advice, and enthusiastic encouragement. 
\appendix
\section{}\label{Appendix}
For a better understanding of the construction of the approximation of the solutions $A_j(y,\xi,\tau)$ and the order of the coefficients $B_{j,k}(y,\xi)$ as well as the proof of Proposition \ref{SolSystSj}, an explicit calculation is presented in this appendix, which aims to obtain an exact form of the solutions $A_j(y,\xi,\tau)$ for $j=1,2$.\\

For $j= 1$, we define $A_{1}(y, \xi, \tau)$ inductively by 
\begin{equation}\label{T1approx}	
	\left\{
\begin{aligned}
&h\partial_{\tau} A_{1}(y, \xi, \tau)	= L_{0}(y, \xi)A_{1}(y, \xi, \tau)+\Big(L_1(y,\xi)+(\alpha\cdot\tilde{\nu}^{\varphi}c_{3})L_{0}(y,\xi) - i \partial_{\xi}L_{0}(y,\xi)\cdot\partial_{y}\Big)A_{0}(y, \xi, \tau),\\
&P_+ A_{1}(y, \xi, \varepsilon)  =  0,  
	\end{aligned}
	\right.
	\end{equation}
we have $\partial_{\xi}L_{0}(y,\xi)\cdot\partial_{y}=i\alpha\cdot\tilde{\nu}^{\varphi}(\alpha\cdot\partial_{y}):=a_{0}(y)(\alpha\cdot\partial_{y})$, with $a_{0}(y)=i\alpha\cdot\tilde{\nu}^{\varphi}.$
	 The solution of the differential system \eqref{T1approx} is 
\begin{align*}
A_{1}(y, \xi, \tau)&=
e^{h^{-1}L_{0}(\tau-\varepsilon)}A_{1}(y, \xi, \varepsilon)\\&+e^{h^{-1}L_{0}\tau}\int_{\varepsilon}^{\tau}e^{-h^{-1}L_{0}(y,\xi)s}  \Big(L_1+(\alpha\cdot\tilde{\nu}^{\varphi}c_{3})L_{0} -i \partial_{\xi}L_{0}(y,\xi)\cdot\partial_{y}\Big)  A_{0}(y, \xi, \tau)\mathrm{d}s\\&=e^{h^{-1}L_{0}(\tau-\varepsilon)}A_{1}(y, \xi, \varepsilon)\\&+e^{h^{-1}L_{0}\tau}\int_{\varepsilon}^{\tau}e^{-h^{-1}L_{0}s}a_{0}(y)\Big(-z+c\cdot\xi -ic_{3}L_{0}-i\alpha\cdot\partial_{y} \Big)A_{0}(y, \xi, \tau)\mathrm{d}s\\&:=I_{1}+I_{2},
\end{align*}
where $I_1$ and $I_2$ have the following quantity:
    \begin{align*}
        I_{1}&= \left(e^{(\tau-\varepsilon)\varrho_{-}(y, \xi)}\Pi_{-}+e^{(\tau-\varepsilon)\varrho_{-}(y, \xi)}\Pi_{+}\right)A_{1}(y, \xi, \varepsilon),
         \\
 I_{2}&=e^{h^{-1}L_{0}(y,\xi)\tau}\int_{\varepsilon}^{\tau}e^{-h^{-1}L_{0}(y,\xi)s}a_{0}(y)\Big(-z+c\cdot\xi -ic_{3}L_{0}-i\alpha\cdot\partial_{y} \Big)A_{0}(y, \xi, s)\mathrm{d}s .
\end{align*}
Now, to obtain an explicit form of $I_2$, let's decompose the quantity $ e^{-h^{-1}L_{0}(y,\xi)s}.$ To do this, we have 
\begin{equation}\label{Expenontielle}
\begin{aligned}
&\int_{\varepsilon}^{\tau}e^{-h^{-1}L_{0}(y,\xi)s}a_{0}(y)\Big(-z+c\cdot\xi -ic_{3}L_{0}-i\alpha\cdot\partial_{y} \Big)A_{0}(y, \xi, s)\mathrm{d}s\\&=
\int_{\varepsilon}^{\tau}\left(e^{-h^{-1}s\varrho_{-}(y, \xi)}\Pi_{-}+e^{-h^{-1}s\varrho_{+}(y, \xi)}\Pi_{+}\right)a_{0}(y)\Big(-z+c\cdot\xi -ic_{3}L_{0}-i\alpha\cdot\partial_{y} \Big)A_{0}(y, \xi, s)\mathrm{d}s
\\&=\int_{\varepsilon}^{\tau}\left(e^{-h^{-1}s\varrho_{-}}\Pi_{-}+e^{-h^{-1}s\varrho_{+}}\Pi_{+}\right)a_{0}(y)\Big(-z+c\cdot\xi -ic_{3}L_{0}-i\alpha\cdot\partial_{y} \Big) \left(e^{h^{-1}(s-\varepsilon)\varrho_{-}}\dfrac{\Pi_{-}P_+ }{k_{-}^{\varphi}}\right)\mathrm{d}s\\&=\underbrace{\int_{\varepsilon}^{\tau}e^{-h^{-1}s\varrho_{-}}\Pi_{-}a_{0}(y)\Big(-z + c\cdot\xi -ic_{3}L_{0}-i\alpha\cdot\partial_{y} \Big) \left(e^{h^{-1}(s-\varepsilon)\varrho_{-}}\dfrac{\Pi_{-}P_+}{k_{-}^{\varphi}}\right)\mathrm{d}s}_{(1)}\\&+\underbrace{\int_{\varepsilon}^{\tau}e^{-h^{-1}s\varrho_{+}}\Pi_{+}a_{0}(y)\Big(-z+c\cdot\xi -ic_{3}L_{0}-i\alpha\cdot\partial_{y} \Big) \left(e^{h^{-1}(s-\varepsilon)\varrho_{-}}\dfrac{\Pi_{-}P_+ }{k_{-}^{\varphi}}\right)\mathrm{d}s}_{(2)}.
\end{aligned}
\end{equation}
First of all,  note that the quantity 
\begin{align*}
(-z+c\cdot\xi -ic_{3}L_0 -i \alpha\cdot\partial_{y}) \Big ( e^{h^{-1}\varrho_- (\tau - \varepsilon)} \mathfrak{M} \Big)= e^{h^{-1}\varrho_- (\tau - \varepsilon)} \big ( a + b\cdot \xi -i h^{-1}(\tau-\varepsilon)\alpha\cdot\partial_{y}\varrho_-\big) \mathfrak{M},
\end{align*}
with $\mathfrak{M}\in \mathscr{M}_{4}(\mathbb{C})$ and 
\begin{equation}\label{Tab}
\begin{aligned}
    a=-z+c_{3}\alpha\cdot\tilde{\nu}^{\varphi}\beta -i \alpha\cdot\partial_{y} \quad \text{ and }\quad  
    b=c+c_{3}\alpha\cdot\tilde{\nu}^{\varphi}\alpha
\end{aligned}
\end{equation}
belong to $\mathscr{M}_{4}(\mathbb{C})$. Note also the term $\alpha\cdot\partial_{y}$  in the quantity $a$ is applies to $\dfrac{\Pi_- P_+}{k_-^{\varphi}}$ in the following calculation. Now, we want to explain the quantities (1) and (2) given in \eqref{Expenontielle}. Let's start with $(1)$: 
\begin{equation}\label{(1)}
\begin{aligned}
(1)&=\int_{\varepsilon}^{\tau}e^{-h^{-1}s\varrho_{-}}\Pi_{-}a_{0}(y)\Big(-z+c\cdot\xi -ic_{3}L_{0}-i\alpha\cdot\partial_{y} \Big) \left(e^{h^{-1}(s-\varepsilon)\varrho_{-}}\dfrac{\Pi_{-}P_+}{k_{-}^{\varphi}}\right)\mathrm{d}s\\&= \int_{\varepsilon}^{\tau}e^{-\varepsilon h^{-1}\varrho_{-}}\Pi_{-}a_{0}(y)\Big(a + (b\cdot \xi) -i h^{-1}(\tau-\varepsilon)\alpha\cdot\partial_{y}\varrho_- \Big) \dfrac{\Pi_{-}P_+}{k_{-}^{\varphi}}\mathrm{d}s\\&= (\tau-\varepsilon)e^{-\varepsilon h^{-1}\varrho_{-}}\Pi_{-}a_{0}(y)\Big(a + b\cdot \xi \Big)B_{0,0} -i h^{-1}(\tau-\varepsilon)^{2}e^{-\varepsilon h^{-1}\varrho_{-}}\Pi_{-}a_{0}(y)\Big (\dfrac{\alpha\cdot\partial_{y}\varrho_-}{2}\Big) B_{0,0},
 \end{aligned}
 \end{equation}
 with $B_{0,0}(y,\xi)= \dfrac{\Pi_{-}P_+}{k_{-}^{\varphi}}\in \mathcal{S}^{0}$.\\
 Similarly, for (2) we get
\begin{equation}\label{(2)}
\begin{aligned}
(2)=&\int_{\varepsilon}^{\tau}e^{-h^{-1}s\varrho_{+}}\Pi_{+}a_{0}(y)\Big(-z+c\cdot\xi -ic_{3}L_{0}-i\alpha\cdot\partial_{y} \Big) \left(e^{h^{-1}(s-\varepsilon)\varrho_{-}}B_{0,0}\right)\mathrm{d}s\\&=
e^{-\varepsilon h^{-1}\varrho_{-}}\int_{\varepsilon}^{\tau}e^{h^{-1}s(\varrho_- -\varrho_{+})}\Pi_{+}a_{0}(y)\Big(a + b\cdot\xi -i h^{-1}(s-\varepsilon)\alpha\cdot\partial_{y}\varrho_-\Big)B_{0,0}\mathrm{ d}s\\&=e^{-\varepsilon h^{-1}\varrho_{-}} h(\varrho_- - \varrho_+)^{-1}\Pi_{+}a_{0}(y) \left(e^{h^{-1}(\varrho_- - \varrho_+)\tau} - e^{h^{-1}(\varrho_- - \varrho_+)\varepsilon}\right)\Big(a + b\cdot \xi \Big)B_{0,0}
\\&+e^{-\varepsilon h^{-1}\varrho_{-}}e^{h^{-1}(\varrho_- - \varrho_+)\tau}\Pi_{+}a_{0}(y) \left[ \dfrac{-i(\tau-\varepsilon)\alpha\cdot\partial_{y}\varrho_-}{\varrho_- - \varrho_+} + \dfrac{h \, i\alpha\cdot\partial_{y}\varrho_-}{(\varrho_- - \varrho_+)^{2}}\right]B_{0,0}\\& + e^{-\varepsilon h^{-1}\varrho_{-}}e^{h^{-1}(\varrho_- - \varrho_+)\varepsilon}\Pi_{+}a_{0}(y) \left[\dfrac{-h \,i\alpha\cdot\partial_{y}\varrho_-}{(\varrho_- - \varrho_+)^{2}}\right]B_{0,0}.
\end{aligned}
\end{equation}
Putting the formula of (1) and (2) as in \eqref{(1)} and \eqref{(2)}, respectively, in $I_2$. Together, with $I_1$, we obtain that
   \begin{align*}
	A_{1}(y, \xi, \tau)&=\left(e^{h^{-1}(\tau-\varepsilon)\varrho_{-}}\Pi_{-}+e^{h^{-1}(\tau-\varepsilon)\varrho_{+}}\Pi_{+}\right)A_{1}(y, \xi, \varepsilon)\\&+e^{h^{-1}\varrho_-(\tau-\varepsilon)}\Pi_{-}a_{0}(y)\left[(\tau-\varepsilon)\Big(a  + (b\cdot \xi )\Big)B_{0,0} -i h^{-1}(\tau-\varepsilon)^{2}\Big (\dfrac{\alpha\cdot\partial
 _{y}\varrho_-}{2}\Big) B_{0,0}\right]\\& + \dfrac{h}{(\varrho_- - \varrho_+)}\Pi_{+}a_{0}(y) e^{h^{-1}\varrho_-(\tau-\varepsilon)} \Big(a + (b\cdot \xi) \Big)B_{0,0} \\&+e^{h^{-1}\varrho_-(\tau-\varepsilon)}\Pi_{+}a_{0}(y) \left[ \dfrac{-i(\tau-\varepsilon)\alpha\cdot\partial_{y}\varrho_-}{\varrho_- - \varrho_+} + \dfrac{h\,i \alpha\cdot\partial_{y}\varrho_-}{(\varrho_- - \varrho_+)^{2}}\right]B_{0,0}\\& + e^{h^{-1}\varrho_+  (\tau - \varepsilon)}\Pi_{+}a_{0}(y) \left[-i\dfrac{h\,\, \alpha\cdot\partial_{y}\varrho_-}{(\varrho_- - \varrho_+)^{2}}\right]B_{0,0} - e^{h^{-1}\varrho_+(\tau-\varepsilon)}\dfrac{h\,\,\Pi_{+}a_{0}(y)}{(\varrho_- - \varrho_+)} \Big(a B_{0,0} + b\cdot\xi B_{0,0}  \Big) .
	\end{align*}
Thanks to the properties of $\varrho_+$ given in \eqref{ellip}, and the fact that $e^{h^{-1}(\tau-\varepsilon)\varrho_{+}}\Pi_+ a_0(y)$ is unbounded in $L^{2}(\lbrace \tau > \varepsilon\rbrace)$, then we look $ A_{1}(\tilde{y}, \xi, \varepsilon)$ such that
	\begin{align}\label{condition}
	   \Pi_{+}A_{1}(y, \xi, \varepsilon)= h\,\,\dfrac{\Pi_{+}a_{0}}{\varrho_{-} - \varrho_{+}}\Big(a  + b\cdot\xi  +\dfrac{i\alpha\cdot\partial_{y}\varrho_-}{\varrho_{-} - \varrho_{+}}  \Big)B_{0,0}.
	\end{align}
	Thus, we obtain 
\begin{equation}\label{PreA1}
\begin{aligned}
&A_{1}(y, \xi, \tau)=e^{h^{-1}(\tau-\varepsilon)\varrho_{-}}\times\Bigg\lbrace \Pi_{-}A_{1}(y, \xi, \varepsilon)+ h\,\,\Pi_{+}a_{0}(y){(\varrho_- - \varrho_+)}  \left[a  + b\cdot\xi + \dfrac{i\alpha\cdot\partial_{y}\varrho_-}{(\varrho_- - \varrho_+)} \right]B_{0,0}\\&+ (\tau-\varepsilon)\left[\Pi_{-}a_{0}(y)\big(a  + b\cdot\xi \big) -\Pi_{+}a_{0}(y) \dfrac{i\alpha\cdot\partial_{y}\varrho_-}{(\varrho_- - \varrho_+)}\right]B_{0,0} + h^{-1}(\tau-\varepsilon)^{2}\Pi_{-}a_{0}(y)\Big( \frac{-i\alpha\cdot\partial_{y}\varrho_{-}}{2} \Big) B_{0,0}\Bigg\rbrace.
\end{aligned}
\end{equation}
Calculate of $\Pi_{-} A_{1}(y,\xi,\varepsilon).$
From \eqref{PreA1}, we get that
\begin{align*}
A_{1}(y,\xi,\varepsilon)= \Pi_-(P_- + P_+)A_{1}(y,\xi,\varepsilon) +  \dfrac{h\,\,\Pi_{+}a_{0}(y)}{(\varrho_- - \varrho_+)}  \left[a  + b\cdot\xi  + \dfrac{i\alpha\cdot\partial_{y}\varrho_-}{(\varrho_- - \varrho_+)} \right]B_{0,0} .
\end{align*}
From \eqref{T1approx} we have $P_{+} A_{1}(y,\xi,\varepsilon)=0$, then
\begin{align*}
P_{-}A_{1}(y,\xi,\varepsilon)= P_{-}\Pi_-P_{-}A_{1}(y,\xi,\varepsilon) +   \dfrac{h\,\,\Pi_{+}a_{0}(y)}{(\varrho_- - \varrho_+)}  \left[a  + b\cdot\xi  + \dfrac{i\alpha\cdot\partial_{y}\varrho_-}{(\varrho_- - \varrho_+)} \right]B_{0,0}. 
\end{align*}
Thanks to the relations \eqref{PPiP}, we obtain 
\begin{align*}
\Pi_{-}P_{-}A_{1}(y,\xi,\varepsilon)=  \dfrac{h\,\Pi_{-}a_{0}P_+}{(\varrho_- - \varrho_+)} \left(I_{4} - \dfrac{ \Theta^{\varphi}}{k^{\varphi_{-}}}\right) \left[a  + b\cdot\xi + \dfrac{i\alpha\cdot\partial_{y}\varrho_-}{(\varrho_- - \varrho_+)} \right]B_{0,0}, 
\end{align*}
and so \eqref{PreA1} becomes as follows
\begin{align*}\label{FA1}
&A_{1}(y, \xi, \tau)=e^{h^{-1}(\tau-\varepsilon)\varrho_{-}}\times\Bigg\lbrace h \left[\Pi_{-}a_{0}\left(P_+ - \dfrac{ P_+\,\Theta^{\varphi}}{k^{\varphi_{-}}}\right) +\Pi_{+}a_0\right]\left[\dfrac{a  + b\cdot\xi}{(\varrho_- - \varrho_+)} + \dfrac{i\alpha\cdot\partial_{y}\varrho_-}{(\varrho_- - \varrho_+)^{2}} \right]B_{0,0}\\&+ (\tau-\varepsilon)\left[\Pi_{-}a_{0}(y)\big (a  + b\cdot\xi \big) - \Pi_{+}a_{0}(y) \dfrac{i\alpha\cdot\partial_{y}\varrho_-}{(\varrho_- - \varrho_+)}\right]B_{0,0} \, + \, h^{-1}(\tau-\varepsilon)^{2}\Pi_{-}a_{0}(y)\Big( \frac{-i\alpha\cdot\partial_{y}\varrho_{-}}{2} \Big) B_{0,0}\Bigg\rbrace.
\end{align*}
Consequently, we get that
\begin{equation}\label{A_1}
\begin{aligned}
A_{1}(y, \xi, \tau)&=e^{h^{-1}\varrho_{-}(\tau-\varepsilon)}\Big\lbrace B_{1, 0}(y,\xi)+\big(h^{-1}(\tau - \varepsilon)(\varrho_- - \varrho_+)\big)B_{1,1}(y,\xi) \\& \hspace{6cm}+\big(h^{-1}(\tau - \varepsilon)(\varrho_- - \varrho_+)\big)^{2} B_{1,2}(y,\xi)\Big\rbrace \\& = e^{h^{-1}(\tau-\varepsilon)\varrho_{-}} \sum_{k=0}^{2} \Big(h^{-1}(\tau - \varepsilon)(\varrho_- - \varrho_+)\Big)^{k} B_{1,k} (y,\xi) ,
\end{aligned}
\end{equation}
where,
\begin{equation*}
\begin{aligned}
B_{1, 0}(y,\xi)&= h \left[\Pi_{-}a_{0}\left(P_+ - \dfrac{ P_+\,\Theta^{\varphi}}{k^{\varphi}_{-}}\right) +\Pi_{+}a_0\right]\left[\dfrac{a  + b\cdot\xi}{(\varrho_- - \varrho_+)} + \dfrac{i\alpha\cdot\partial_{y}\varrho_-}{(\varrho_- - \varrho_+)^{2}} \right]B_{0, 0},\\
B_{1, 1}(y,\xi)&= h \left[\Pi_{-}a_{0}(y)\dfrac{\big(a  + b\cdot\xi\big)}{(\varrho_- - \varrho_+)}  - \Pi_{+}a_{0}(y) \dfrac{i\alpha\cdot\partial_{y}\varrho_-}{(\varrho_- - \varrho_+)^{2}}\right]B_{0, 0},\\
B_{1, 2}(y,\xi)&=-h\,\,\Pi_{-}a_{0}(y)\Big( \frac{i\alpha\cdot\partial_{y}\varrho_{-}}{2(\varrho_- - \varrho_+)^{2}} \Big)B_{0, 0} ,
\end{aligned}
\end{equation*}
with $\varrho_{-}-\varrho_{+}=-2\lambda(y,\xi)\in \mathcal{S}^{1}$ and $\partial_{y}\varrho_-\in \mathcal{S}^{1}$ identified with $\langle\xi\rangle$, then  $B_{1,k}\in h\,\mathcal{S}^{0}$ for $k=0,1$, and $B_{1,2}\in h\,\mathcal{S}^{-1}.$\\\\

Let's look at the form of $A_j$ for $j=2$.  To do it, we define $A_{2}(y, \xi, \tau)$ inductively by 
\begin{equation*}
	\left\{
\begin{aligned}\label{T2approx}	
& h\partial_{\tau} A_{2}(y, \xi, \tau)  =  L_{0}(y, \xi)A_{2}(y, \xi, \tau)+\Big(\widetilde{L}_1(y,\xi) -i  \partial_{\xi}L_{0}(y,\xi)\cdot\partial_{y}\Big)A_{1}(y, \xi, \tau)   \\& \hspace{6cm}+\Big((\alpha\cdot\tilde{\nu}^{\varphi}c_{3})\widetilde{L}_{1}(y,\xi) -i\partial_{\xi}\widetilde{L}_1(y,\xi)\cdot\partial_{y}\Big)A_{0}(y, \xi, \tau),\\
	& P_+ A_{2}(y, \xi, \varepsilon)  =  0,  
	\end{aligned}
	\right.
	\end{equation*}
where,
\begin{equation}\label{TL2}
\begin{aligned}
&\big(\widetilde{L}_{1} - i\partial_{\xi}L_{0}\cdot\partial_{y}\big)\Big ( e^{h^{-1}\varrho_- (\tau - \varepsilon)} \mathfrak{M} \Big)=e^{h^{-1}\varrho_- (\tau - \varepsilon)} a_{0}(y)\big(a + b\cdot \xi - ih^{-1}(\tau-\varepsilon)\alpha\cdot\partial_{y}\varrho_-\big) \mathfrak{M},\\ 
&\big(-i\partial_{\xi}\widetilde{L}_1\cdot\partial_{y}+(\alpha\cdot\tilde{\nu}^{\varphi}c_{3})\widetilde{L}_{1}\big)\Big ( e^{h^{-1}\varrho_- (\tau - \varepsilon)} \mathfrak{M} \Big) =  e^{h^{-1}\varrho_- (\tau - \varepsilon)} a_{0}(y)\Big ( d + e\cdot \xi -i h^{-1}(\tau-\varepsilon) f \cdot\partial_{y}\varrho_- \Big) \mathfrak{M} ,
\end{aligned}
\end{equation}
 with $\mathfrak{M}$, $a,\,b $ were noted in \eqref{Tab}, $d,\,e,\,f $ belong to $\mathscr{M}_{4}(\mathbb{C})$, where $d$, $e$ and $f$ are the following
 \begin{equation}\label{Tdf}
     \begin{aligned}
         &d=  (c_{3}\alpha\cdot\tilde{\nu}^{\varphi})^{2}\beta - c_{3}\alpha\cdot\tilde{\nu}^\varphi z - i(c+  c_{3}\alpha\cdot\tilde{\nu}^{\varphi}\alpha)\cdot\partial_{y},\\ &e = \big ( c_{3}\alpha\cdot\tilde{\nu}^{\varphi} + (c_{3}\alpha\cdot\tilde{\nu}^{\varphi})^{2}\alpha \big)\cdot\xi \quad \text{ and } \quad 
         f =  c+c_{3}\alpha\cdot\tilde{\nu}^{\varphi} \alpha.
     \end{aligned}
 \end{equation}
 Then, after a many calculation, we arrive at the following formula 
 \begin{equation}\label{A2}
\begin{aligned}
A_{2}(y, \xi, \tau)&=e^{h^{-1}\varrho_{-}(\tau-\varepsilon)}\Big\lbrace \big(h^{-1}(\tau - \varepsilon)\langle\xi\rangle\big)^{0}B_{2, 0}(y,\xi)+\big(h^{-1}(\tau - \varepsilon)\langle\xi\rangle\big)^{1}B_{2,1}(y,\xi)\\&+\big(h^{-1}(\tau - \varepsilon)\langle\xi\rangle\big)^{2}B_{2,2}(y,\xi) +\big(h^{-1}(\tau - \varepsilon)\langle\xi\rangle\big)^{3}B_{2,3}(y,\xi) + \big(h^{-1}(\tau - \varepsilon)\langle\xi\rangle\big)^{4}B_{2,4}(y,\xi)\Big\rbrace \\& =: e^{h^{-1}\varrho_{-}(\tau-\varepsilon)} \sum_{k=0}^{4} \Big(h^{-1}(\tau - \varepsilon)\langle\xi\rangle\Big)^{k} B_{1,k} (y,\xi),
\end{aligned}
\end{equation}
where,
\begin{equation*}
\begin{aligned}
B_{2, 0}(y,\xi)&=  h \,\,\Bigg[\Pi_{-}a_{0}\left(P_+ - \dfrac{ P_+\,\Theta^{\varphi}}{k^{\varphi}_{-}}\right)+\Pi_{+}a_{0}\Bigg]\Bigg[\dfrac{a B_{1,0} + (b\cdot \xi) B_{1,0} + d B_{0,0} + (e\cdot \xi)  B_{0,0} }{(\varrho_- - \varrho_+)}\\&\hspace{3cm}-\dfrac{(\alpha\cdot\partial_{y}\varrho_-) B_{1,0}+ \langle\xi\rangle a B_{1,1} + \langle\xi\rangle(b\cdot \xi) B_{1,1}+ f\cdot\partial_{y}\varrho_- B_{0,0}}{(\varrho_- - \varrho_+)^{2}}  \\& +\dfrac{ 2 \langle\xi\rangle \alpha\cdot\partial_y \varrho_- B_{1,1} + 2 \langle\xi\rangle (b\cdot \xi) B_{1,0} + 2 \langle\xi\rangle^{2} (b\cdot \xi) B_{1,2}      }{(\varrho_- - \varrho_+)^{3}} - \dfrac{6  \langle\xi\rangle^{2} \alpha\cdot\partial_y \varrho_- B_{1,2}}{(\varrho_- - \varrho_+)^{4}}\Bigg],
\end{aligned}
\end{equation*}
\begin{align*}
B_{2, 1}(y,\xi)&= h \,\,\Pi_{-}a_{0}(y)\left[\dfrac{a B_{1,0} + d B_{0,0}}{(\varrho_- - \varrho_+)} + b B_{1,0} + e B_{0,0}\right] \\& \hspace{2cm}+ h\,\,\Pi_{+}a_{0}(y) \Bigg[\dfrac{f\cdot\partial_{y}\varrho_-B_{0,0} + a B_{1,1} +\langle\xi\rangle B_{1,1}}{(\varrho_- - \varrho_+)} + \dfrac{\alpha\cdot\partial_{y}\varrho_- B_{1,0}}{(\varrho_- - \varrho_+) \langle\xi\rangle}\\& \hspace{3cm} - \dfrac{2\alpha\cdot\partial_y \varrho_- B_{1,1} + (2\langle\xi\rangle a + 2\langle\xi\rangle^{2}) B_{1,2}} {(\varrho_- - \varrho_+)^{2}}+ \dfrac{6\langle\xi\rangle \alpha\cdot\partial_y \varrho_- B_{1,2}}{(\varrho_- - \varrho_+)^{3}}\Bigg],
\end{align*}
\begin{align*}
B_{2, 2}(y,\xi)&= h \,\,\Pi_{-}a_{0}(y)\left[\dfrac{a B_{1,1}}{2(\varrho_- - \varrho_+)}+ \dfrac{b B_{1,1}}{2} + \dfrac{(\alpha\cdot\partial_{y}\varrho_-)B_{1,0}}{2(\varrho_- - \varrho_+)^{2}}+ \dfrac{(f\cdot\partial_{y}\varrho_-)B_{0,0}}{2(\varrho_- - \varrho_+)^{2}}\right] \\&+ h \,\,\Pi_{+}a_{0}(y) \Bigg[\dfrac{(\alpha\cdot\partial_{y}\varrho_-)B_{1,1}}{(\varrho_- - \varrho_+)\langle\xi\rangle} - \dfrac{a B_{1,2}}{(\varrho_- - \varrho_+)} + \dfrac{(b\cdot \xi) B_{1,2}}{(\varrho_-  -  \varrho_+)} - \dfrac{3(\alpha\cdot\partial_{y}\varrho_-)B_{1,2}}{(\varrho_- - \varrho_+)^{2}}\Bigg],
\end{align*}
\begin{align*}
B_{2, 3}(y,\xi)&=h\,\,\Pi_{-}a_{0}(y)\left[\dfrac{a B_{1,2}}{3\langle\xi\rangle}+ \dfrac{b B_{1,2}}{3} +\dfrac{(\alpha\cdot\partial_{y}\varrho_-)B_{1,1}}{3 (\varrho_- - \varrho_+)^{2}}\right] + h \,\,\Pi_{+}a_{0}(y) \Bigg[\dfrac{(\alpha\cdot\partial_{y}\varrho_-)\big(B_{1,2}\big)}{(\varrho_- - \varrho_+)\langle\xi\rangle}\Bigg],
\\
B_{2, 4}(y,\xi)& =h\,\,\Pi_{-}a_{0}(y)\Bigg[ \frac{(\alpha\cdot\partial_{y}\varrho_{-}) \, B_{1,2}}{4(\varrho_- - \varrho_+)^{2}} \Bigg],
\end{align*}
with $\varrho_{-}-\varrho_{+}=-2\lambda(y,\xi)\in \mathcal{S}^{1}$ and $\partial_{y}\varrho_-\in \mathcal{S}^{1}$. Then  $B_{2,k}\in h\,\mathcal{S}^{0}$ for $k=0,1,2$, $B_{2,3}\in h^{2}\,\mathcal{S}^{-1}$, and $B_{2,4}\in h^{2}\,\mathcal{S}^{-2}$.\qed
\begin{remark}
Using \eqref{deflambdapro} and \eqref{PPiP}, then the boundary condition associated with $A_{2}(y,\xi,\varepsilon)$ is the following
\begin{align*}
&\Pi_{+}A_{2}(y,\xi,\varepsilon)=  h \,\Pi_{+}a_{0} \Bigg[\dfrac{a B_{1,0} + (\xi\cdot b) B_{1,0} + d B_{0,0} + (e\cdot \xi)  B_{0,0} }{(\varrho_- - \varrho_+)}\\&\hspace{3cm}-\dfrac{(\alpha\cdot\partial_{y}\varrho_-) B_{1,0}+ \langle\xi\rangle a B_{1,1} + \langle\xi\rangle (b\cdot \xi) B_{1,1}+ f\cdot\partial_{y}\varrho_- B_{0,0}}{(\varrho_- - \varrho_+)^{2}}  \\& \hspace{4cm}+\dfrac{ 2 \langle\xi\rangle \alpha\cdot\partial_y \varrho_- B_{1,1} + 2 \langle\xi\rangle (b\cdot \xi) B_{1,0} + 2 \langle\xi\rangle^{2} (b\cdot \xi) B_{1,2}      }{(\varrho_- - \varrho_+)^{3}} - \dfrac{6  \langle\xi\rangle^{2} \alpha\cdot\partial_y \varrho_- B_{1,2}}{(\varrho_- - \varrho_+)^{4}}\Bigg].
\end{align*}
\end{remark}

\end{document}